\documentclass[11pt]{article}
\usepackage{xcolor}
\usepackage{amsthm}
\usepackage{xcolor}
\usepackage{textcomp}
\usepackage{enumerate}      
\usepackage{graphicx}        
\usepackage{caption, subcaption}

\usepackage{url}

\usepackage{tabu}

\usepackage[export]{adjustbox}

\usepackage{mathrsfs}

\usepackage{url}

\usepackage{array}
\newcommand{\PreserveBackslash}[1]{\let\temp=\\#1\let\\=\temp}
\newcolumntype{C}[1]{>{\PreserveBackslash\centering}p{#1}}
\newcolumntype{R}[1]{>{\PreserveBackslash\raggedleft}p{#1}}
\newcolumntype{L}[1]{>{\PreserveBackslash\raggedright}p{#1}}

\renewcommand{\setminus}{{\smallsetminus}}

\newcommand{\cp}[1]{\vcenter{\hbox{#1}}}

\usepackage{amssymb,amsmath,bm}
\usepackage{amsthm}
\usepackage{hyperref}
\usepackage{mathrsfs}
\usepackage{textcomp}
\usepackage{enumerate}
\usepackage{graphicx}
\usepackage{tikz, tikz-cd}
\usepackage{verbatim}

\newtheorem{theorem}{Theorem}[section]
\newtheorem{lemma}[theorem]{Lemma}
\newtheorem{proposition}[theorem]{Proposition}
\newtheorem{definition}[theorem]{Definition}

\newtheorem{conjecture}[theorem]{Conjecture}

\theoremstyle{remark}
\newtheorem{remark}[theorem]{Remark}
\newcommand{\R}{\mathbb{R}}

\theoremstyle{remark}

\numberwithin{equation}{section}

\title{\bf On the Volume Conjectures for Dehn-filled $3$-manifolds along fundamental shadow link complements with big volume}
\author{Yunpeng Meng, Baojun Wu, Tian Yang, and Yuxuan Yang}

\date{}

\begin{document}
\maketitle
\begin{abstract}
 We prove the Volume Conjecture and the Asymptotic Expansion Conjecture for closed oriented hyperbolic $3$-manifolds obtained by integral Dehn-fillings along a fundamental shadow link complement, provided that their volumes are sufficiently close to that of the complement. More generally, we show that the conjectured asymptotic expansions hold for pairs obtained by change-of-pair operations on a fundamental shadow link complement whenever the volume of the corresponding hyperbolic cone manifold is sufficiently close to that of the complement. A key ingredient in the proof is the introduction of extended quantum $6j$-symbols and a Weak Reflection Symmetry satisfied by them. The Weak Reflection Symmetry is established by relating the extended quantum $6j$-symbols to the $b$-$6j$ symbols through  a real Gevrey-$1$ estimate and a truncated Borel-Laplace transform.
\end{abstract}

\tableofcontents

\section{Introduction}

A central theme in quantum topology is the principle that the asymptotic behavior of quantum invariants reflects the geometry of $3$-manifolds. This idea first emerged in Kashaev's volume conjecture, which relates the exponential growth of a knot invariant constructed from the quantum dilogarithm to the hyperbolic volume of the knot complement \cite{Kashaev}. Murakami and Murakami subsequently showed that Kashaev's invariant can be identified with a particular specialization of the colored Jones polynomial \cite{MurakamiMurakami}. This perspective was later extended to hyperbolic $3$-manifolds through analogous volume conjectures for Reshetikhin-Turaev and relative Reshetikhin-Turaev invariants \cite{CY,O,GRY,WY1,WY2}.
\\

Let $M$ be a closed oriented hyperbolic $3$-manifold. For any integer $r\geqslant 3,$ let $\mathrm{RT}_r(M)$ be the $r$-th Reshetikhin-Turaev invariant of $M$ evaluated at the root of unity $q=e^{\frac{2\pi\mathbf i}{r}}.$ 

\begin{conjecture}\cite[Conjecture 1.2]{CY}\label{VC} As $r$ varies over all positive odd integers, $$\lim_{r\to+\infty}\frac{4\pi }{r}\log\big( \mathrm{RT}_{r}(M)\big)= \mathrm{Vol}(M)+\mathbf i\mathrm{CS}(M) \mod \mathbf i\pi^2\mathbb Z,$$
where $\mathrm{Vol}(M)$ and $\mathrm{CS}(M)$ are respectively the volume and the Chern-Simons invariant of $M$ in the hyperbolic metric.
\end{conjecture}

\begin{conjecture}\cite{O,GRY}\label{Asymp} As $r$ varies over all positive odd integers,
$$\mathrm{RT}_r(M)= \frac{C_r}{\sqrt{\pm\mathrm{Tor}(M)}}e^{\frac{r}{4\pi}\big(\mathrm{Vol}(M)+\mathbf i\mathrm{CS}(M)\big)}\Big(1+O\Big(\frac{1}{r}\Big)\Big),$$
where $C_r$ is a complex number of norm $1$ independent of the geometric structure on $M,$ and $\mathrm{Tor}(M)$ is the adjoint twisted Reidemeister torsion of $M.$
\end{conjecture}

 Conjectures \ref{VC} and \ref{Asymp} are respectively proved in \cite{O, WY, CZ, GMWY} for hyperbolic $3$-manifolds obtained by Dehn-fillings along the figure-$8$ knot, and by long Dehn-fillings along the family of twist knots. 
\medskip

Let $M$ be a closed oriented $3$-manifold and let $L$ be a framed hyperbolic link in $M$ with $n$ components. For an odd integer $r\geqslant 3,$ let  $\boldsymbol  b=(b_1,\dots,b_n)\in\{0,2,\dots,r-3\}^n$ and let $\mathrm{RT}_r(M,L,\boldsymbol  b)$ be the $r$-th relative Reshetikhin-Turaev invariant of the pair $(M,L)$ colored by $\boldsymbol  b$ and evaluated at the root of unity $q=e^{\frac{2\pi\mathbf i}{r}}.$ For a sequence $\boldsymbol  b^{(r)}=(b^{(r)}_1,\dots,b^{(r)}_n),$ let 
$$\theta_k=\bigg|\lim_{r\to\infty}\frac{4\pi b^{(r)}_k}{r}-2\pi\bigg|$$
for $k\in\{1,\dots,n\},$ and let $\boldsymbol\theta=(\theta_1,\dots, \theta_n).$

\begin{conjecture}\cite[Conjecture 1.1]{WY1}\label{CVC}  
If $M_{L_{\boldsymbol\theta}}$ is a hyperbolic cone manifold consisting of $M$ and a hyperbolic cone metric on $M$ with singular locus $L$ and cone angles $\boldsymbol\theta,$ then as $r$ varies over all positive odd integers,
$$\lim_{r\to\infty} \frac{4\pi }{r}\log \mathrm{RT}_r(M,L,\boldsymbol  b^{(r)})=\mathrm{Vol}(M_{L_{\boldsymbol\theta}})+\mathbf i\mathrm{CS}(M_{L_{\boldsymbol\theta}}) \mod \mathbf i\pi^2\mathbb Z,$$
where $\mathrm{Vol}(M_{L_{\boldsymbol\theta}})$ and $\mathrm{CS}(M_{L_{\boldsymbol\theta}})$ are respectively the volume and the Chern-Simons invariant of $M_{L_{\boldsymbol \theta}}.$
\end{conjecture}
Suppose further  that for each $k\in\{1,\dots,n\},$ either $b_k^{(r)}> \frac{r}{2}$ for  all $r$ sufficiently large or $b_k^{(r)}< \frac{r}{2}$ for all $r$ sufficiently large. In the former case we let $\mu_k=1$ and in the latter case we let $\mu_k=-1,$ and we let 
$$\theta^{(r)}_k=\mu_k\bigg(\frac{4\pi b^{(r)}_k}{r}-2\pi\bigg).$$
Let $\boldsymbol \theta^{(r)}=(\theta^{(r)}_1,\dots,\theta^{(r)}_{n}),$ and suppose that $\{\boldsymbol \theta^{(r)}\}$ converges as $r$ tends to infinity. Also suppose that for all $r$ sufficiently large, a hyperbolic cone metric on $M$ with singular locus $L$ and cone angles $\boldsymbol \theta^{(r)}$ exists. We denote $M$ with such a hyperbolic cone metric by $M^{(r)},$ and let $\mathrm{Vol}(M^{(r)})$ and $\mathrm{CS}(M^{(r)})$ respectively be the volume and the Chern-Simons invariant of $M^{(r)}.$ For the parallel copies $\{\gamma_k\}_{k\in\{1,\dots,n\}} $ of the core curves of $L$ given by the framing, let $\mathrm H^{(r)}(\gamma_k)$ be the logarithmic holonomy in $M^{(r)}.$ 
Let $\rho_{M^{(r)}}:\pi_1(M\setminus L)\to \mathrm{PSL}(2;\mathbb C)$ be the holonomy representation of  the restriction of $M^{(r)}$ to $M\setminus L,$ and let $\mathbb{T}_{(M\setminus L,\boldsymbol m)}([\rho_{M^{(r)}}])$ be the Reidemeister torsion of $M\setminus L$  twisted by the adjoint action of  $\rho_{M^{(r)}}$ with respect to the system of meridians $\boldsymbol m$ of a tubular neighborhood of the core curves of $L.$

\begin{conjecture}\cite[Conjecture 1.1]{WY2}\label{CAsymp} 
As $r$ varies over all positive odd integers,
$$ \mathrm{RT}_r(M,L,\boldsymbol b^{(r)})=C_r\frac{e^{\sum_{k=1}^{n}\mu_k\frac{\mathrm H^{(r)}(\gamma_k)}{2}}}{\sqrt{\pm\mathbb{T}_{(M\setminus L,\boldsymbol m)}([\rho_{M^{(r)}}])}}e^{\frac{r}{4\pi}\big(\mathrm{Vol}(M^{(r)})+\mathbf i\mathrm{CS}(M^{(r)})\big)}\bigg(1+O\Big(\frac{1}{r}\Big)\bigg),$$
where $C_r$ is a complex number of norm $1$ independent of the geometric structure on $(M,L).$
\end{conjecture}

\begin{remark}\label{rmk} 
If we let 
$$\widetilde \theta^{(r)}_k=\mu_k\bigg(\frac{4\pi (b^{(r)}_k+1)}{r}-2\pi\bigg)$$
for each $k\in\{1,\dots,n\},$ let $\widetilde {\boldsymbol\theta}^{(r)}=(\widetilde \theta_1^{(r)},\dots,\widetilde \theta_n^{(r)}),$ and by abuse of notation let $M^{(r)}=M_{L_{\widetilde {\boldsymbol\theta}^{(r)}}}$ be the corresponding hyperbolic cone manifold, 
then Conjecture \ref{CAsymp} can be rephrased as: 
As $r$ varies over all positive odd integers,
$$ \mathrm{RT}_r(M,L,\boldsymbol b^{(r)})=\frac{C_r}{\sqrt{\pm\mathbb{T}_{(M\setminus L,\boldsymbol m)}([\rho_{M^{(r)}}])}}e^{\frac{r}{4\pi}\big(\mathrm{Vol}(M^{(r)})+\mathbf i\mathrm{CS}(M^{(r)})\big)}\bigg(1+O\Big(\frac{1}{r}\Big)\bigg).$$
The reason is due to the Complexified Schl\"afli Formula of cone metrics (see, \cite[Theorem 3.20]{CHK} for the real part and   \cite[Theorem 3.5]{HLM} for the imaginary part) that, if $\widetilde  {\boldsymbol \theta}^{(r)}= \boldsymbol \theta^{(r)} + \frac{4\pi }{r} {\boldsymbol \mu},$ where $\boldsymbol \mu = (\mu_k)_{k\in\{1,\dots,n\}} \in \{-1,1\}^n,$ then 
$$\mathrm{Vol}(M_{L_{\widetilde{\boldsymbol \theta}^{(r)}}})+\mathbf i \mathrm{CS}(M_{L_{\widetilde{\boldsymbol \theta}^{(r)}}}) = \mathrm{Vol}(M_{L_{\boldsymbol \theta^{(r)}}})+\mathbf i \mathrm{CS}(M_{L_{\boldsymbol \theta^{(r)}}}) -\frac{1}{2} \sum_{k=1}^n{\mu_k}\mathrm{H}^{(r)}(\gamma_k)\cdot\frac{4\pi}{r} + O\Big(\frac{1}{r^2}\Big).$$
\end{remark}

Conjectures \ref{CVC} and \ref{CAsymp} were proved in \cite{WY1, WY2, PW} for the pairs $(M,L)$ such that  $M\setminus L$ is homeomorphic to a fundamental shadow link complement (see Section \ref{fundsl}) and all the cone angles are sufficiently small. They were also proved  in \cite{PW} for pairs $(M,L)$ such that  $M$ is obtained by doing a Dehn-filling along a fundamental shadow link complement with sufficiently large surgery coefficients, $L$ consists of the core of the filled solid tori, and all the cone  angles are strictly less than $\pi.$ Topologically, $(M,L)$ is obtained by doing a change-of-pair operation from the fundamental shadow link complement (see Section \ref{COP}.)
\\

As shown in \cite{CT}, every oriented 3-manifold with toroidal or empty boundary can be obtained from a suitable fundamental shadow link complement by performing  integral Dehn-fillings on some of its boundary
components. Notice also that, on the one hand, the relative Reshetikhin-Turaev invariants of the triple $(M,L,\boldsymbol 0)$ coincide with the Reshetikhin-Turaev invariants of $M,$ while, on the other hand, the coloring $\boldsymbol 0$ corresponds to the cone metric on $(M,L)$ for which all the cone angles are equal to $2\pi;$ this is precisely the hyperbolic metric on $M.$ 

Therefore, as proposed in \cite{WY1}, to prove Conjectures \ref{VC} and \ref{Asymp} for all the closed oriented hyperbolic $3$-manifolds, it suffices to prove Conjectures \ref{CVC} and \ref{CAsymp} for pairs $(M,L)$ such that  $M$ is obtained by doing an integral Dehn-filling along a fundamental shadow link complement, $L$ consists of the cores of the filled solid tori, and the coloring is $\boldsymbol 0.$ In view of  the result of \cite{PW} mentioned above, it therefore suffices to 
\begin{enumerate}[(1)]
\item extend the range of the cone angles from $[0,\pi)$ to $[0,2\pi],$ and 
\item remove the assumption that the Dehn-filling  coefficients are sufficiently large. 
\end{enumerate}

The main results of this article accomplish step (1). More precisely, we have

\begin{theorem}\label{main1}
For any fundamental shadow link complement $M_\text{FSL},$ there exists a $\delta >0$ such that if the pair $(M,L)$ obtained from $M_\text{FSL}$ by doing a change-of-pair operation and the cone angles $\boldsymbol \theta$ satisfy 
$$\mathrm{Vol}(M_{L_{\boldsymbol \theta}})> \mathrm{Vol}(M_\text{FSL})-\delta,$$ 
then Conjecture \ref{CVC} and Conjecture \ref{CAsymp} hold for $(M,L).$  
\end{theorem}

The change-of-pair operation will be recalled in Section \ref{COP}. As part of the construction, the closed $3$-manifold $M$ in Theorem \ref{main1} is obtained by doing an integral Dehn-filling along a fundamental shadow link complement $M_\text{FSL},$ and as a consequence, we have the following

\begin{theorem}\label{main2}
For any fundamental shadow link complement $M_\text{FSL},$ there exists a $\delta >0$ such that if the closed oriented hyperbolic $3$-manifold $M$ obtained by doing an integral Dehn-filling along $M_\text{FSL}$ satisfies 
$$\mathrm{Vol}(M)> \mathrm{Vol}(M_\text{FSL})-\delta,$$ then Conjecture \ref{VC} and Conjecture \ref{Asymp} hold for $M.$
\end{theorem}

\begin{remark}
With techniques in \cite{WY, PW, GMWY}, one can generalize Theorem \ref{main2} to any closed or  cusped hyperbolic $3$-manifold $M$ obtained by doing a rational Dehn-filling along some components of $M_\text{FSL}$ that satisfies $\mathrm{Vol}(M)> \mathrm{Vol}(M_\text{FSL})-\delta.$
\end{remark}

The main technical difficulty we overcome in this paper is the following. In a direct application of Ohtsuki's method that writes the  summation in the definition of the Reshetikhin-Turaev invariant into the summation of Fourier coefficients (each of which is a contour integral) then estimates each of these contour integrals using the saddle point approximation, there are certain problematic Fourier coefficients  that cannot be estimated by the usual method. That is essentially why the assumption that the cone angles are less than $\pi$ is required in \cite{PW}. 

To resolve this issue, we introduce  the \emph{extended quantum $6j$-symbols} (see Definition \ref{eq6j}) and establish a \emph{Weak Reflection Symmetry} for them (see Proposition \ref{key} and Theorem \ref{wrs}). We then carefully add suitable extra terms arising from the 
extended $6j$-symbols into the summation defining the Reshetikhin-Turaev invariant (see Formula (\ref{rt3})). These additional terms, on the one hand, do not change the asymptotic behavior of the invariant (see Propositions \ref{rt4} and \ref{rt5}) and, on the other hand, annihilate those problematic Fourier coefficients (see Propositions \ref{Poisson} and \ref{Poisson2}). 

To prove the Weak Reflection Symmetry of the extended quantum $6j$-symbols, we relate them to the $b$-$6j$ symbols, objects that have recently been intensively studied by the second and third authors and their collaborators in \cite{LMSWY,LMSWY2}  and  satisfy a reflection symmetry. This relation is established using a  \emph{real Gevrey-$1$ estimate} and a truncated Borel-Laplace transform developed in this paper (see Section \ref{pfKey} and Appendix \ref{Gevrey-1}). 
\\

\noindent{\bf Acknowledgment.} Y.M. thanks Shanzhong Sun for discussions, and  B.W. and T.Y. thank Xin Sun for discussions. 
B.W. is supported by National Key R\&D Program of China (No. 2023YFA 1010700) and NSF China Grant (No. 12526204). T.Y. is supported, and Y.M. is partially supported, by NSF Grant DMS-2505908. 
\\

\noindent{\bf AI Usage Disclosure.} The authors use AI assistance for the proof of Proposition \ref{EST} (2), Proposition \ref{shGevrey} and the proofs in Appendix \ref{Gevrey-1}.


\section{Preliminaries}

Most of the material in this section, except Proposition \ref{EST} (2), is taken from \cite{WY1,LMSWY,LMSWY2} and is included here for the reader's convenience. 

\subsection{Relative Reshetikhin-Turaev invariants}\label{RRT}

We recall the relative Reshetikhin-Turaev invariants\,\cite{BHMV, Li} following the skein theoretical approach and focusing on the $\mathrm{SO}(3)$-theory. In particular, we consider the values at the root of unity  $q=e^{\frac{2\pi\mathbf i}{r}}$ for odd integers $r\geqslant 3.$ To be more precise, we will follow the normalization of the invariants as described in \cite{Li}.

A framed link in an oriented $3$-manifold $M$ is a smooth embedding $L$ of a disjoint union of finitely many thickened circles $\mathrm S^1\times [0,1]$ into $M.$ The Kauffman bracket skein module $\mathrm K_r(M)$ of $M$ is the $\mathbb C$-module generated by the isotopic classes of framed links in $M$  modulo the following two relations: 

\begin{enumerate}[(1)]
\item  \emph{Kauffman Bracket Skein Relation:} \ $\cp{\includegraphics[width=1cm]{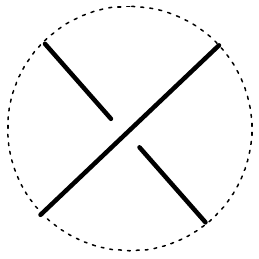}}\ =\ e^{\frac{\pi\mathbf i}{r}}\ \cp{\includegraphics[width=1cm]{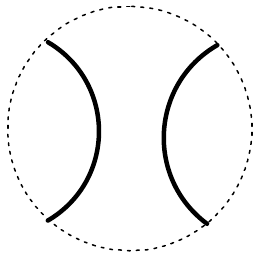}}\  +\ e^{-\frac{\pi\mathbf i}{r}}\ \cp{\includegraphics[width=1cm]{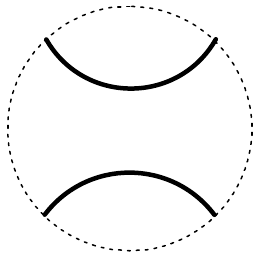}}.$ 

\item \emph{Framing Relation:} \ $L \cup \cp{\includegraphics[width=0.8cm]{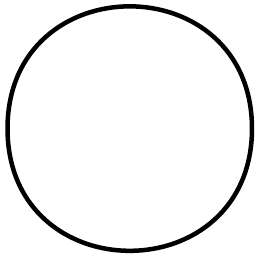}}=\Big(-e^{\frac{2\pi\mathbf i}{r}}-e^{-\frac{2\pi\mathbf i}{r}}\Big)\ L.$ 
\end{enumerate}

There is a canonical isomorphism 
$$\langle\ \ \rangle:\mathrm K_r(\mathrm S^3)\to\mathbb C$$
defined by sending the empty link to $1.$ The image $\langle L\rangle$ of the framed link $L$ is called the Kauffman bracket of $L.$

Let $\mathrm K_r(A\times [0,1])$ be the Kauffman bracket skein module of the product of an annulus $A$ with a closed interval. For any link diagram $D$ in $\mathbb R^2$ with $k$ ordered components and $b_1, \dots, b_k\in \mathrm K_r(A\times [0,1]),$ let 
$$\langle b_1,\dots, b_k\rangle_D$$
be the complex number obtained by cabling $b_1,\dots, b_k$ along the components of $D$ considered as an element of $K_r(\mathrm S^3)$ then taking the Kauffman bracket $\langle\ \rangle.$

On $\mathrm K_r(A\times [0,1])$ there is a commutative multiplication induced by the juxtaposition of annuli, making it a $\mathbb C$-algebra; and as a $\mathbb C$-algebra $\mathrm K_r(A\times [0,1])  \cong \mathbb C[z],$ where $z$ is the core curve of $A.$ For an integer $n\geqslant 0,$ let $e_n(z)$ be the $n$-th Chebyshev polynomial defined recursively by
$e_0(z)=1,$ $e_1(z)=z$ and $e_n(z)=ze_{n-1}(z)-e_{n-2}(z).$ Let $\mathrm{I}_r=\{0,2,\dots,r-3\}$ be the set of even integers in between $0$ and $r-2.$ Then the Kirby coloring $\Omega_r\in\mathrm K_r(A\times [0,1])$ is defined by 
$$\Omega_r=\mu_r\sum_{n\in \mathrm{I}_r}[n+1]e_{n},$$
where
 $$\mu_r=\frac{2\sin\frac{2\pi}{r}}{\sqrt r}$$ 
 and $[n]$ is the quantum integer defined by
$$[n]=\frac{q^n-q^{-n}}{q-q^{-1}}=\frac{\sin\frac{2n\pi}{r}}{\sin\frac{2\pi}{r}}.$$

Let $M$ be a closed oriented $3$-manifold and let $L$ be a framed link in $M$ with $n$ components. Suppose $M$ is obtained from $S^3$ by doing a surgery along a framed link $L',$ $D(L')$ is a standard diagram of $L'$ (i.e., the blackboard framing of $D(L')$ coincides with the framing of $L'$). Then $L$ adds extra components to $D(L')$ forming a linking diagram $D(L\cup L')$ with $D(L)$ and $D(L')$ linking in possibly a complicated way. Let
$U_+$ be the diagram of the unknot with framing $1,$ $\sigma(L')$ be the signature of the linking matrix of $L'$ and $\boldsymbol b=(b_1,\dots,b_n)$ be a multi-element of $I_r.$ Then the $r$-th \emph{relative Reshetikhin-Turaev invariant of $M$ with $L$ colored by $\boldsymbol b$} is defined as
\begin{equation}\label{RT}
\mathrm{RT}_r(M,L,\boldsymbol b)=\mu_r \langle e_{b_1},\dots,e_{b_n}, \Omega_r, \dots, \Omega_r\rangle_{D(L\cup L')}\langle \Omega_r\rangle _{U_+}^{-\sigma(L')}.
\end{equation}

Note that if $L=\emptyset$ or $b_1=\dots=b_n=0,$ then $\mathrm{RT}_r(M,L,\boldsymbol b)=\mathrm{RT}_r(M),$ the $r$-th Reshetikhin-Turaev invariant of $M.$ 


\subsection{Hyperbolic cone manifolds}\label{CCS}

According to \cite{CHK}, a $3$-dimensional \emph{hyperbolic cone-manifold} is a $3$-manifold $M,$ which can be triangulated so that the link of each simplex is piecewise linear homeomorphic to a standard sphere and $M$ is equipped with a complete path metric such that the restriction of the metric to each simplex is isometric to a hyperbolic geodesic simplex. The \emph{singular locus} $L$ of a cone-manifold $M$ consists of the points with no neighborhood isometric to a ball in a Riemannian manifold. It follows that
\begin{enumerate}[(1)]
\item $L$ is a link in $M$ such that each component is a closed geodesic.

\item At each point of $L$ there is a \emph{cone angle} $\theta$ which is the sum of dihedral angles of $3$-simplices containing the point.

\item The restriction of the metric on $M\setminus L$ is a smooth hyperbolic metric, but is incomplete if $L\neq \emptyset.$
\end{enumerate}

Hodgson-Kerckhoff\,\cite{HK} proved that hyperbolic cone metrics on $M$ with singular locus $L$ are locally parametrized by the cone angles provided all the cone angles are less than or equal to $2\pi,$ and Kojima\,\cite{Ko} proved that hyperbolic cone manifolds $(M,L)$ are globally rigid provided all the cone angles are less than or equal to $\pi.$ It is expected to be  globally rigid if all the cone angles are less than or equal to $2\pi.$
\\

For a $3$-manifold $N$ with boundary consisting of tori $T_1,\dots,T_n,$ with a choice of generators $(u_i,v_i)$ for $\pi_1(T_i)$ and an integer $p_i$ for each $i\in\{1,\dots,n\},$  the $(p_1,\dots, p_n)$-Dehn filling on $N$ attaches a solid torus to each $T_i$ so that $p_iu_i+v_i$ bounds a disk.  If $\mathrm H(u_i)$ and  $\mathrm H(v_i)$ are respectively the logarithmic holonomy for $u_i$ and $v_i,$ then a solution to  
\begin{equation}\label{DF}
p_i\mathrm H(u_i)+\mathrm H(v_i)=\mathbf i\theta_i
\end{equation}
 near the complete structure gives a cone-manifold structure on the resulting manifold $M$ with the cone angle $\theta_i$ along the core curve $L_i$ of the solid torus attached to $T_i;$ it is a smooth structure if $\theta_1=\dots=\theta_n= 2\pi.$

In this setting, the Chern-Simons invariant for a hyperbolic cone manifold $(M,L)$ can be defined by using the Neumann-Zagier potential function\,\cite{NZ}. To do this, we need a framing on each component, namely, a choice of a curve $\gamma_i$ on $T_i$ that is isotopic to the core curve $L_i$ of the solid torus attached to $T_i.$ If we choose the orientation of $\gamma_i$ so that $(p_iu_i+v_i)\cdot \gamma_i=1,$ then we can let $\gamma_i=-u_i.$ We consider the following function 
$$\frac{\Phi(\mathrm{H}(u_1),\dots,\mathrm{H}(u_n))}{\mathbf i}-\sum_{i=1}^n\frac{\mathrm H(u_i)\mathrm H(v_i)}{4\mathbf i}+\sum_{i=1}^n\frac{\theta_i\mathrm H(\gamma_i)}{4},$$
where $\Phi$ is the Neumann-Zagier potential function (see \cite{NZ}) defined on the deformation space of hyperbolic structures on $M\setminus L$ parametrized by the holonomy of the meridians $\{\mathrm H(u_i)\},$ characterized by 
\begin{equation}\label{char}
\left \{\begin{array}{l}
\frac{\partial \Phi(\mathrm{H}(u_1),\dots,\mathrm{H}(u_n))}{\partial \mathrm{H}(u_i)}=\frac{\mathrm H(v_i)}{2},\\
\\
\Phi(0,\dots,0)=\mathbf i\bigg(\mathrm{Vol}(M\setminus L)+\mathbf i\mathrm{CS}(M\setminus L)\bigg)\quad\quad\mathrm{mod}\ \pi^2\mathbb Z,
\end{array}\right.
\end{equation} 
where $M\setminus L$ is with the complete hyperbolic metric. Another important feature of $\Phi$ is that it is even in each of its variables $\mathrm H(u_i).$

Following the argument in \cite[Sections 4 \& 5]{NZ}, one can prove that if  the cone angles of components of $L$ are $\theta_1,\dots,\theta_n,$ then
\begin{equation}\label{VOL}
\mathrm{Vol}(M_{L_{\boldsymbol\theta}})=\mathrm{Re}\bigg(\frac{\Phi(\mathrm{H}(u_1),\dots,\mathrm{H}(u_n))}{\mathbf i}-\sum_{i=1}^n\frac{\mathrm H(u_i)\mathrm H(v_i)}{4\mathbf i}+\sum_{i=1}^n\frac{\theta_i\mathrm H(\gamma_i)}{4}\bigg);
\end{equation}
and in \cite{Y}, Yoshida proved that when $\theta_1=\dots=\theta_n=2\pi,$
\begin{equation}\label{Yo}
 \mathrm{Vol}(M)+\mathbf i\mathrm{CS}(M)=\frac{\Phi(\mathrm{H}(u_1),\dots,\mathrm{H}(u_n))}{\mathbf i}-\sum_{i=1}^n\frac{\mathrm H(u_i)\mathrm H(v_i)}{4\mathbf i}+\sum_{i=1}^n\frac{\theta_i\mathrm H(\gamma_i)}{4}\quad\quad\text{mod }\mathbf i\pi^2\mathbb Z.   
\end{equation}

\begin{definition}\cite[Definition 2.1]{WY1}\label{CS} The \emph{Chern-Simons invariant} of a hyperbolic cone manifold $M_{L_{\boldsymbol\theta}}$ with the framing $(\gamma_1,\dots,\gamma_n)$ 
is defined as
$$\mathrm{CS}(M_{L_{\boldsymbol\theta}})=\mathrm{Im}\bigg(\frac{\Phi(\mathrm{H}(u_1),\dots,\mathrm{H}(u_n))}{\mathbf i}-\sum_{i=1}^n\frac{\mathrm H(u_i)\mathrm H(v_i)}{4\mathbf i}+\sum_{i=1}^n\frac{\theta_i\mathrm H(\gamma_i)}{4}\bigg)\mod \pi^2\mathbb Z.$$
\end{definition}

Then together with (\ref{VOL}), we have
\begin{equation}\label{VCS}
\mathrm{Vol}(M_{L_{\boldsymbol\theta}})+\mathbf i\mathrm{CS}(M_{L_{\boldsymbol\theta}})=\frac{\Phi(\mathrm{H}(u_1),\dots,\mathrm{H}(u_n))}{\mathbf i}-\sum_{i=1}^n\frac{\mathrm H(u_i)\mathrm H(v_i)}{4\mathbf i}+\sum_{i=1}^n\frac{\theta_i\mathrm H(\gamma_i)}{4}\quad\quad\text{mod }\mathbf i\pi^2\mathbb Z.
\end{equation}


\subsection{Fundamental shadow link complements}\label{fundsl}

\begin{figure}[htbp]
\centering
\includegraphics[scale=0.25]{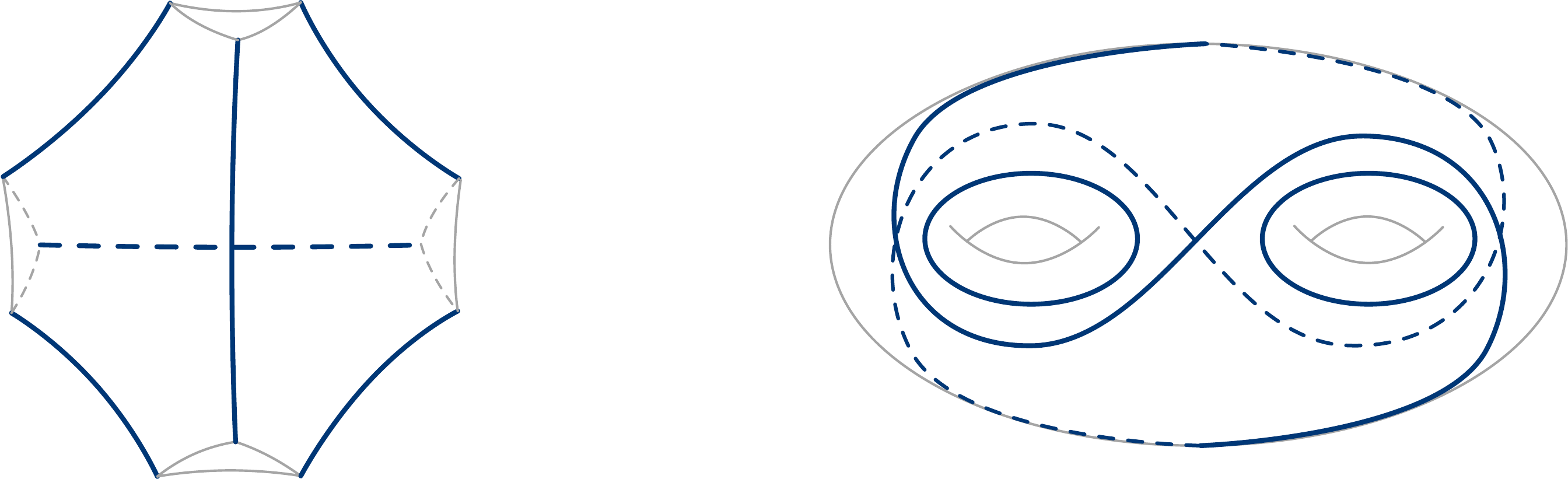}
\caption{The handlebody on the right is obtained from the truncated tetrahedron on the left by identifying the triangles on the top and the bottom by a horizontal reflection and the triangles on the left and the right by a vertical reflection.}
\label{bb}
\end{figure}

If we take $c$ truncated tetrahedra $\Delta_1,\dots, \Delta_c$ and glue them together along the triangles of truncation, we obtain a (possibly non-orientable) handlebody of genus $c+1$ with a link in its boundary consisting of the edges of the truncated tetrahedra. By taking the orientable double (the orientable double
covering with the boundary quotient out by the deck involution) of this handlebody, we obtain a link $L_{\text{FSL}}$ inside $M_c=\#^{c+1}(S^2\times S^1).$ We call a link obtained this way a \emph{fundamental shadow link}, and its complement in $M_c$ a \emph{fundamental shadow link complement}.  Alternatively, to construct a fundamental shadow link complement, we can first  take the double of each tetrahedron along the hexagonal faces and then glue the resulting pieces together along homeomorphisms between the $3$-punctured spheres coming from the double of the triangles of truncation. See Figure \ref{bb}.

The fundamental importance of the family of the fundamental shadow link complements  is the following 

\begin{theorem}\cite{CT} Every compact oriented $3$-manifold with toroidal or empty boundary can be obtained from a suitable fundamental shadow link complement by doing an integral Dehn-filling to some of the boundary components.
\end{theorem}

\begin{figure}[htbp]
\centering
\includegraphics[scale=0.5]{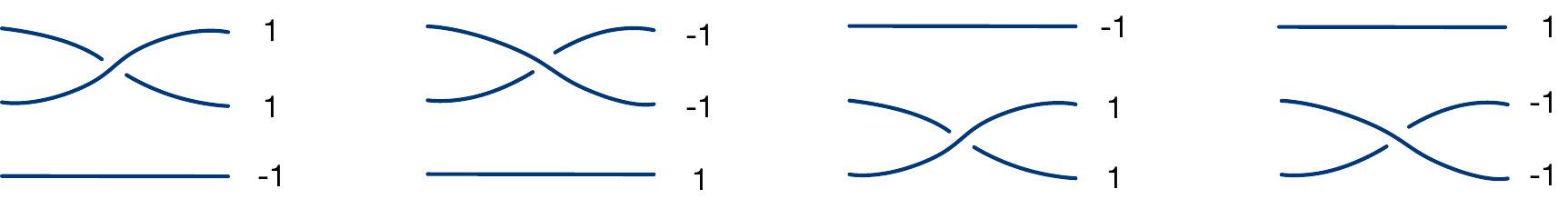}
\caption{ }
\label{mutation}
\end{figure}

A hyperbolic cone metric on $M_c$ with singular locus $L_{\text{FSL}}$ and with cone angles $\theta_1,\dots,\theta_n$ can be constructed as follows. For each $s\in \{1,\dots, c\},$ let $e_{s_1},\dots,e_{s_6}$ be the edges of the building block $\Delta_s,$ and let $\theta_{s_j}$ be half of the cone angle of the component of $L_{\text{FSL}}$ containing $e_{s_j}.$ If $\{{\theta_{s_1}},\dots,{\theta_{s_6}}\}$ form the set of dihedral angles of a truncated hyperideal tetrahedron, by abuse of notation still denoted by $\Delta_s,$ then the hyperbolic cone manifold $M_c$ with singular locus $L_{\text{FSL}}$ and cone angles $\theta_1,\dots,\theta_n,$ where $n$ is the number of components of $L_{\text{FSL}},$ is obtained by gluing $\Delta_s$'s together along isometries of the triangles of truncation, and taking the double. In this metric, the logarithmic holonomy of the meridian $u_i$ of a tubular neighborhood $N(L_i)$ of $L_i$ satisfies 
\begin{equation}\label{m}
\mathrm{H}(u_i)=\mathbf i\theta_i.
\end{equation}
A preferred longitude $v_i$ on the boundary of $N(L_i)$ can be chosen as follows. Recall that a fundamental shadow link is obtained from the double of a set of truncated tetrahedra (along the hexagonal faces) glued together by orientation preserving homeomorphisms between the three-punctured spheres coming from the double of the triangles of truncation, and recall also that 
 the mapping class group of the three-punctured sphere is generated by mutations, which could be represented by the four $3$-braids in Figure \ref{mutation}. For each mutation, we assign an integer $\pm 1$ to each component of the braid as in Figure \ref{mutation}; and for a composition of a sequence of mutations, we assign the sum of the  $\pm 1$ assigned by the mutations to each component of the $3$-braid.
In this way, each orientation-preserving homeomorphism between the three-punctured spheres assigns three integers to three of the components of $L_{\text{FSL}},$ one for each. For each $i\in \{1,\dots, n\},$ let $\iota_i$ be the sum of all the integers on $L_i$ assigned by the homeomorphisms between the three-punctured spheres. Then we can choose a preferred longitude $v_i$ on $N(L_i)$ such that $u_i\cdot v_i=1$ and the logarithmic holonomy satisfies
\begin{equation}\label{l}
\mathrm{H}(v_i)=-l_i+\mathbf i\frac{\iota_i\theta_i}{2},
\end{equation}
where $l_i$ is the length of the closed geodesic $L_i.$  In this way, a framing on $L_i$ gives an integer $p_i$ in the way that the parallel copy of $L_i$ on $N(L_i)$ is isotopic to the curve representing $p_iu_i+v_i.$

 In the complete hyperbolic metric, since $M_c \setminus L_{\text{FSL}}$ is the union of $2c$ regular ideal octahedra, we have the volume
\begin{equation} \label{VolFSL}
\mathrm{Vol}(M_c \setminus L_{\text{FSL}}) = 2c v_8,
\end{equation} 
where $v_8$ is the volume of the regular ideal octahedron; and by \cite[Corollary 2.5]{MO} and \cite[Theorem 2.4]{MR}, the Chern-Simons invariant
\begin{equation}\label{CSFSL}
\mathrm {CS}(M_c \setminus L_{\text{FSL}}) = \Big(\sum_{i=1}^n \frac{\iota_i}{2}\Big) \pi^2  \quad\quad\text{mod }\pi^2\mathbb Z.
\end{equation}

\begin{proposition}[\cite{C1,C2}]\label{FSL} Let $L_{\text{FSL}}=L_1\cup\dots\cup L_n$  be a framed fundamental shadow link in $M_c=\#^{c+1}(S^2\times S^1)$ with framing $p_i$ on the component $L_i,$ and with a coloring $\boldsymbol a=(a_1,\dots,a_n)$ of even integers in $\{0,2,\dots,r-3\}.$ Then the $r$-th relative Reshetikhin-Turaev invariant 
$$\mathrm{RT_r}(M_c,L_{\text{FSL}},\boldsymbol a)=\Bigg(\frac{2\sin\frac{2\pi}{r}}{\sqrt{r}}\Bigg)^{-c}\prod_{i=1}^n(-1)^{\frac{\iota_ia_i}{2}}q^{(p_i+\frac{\iota_i}{2})\frac{a_i(a_i+2)}{2}}\prod_{s=1}^c\bigg|\begin{matrix}
        a_{s_1} & a_{s_2} & a_{s_3} \\
        a_{s_4} & a_{s_5} & a_{s_6} 
      \end{matrix} \bigg|,$$
where $a_{s_1},\dots,a_{s_6}$ are the colors of the edges of the building block $\Delta_s$ inherited  from the color $\boldsymbol a$ on $L_{\text{FSL}},$ and $\footnotesize{\bigg|\begin{matrix}
        a_{s_1} & a_{s_2} & a_{s_3} \\
        a_{s_4} & a_{s_5} & a_{s_6} 
      \end{matrix} \bigg|}$ is the quantum $6j$-symbol of $(a_{s_1},\dots,a_{s_6})$ to be recalled in Section \ref{q6j}.
\end{proposition}

\subsection{The change-of-pair operations}\label{COP}

Let $M$ be a closed oriented $3$-manifold and let $L$ be a framed link in $M.$ The change-of-pair operation is a topological operation that changes the pair $(M,L)$ without changing the complement $M\setminus L,$ and these operations connect all such pairs that have homeomorphic complements.  

\begin{figure}[htbp]
\centering
\includegraphics[scale=0.3]{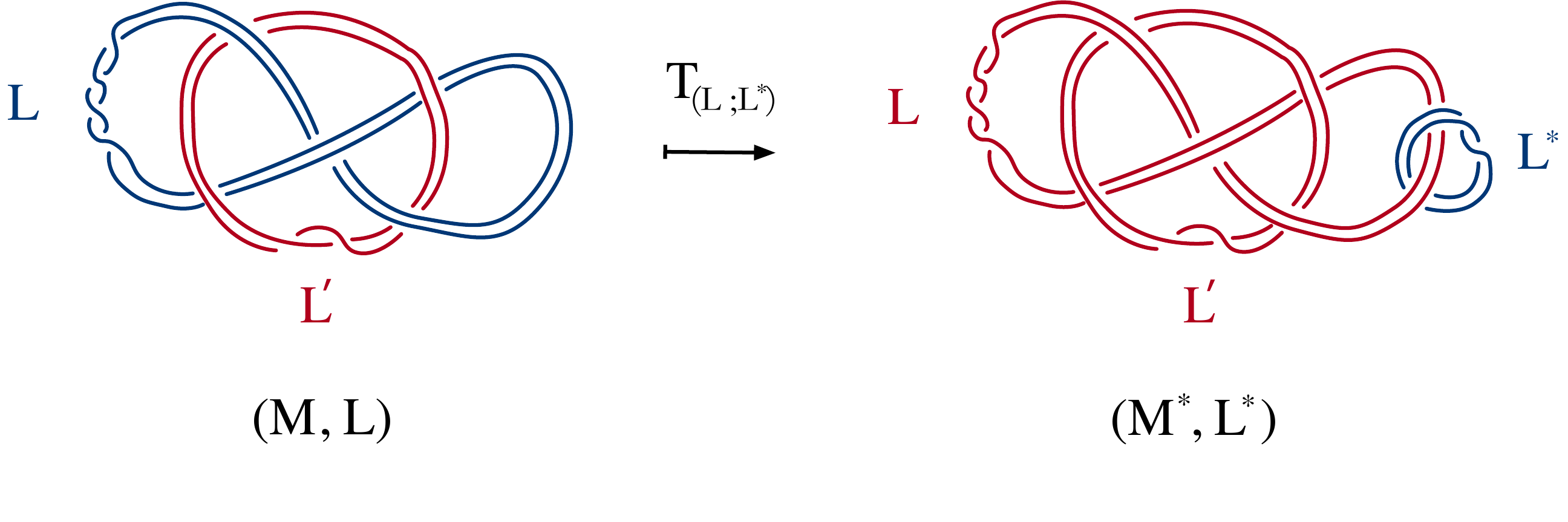}
\caption{$M$ is obtained from $S^3$ by the surgery along $L'.$  $M^*$ is obtained from $M$ by doing the surgery along $L,$ hence is obtained from $S^3$ by doing the surgery along $L'\cup L.$ }
\label{TFT}
\end{figure}

Precisely, suppose $L=L_1\cup\dots\cup L_n.$ For each component  $L_i:[0,1]\times S^1\to M$ of $L,$ recall that $L_i(\{0\}\times S^1)\subset M$ is the core curve of $L_i$ and $L_i(\{1\}\times S^1)\subset M$ is the parallel copy. For each $i,$ let $L_i^*$ be the framed unknot in $M\setminus L$ whose core curve is isotopic to the meridian of the tubular neighborhood of $L_i,$ and let $L^*=L^*_1\cup\dots\cup L^*_n.$ Let $M^*$ be the closed $3$-manifold obtained from $M$ by doing the surgery along $L$ and consider $L^*$ as a framed link in $M^*=M_{L}$ via the inclusions $L^*\subset M\setminus L\subset M_L.$  The \emph{change-of-pair operation} $T_{(L;L^*)}$  is defined by sending $(M,L)$ to  $(M^*,L^*).$ See Figure \ref{TFT}.

\begin{proposition}\cite[Propositions 1.3 and 1.4]{WY1}
\begin{enumerate}[(1)]
\item If $(M^*,L^*)$ is obtained from $(M,L)$ by doing a change-of-pair operation, then $M^*\setminus L^*$ is homeomorphic to $M\setminus L.$

\item Every two pairs $(M,L)$ and $(M^*,L^*)$ such that $M\setminus L$ is homeomorphic to $M^*\setminus L^*$ are related by a sequence of the change-of-pair operations.
\end{enumerate}
\end{proposition}

For the change-of-pair operations in full generality and the relationship with the discrete Fourier transforms, see \cite[Sections 1.2 and 1.3]{WY1}.


\subsection{Quantum $6j$-symbols}\label{q6j}
For $r\geqslant 3$ odd, let $q=e^{\frac{2\pi \mathbf i}{r}}.$ For $n\in \mathbb Z_{\geqslant 0},$ let 
$$[n]=\frac{q^{n}-q^{-n}}{q-q^{-1}}=\frac{\sin \frac{2n\pi}{r}}{\sin \frac{2\pi}{r}},$$
and let 
$$\{n\}!=\prod_{k=1}^n\{k\}\quad\text{and}\quad \{0\}!=0.$$

A triple $(a_1,a_2,a_3)$ of even integers in $\{0,2,\dots,r-3\}$ is \emph{$r$-admissible} if 
\begin{enumerate}[(1)]
\item  $a_i+a_j-a_k\geqslant 0$ for $\{i,j,k\}=\{1,2,3\},$ and
\item $a_1+a_2+a_3\leqslant 2(r-2).$
\end{enumerate}

A  six-tuple $(a_1,\dots,a_6)$  is \emph{$r$-admissible} if the triples $(a_1,a_2,a_3),$ $(a_1,a_5,a_6),$ $(a_2,a_4,a_6)$ and $(a_3,a_4,a_5)$ are $r$-admissible. In the rest of this article, we let 
\begin{equation*}
\begin{split}
   & T_1=\frac{a_1+a_2+a_3}{2},\quad T_2=\frac{a_1+a_5+a_6}{2},\quad T_3=\frac{a_2+a_4+a_6}{2},\quad T_4=\frac{a_3+a_4+a_5}{2},\\
   & Q_1=\frac{a_1+a_2+a_4+a_5}{2},\quad Q_2=\frac{a_1+a_3+a_4+a_6}{2}\quad\text{and}\quad   Q_3=\frac{a_2+a_3+a_5+a_6}{2}.
\end{split}
\end{equation*}

\begin{definition}
The \emph{quantum $6j$-symbol} of an $r$-admissible 6-tuple $(a_1,\dots,a_6)$ is 
\begin{equation*}
\bigg|\begin{matrix} a_1 & a_2 & a_3 \\ a_4 & a_5 & a_6 \end{matrix} \bigg|
= \mathbf i^{-\sum_{k=1}^6a_k}\frac{\prod_{i=1}^4\prod_{j=1}^3\sqrt{[Q_j-T_i]!}}{\prod_{i=1}^4\sqrt{[T_i+1]!}}
\sum_{k=\max \{T_1,T_2,T_3,T_4\}}^{\min\{ Q_1,Q_2,Q_3\}}\frac{(-1)^k[k+1]!}{\prod_{i=1}^4[k-T_i]!\prod_{j=1}^3[Q_j-k]!},
\end{equation*}
where, as a convention,  $\sqrt{-x}\doteq \mathbf i\sqrt{x}$ for $x>0.$
\end{definition}

The quantum $6j$-symbols satisfy the following Tetrahedral Symmetry and Reflection Symmetry, where Proposition \ref{tetra} follows directly from the definition.

\begin{proposition}\label{tetra}
    \begin{align*}
\bigg|\begin{matrix} 
      a_1 & a_2 & a_3 \\
      a_4 & a_5 & a_6
   \end{matrix}\bigg |=\bigg|\begin{matrix} 
      a_2 & a_1 & a_3 \\
      a_5 & a_4 & a_6
   \end{matrix}\bigg|=\bigg|\begin{matrix} 
      a_1 & a_3 & a_2 \\
      a_4 & a_6 & a_5
   \end{matrix}\bigg|=\bigg|\begin{matrix} 
      a_1 & a_5 & a_6 \\
      a_4 & a_2 & a_3
   \end{matrix}\bigg|.
\end{align*}
\end{proposition}

\begin{proposition}\cite[Lemma A.3 (c)]{DKY} 
If the six-tuple $(a_1,a_2,a_3,a_4,a_5,a_6)$ is $r$-admissible, then both of  the six-tuples $$(a_1,a_2,a_3, r-2-a_4, r-2-a_5, r-2-a_6)\quad\text{and}\quad(r-2-a_1,r-2-a_2, a_3, r-2-a_4, r-2-a_5, a_6)$$
are $r$-admissible, and   
$$\bigg|\begin{matrix} a_1 & a_2 & a_3 \\ a_4 & a_5 & a_6 \end{matrix} \bigg|=\bigg|\begin{matrix} a_1 & a_2 & a_3 \\ r-2-a_4 & r-2-a_5 & r-2-a_6 \end{matrix} \bigg|$$
and 
$$\bigg|\begin{matrix} a_1 & a_2 & a_3 \\ a_4 & a_5 & a_6 \end{matrix} \bigg|=\bigg|\begin{matrix} r-2-a_1 & r-2-a_2 & a_3 \\ r-2-a_4 & r-2-a_5 & a_6 \end{matrix} \bigg|.$$
\end{proposition}

\begin{proposition}\cite[Theorem 1.2]{BDKY}
 For any $r$-admissible six-tuple $(a_1,\dots,a_6),$   $$\bigg|\begin{matrix} 
      a_1 & a_2 & a_3 \\
      a_4 & a_5 & a_6
   \end{matrix}\bigg |\leqslant O\Big(re^{\frac{r}{2\pi}v_8}\Big),$$
    where $v_8$ is the volume of the regular ideal octahedron.
\end{proposition}

   \begin{definition} An $r$-admissible six-tuple $(a_1,\dots,a_6)$  is called $r$-\emph{hyperideal} if for $\{i,j,k\}=\{1,2,3\},$ $\{1,5,6\},$ $\{2,4,6\}$ and $\{3,4,5\},$
\begin{enumerate}[(1)]
\item $0\leqslant a_i+a_j-a_k\leqslant r-2,$ and
\item $r-2\leqslant a_i+a_j+a_k\leqslant 2(r-2).$
\end{enumerate}
\end{definition}

Using the same argument as in the proof of \cite[Theorem 1.2]{BDKY}, one can prove the following

   \begin{proposition}\label{non-hyper}\label{bdky} 
   There exists a $\delta_1>0$ such that for any $r$-admissible six-tuple $(a_1,\dots,a_6)$ that is not $r$-hyperideal, 
        $$ \bigg|\begin{matrix} 
      a_1 & a_2 & a_3 \\
      a_4 & a_5 & a_6
 \end{matrix}\bigg| < O\Big(e^{\frac{r}{2\pi}(v_8-\delta_1)}\Big).$$
\end{proposition}


\subsection{$b$-$6j$ symbols}\label{b6jsymbols}

In \cite{LMSWY2}, $b$-$6j$ symbols are defined as a variant of the $6j$-symbols of the positive representations of the modular double of $\mathrm U_q\mathfrak{sl}(2;\mathbb R),$ which are the main building blocks of the $\mathrm U_{q\widetilde q}\mathfrak{sl}(2;\mathbb R)$ Turaev-Viro invariants\,\cite{LMSWY, LMSWY3} for hyperbolic $3$-manifolds with cusp ends or with totally geodesic boundary.

For $b\in (0,1),$ let $Q=b+b^{-1}.$ Let $S_b$ be the  \emph{double sine function} defined for $z\in \mathbb C$ with $0<\mathrm{Re}(z)<Q$ by
\begin{equation}\label{eq:def-S}
S_b(z)=\exp\Bigg(\int_\Omega\frac{\sinh\Big(\big(\frac{Q}{2}-z\big)t\Big)}{4t\sinh(\frac{bt}{2})\sinh(\frac{t}{2b})}dt\Bigg),
\end{equation}
where the contour $\Omega$ goes along the real line and passes above the pole at the origin. For a six-tuple $(a_1,\dots,a_6)\in \mathbb R^6,$ in the rest of this article we let \begin{equation*}
\begin{split}
   & t_1=a_1+a_2+a_3,\quad  t_2=a_1+a_5+a_6,\quad t_3=a_2+a_4+a_6,\quad  t_4=a_3+a_4+a_5,\\
   & q_1=a_1+a_2+a_4+a_5,\quad q_2=a_1+a_3+a_4+a_6,\quad  
q_3=a_2+a_3+a_5+a_6\quad\text{and}\quad q_4=2Q.
\end{split}
\end{equation*}
A six-tuple $(a_1,\dots,a_6)\in\mathbb R^6$ is called \emph{$b$-admissible} if 
$$0< q_j- t_i<Q$$
for all $i,j\in\{1,2,3,4\}.$ 
In particular, if $(a_1,\dots,a_6)$ is $b$-admissible, then 
$$\max\{t_1,t_2,t_3,t_4\}<\min\{q_1,q_2,q_3,q_4\}.$$

\begin{figure}[htbp]
\centering
\includegraphics[scale=0.4]{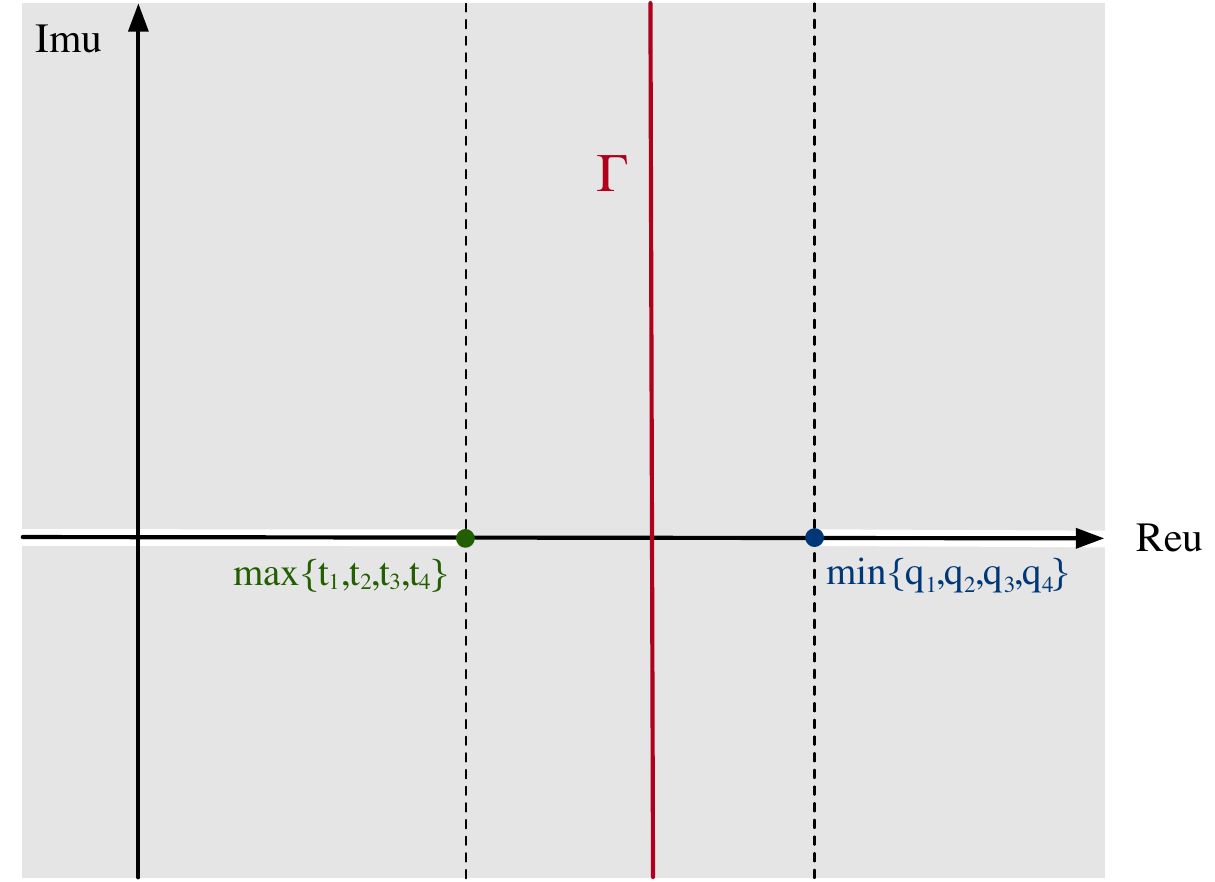}
\caption{The contour $\Gamma.$ }
\label{Gamma}
\end{figure}

\begin{definition}\cite[Definition 1.1]{LMSWY2}.
The \emph{$b$-$6j$ symbol} for a $b$-admissible six-tuple $(a_1,\dots,a_6)\in\mathbb R^6$ is given by 
\begin{equation}\label{b6j}
\bigg\{\begin{matrix} a_1 & a_2 & a_3 \\ a_4 & a_5 & a_6 \end{matrix} \bigg\}_b=\int_\Gamma \frac{\prod_{i=1}^4S_b(u-t_i)\prod_{j=1}^4S_b(q_j-u)}{\prod_{i=1}^4\prod_{j=1}^4\sqrt{S_b(q_j-t_i)}}d\mathrm{Im}u,
\end{equation}
where, as a convention,  the square root of the double sine function is defined as
$$\sqrt{S_b(z)}=\exp\Bigg(\frac{1}{2}\int_\Omega\frac{\sinh\Big(\big(\frac{Q}{2}-z\big)t\Big)}{4t\sinh(\frac{bt}{2})\sinh(\frac{t}{2b})}dt\Bigg),$$
and the contour $\Gamma$ is any vertical line passing the interval $(\max\{t_i\},\min\{q_j\}).$ See Figure \ref{Gamma}.
\end{definition}

\begin{proposition}\cite[Proposition 3.1]{LMSWY2}\label{lem: ab conver}
   The integral in (\ref{b6j}) converges absolutely, and is independent of the choice of the vertical line $\Gamma.$
\end{proposition}

The b-6j symbols satisfy the following Tetrahedral Symmetry and Reflection Symmetry.

\begin{proposition}\cite[Proposition 3.3]{LMSWY2}\label{tetr}
    \begin{align*}
\begin{Bmatrix} 
      a_1 & a_2 & a_3 \\
      a_4 & a_5 & a_6
   \end{Bmatrix}_b=\begin{Bmatrix} 
      a_2 & a_1 & a_3 \\
      a_5 & a_4 & a_6
   \end{Bmatrix}_b=\begin{Bmatrix} 
      a_1 & a_3 & a_2 \\
      a_4 & a_6 & a_5
   \end{Bmatrix}_b=\begin{Bmatrix} 
      a_1 & a_5 & a_6 \\
      a_4 & a_2 & a_3
   \end{Bmatrix}_b.
\end{align*}
\end{proposition}

\begin{proposition}\cite[Propositions 3.4]{LMSWY2}\label{b6jreflection}
If $(a_1,\dots,a_6)$ is $b$-admissible, then $(Q-a_1,a_2,\dots,a_6)$ is also $b$-admissible, and   
$$\bigg\{\begin{matrix} a_1 & a_2 & a_3 \\ a_4 & a_5 & a_6 \end{matrix} \bigg\}_b=\bigg\{\begin{matrix} Q-a_1 & a_2 & a_3 \\ a_4 & a_5 & a_6 \end{matrix} \bigg\}_b.$$
Moreover, the same result holds if we change $a_k\mapsto Q-a_k$ for any $k\in\{1,\dots,6\}.$
\end{proposition}

For a $b$-admissible six-tuple $(a_1,\dots,a_6),$ we will consider the following new set of variables, which will be intensively used later.
\begin{align}
&\alpha_k=2\pi b a_k-\pi b^2 \textrm{ for }k\in\{1,\dots, 6\},\quad \xi=\pi b u, \label{eq:alpha-a}\\
 &\tau_i=\pi b t_i-\frac{3\pi b^2}{2}\textrm{ for }i\in\{1,2,3,4\},\quad\text{and}\quad \eta_j=\pi b q_j-2\pi b^2 \textrm{ for }j\in\{1,2,3,4\}.
\end{align}
Then 
\begin{equation*}
\begin{split}
&\tau_1=\frac{\alpha_1+\alpha_2+\alpha_3}{2}, \quad   \tau_2=\frac{\alpha_1+\alpha_5+\alpha_6}{2}, \quad 
\tau_3=\frac{\alpha_2+\alpha_4+\alpha_6}{2}, \quad \tau_4=\frac{\alpha_3+\alpha_4+\alpha_5}{2},\\    
&\eta_1=\frac{\alpha_1+\alpha_2+\alpha_4+\alpha_5}{2},\quad \eta_2=\frac{\alpha_1+\alpha_3+\alpha_4+\alpha_6}{2}, \quad
\eta_3=\frac{\alpha_2+\alpha_3+\alpha_5+\alpha_6}{2} \quad\text{and}\quad \eta_4=2\pi.
\end{split}    
\end{equation*}
Let $\boldsymbol\alpha=(\alpha_1,\dots,\alpha_6)$ and let
\begin{equation}\label{Ub}
\begin{split}
U_{b}(\boldsymbol \alpha,\xi)= &-\pi \mathbf i b^2 \sum_{i=1}^4\sum_{j=1}^4 \log S_b\bigg(\frac{\eta_j-\tau_i}{\pi b}+\frac{b}{2}\bigg)\\ & + 2\pi \mathbf i b^2 \sum_{i=1}^4  \log S_b\bigg(\frac{\xi-\tau_i-2\pi b^2}{\pi b}+\frac{b}{2}\bigg) + 2\pi \mathbf i b^2 \sum _{j=1}^4\log S_b\bigg(\frac{\eta_j-\xi+\frac{3\pi b^2}{2}}{\pi b}+\frac{b}{2}\bigg).    
\end{split}
\end{equation}
Then we have  
\begin{equation}\label{6jint}
\bigg\{\begin{matrix} a_1 & a_2 & a_3 \\ a_4 & a_5 & a_6 \end{matrix} \bigg\}_b=\frac{1}{\pi \mathbf i b}\int_{\Gamma}\exp\bigg(\frac{U_{b}(\boldsymbol \alpha,\xi)}{2\pi \mathbf ib^2} \bigg)d\xi,
\end{equation}
where $\Gamma$ is a vertical straight line that intersects the real axis in the interval  $\big(\max\{\tau_i\}+\frac{3\pi b^2}{2},\min\{\eta_j\}+2\pi b^2\big)$.


\subsection{Double sine function and dilogarithm function}

By the functional equation (see e.g. \cite[A.15]{TesVar})
\begin{equation}\label{FE1}
S_b\big(z+b^{\pm 1}\big)=2\sin \big(\pi b ^{\pm 1} z\big)S_b(z),
\end{equation}
the double sine function 
$S_b$ extended to a meromorphic function on $\mathbb C$ with the set of poles $\{ -nb-mb^{-1} \ |\ m, n\in \mathbb Z_{\geqslant 0}\}$ and the set of  zeros  $\{ Q + nb + mb^{-1} \ |\ m, n\in \mathbb Z_{\geqslant 0}\}.$ In the rest of this paper, we will intensively use the function $\log S_b\Big(\frac{x}{\pi b}+\frac{b}{2}\Big).$ 
For  $-\frac{\pi b^2}{2}<\mathrm{Re}x<\pi + \frac{\pi b^2}{2}$ so that $0<\mathrm{Re}\Big(\frac{x}{\pi b}+\frac{b}{2}\Big)<Q,$ we let
$$\log S_b\Big(\frac{x}{\pi b}+\frac{b}{2}\Big)=\int_\Omega\frac{\sinh\Big(\big(\frac{1}{2b}-\frac{x}{\pi b}\big)t\Big)}{4t\sinh(\frac{bt}{2})\sinh(\frac{t}{2b})}dt.$$
The function $\log S_b\Big(\frac{x}{\pi b}+\frac{b}{2}\Big)$ can be holomorphically extended to the region $\mathbb C\setminus (-\infty,0]\cup[\pi,\infty).$ 

Let $\log:\mathbb C\setminus (-\infty, 0]\to\mathbb C$ be the standard logarithm function defined by
$$\log z=\log|z|+\mathbf i\arg z$$
with $-\pi<\arg z<\pi.$ The dilogarithm function $\mathrm{Li}_2: \mathbb C\setminus (1,\infty)\to\mathbb C$ is defined by
$$\mathrm{Li}_2(z)=-\int_0^z\frac{\log (1-u)}{u}du$$
where the integral is along any path in $\mathbb C\setminus (1,\infty)$ connecting $0$ and $z,$ which is holomorphic in $\mathbb C\setminus [1,\infty)$ and continuous in $\mathbb C\setminus (1,\infty).$ The following closely related complexified Lobachevsky function 
\begin{equation}\label{eq:Lx}
L(x)=x^2-\pi x +\frac{\pi^2}{6}-\mathrm{Li}_2\big(e^{2\mathbf ix}\big)
\end{equation}
will be used intensively in the rest of the paper, which  is continuous on $\mathbb C\setminus (-\infty,0)\cup(\pi,\infty)$ and holomorphic on $\mathbb C\setminus (-\infty,0]\cup[\pi,\infty).$ The name of the function comes from the fact that for $x\in [0,\pi],$ 
\begin{equation}\label{LL}
L(x)=-2\mathbf i\Lambda(x),    
\end{equation}
where $\Lambda$ is the Lobachevsky function defined by
$$\Lambda(x)=-\int_0^x\log|2\sin t|dt.$$

\begin{proposition}\label{EST} 
\begin{enumerate}[(1)]
\item \cite[Proposition 2.1]{LMSWY}  Let $\nu_b(x)$ be such that
$$2\pi \mathbf i b^2 \log S_b\Big(\frac{x}{\pi b} +\frac{b}{2}\Big)-L(x)=\nu_b(x)b^4.$$ 
Then for any $\epsilon>0,$ there exists a $B_{\epsilon}>0$  independent of $b$ such that for all $b\in(0,1)$ and for all  $x\in \mathbb C$  with $\epsilon<\mathrm{Re}x<\pi-\epsilon,$
$$|\nu_b(x)|<B_\epsilon.$$

\item 
There exists a $B>0$ such that for all $b\in(0,1)$ and for all  $x\in (0,\pi),$ 
$$\bigg|2\pi \mathbf i b^2 \log S_b\Big(\frac{x}{\pi b} +\frac{b}{2}\Big)-L(x)\bigg|<Bb^2.$$
\end{enumerate}
\end{proposition}

\begin{proof} 
We only need to prove (2). By \cite[(B.33)]{Eberhardt}, for $z\in\mathbb C$ with $0<\mathrm{Re}z<Q,$ we have
$$
\log S_b(z)
=\int_0^\infty
\left(
\frac{\sinh\big((\frac Q2-z)t\big)}
{2\sinh\big(\frac{bt}{2}\big)\sinh\big(\frac{t}{2b}\big)}
-\frac{Q-2z}{t}
\right)\frac{dt}{t}.
$$
For $x\in \mathbb C$ with $0<\mathrm{Re}x<\pi,$ we have $0<\mathrm{Re}\big(\frac{x}{\pi b}+\frac{b}{2}\big)<Q;$ and by the change of variable $s=\frac{t}{2\pi b},$ we have 
\begin{equation}\label{216}
2\pi\mathbf i b^2\log S_b\Big(\frac{x}{\pi b}+\frac{b}{2}\Big)
= 2\pi \mathbf i\int_0^\infty
\Bigg(\frac{b^2\sinh\big((\pi-2x)s\big)}
{2\sinh\big(\pi b^2 s\big)\sinh\big(\pi s\big)}
-\frac{(\pi-2x)}{2\pi^2 s}
\Bigg)\frac{ds}{s}.
\end{equation}
Then for any $\epsilon>0$ and $x\in \mathbb C$ with $\epsilon<\mathrm{Re}x<\pi-\epsilon,$ by (1) and letting  $b\to 0,$ we have 
$$L(x)=2\pi \mathbf i\int_0^\infty
\Bigg(\frac{\sinh\big((\pi-2x)s\big)}
{2\pi s\sinh\big(\pi s\big)}
-\frac{(\pi-2x)}{2\pi^2 s}
\Bigg)\frac{ds}{s}.$$
Letting $\epsilon\to 0,$ we have 
\begin{equation}\label{217}
 L(x)=2\pi \mathbf i\int_0^\infty
\Bigg(\frac{\sinh\big((\pi-2x)s\big)}
{2\pi s\sinh\big(\pi s\big)}
-\frac{(\pi-2x)}{2\pi^2 s}
\Bigg)\frac{ds}{s}   
\end{equation}
for all $x\in \mathbb C$ with $0<\mathrm{Re}x<\pi.$
By (\ref{216}) and (\ref{217}), we have 
\begin{equation*}
2\pi\mathbf i b^2\log S_b\Big(\frac{x}{\pi b}+\frac{b}{2}\Big)-L(x)=\mathbf i \int_0^\infty \Bigg(\frac{\pi b^2 s}{\sinh (\pi b^2 s)}-1\Bigg)\frac{\sinh \big((\pi-2x)s\big)}{ s^2\sinh \big(\pi s\big)}ds. 
\end{equation*}
Now, since 
$$0<\frac{\pi b^2 s}{\sinh (\pi b^2 s)}<1$$
for all $s\in (0,+\infty),$
and $$\bigg|\frac{\sinh \big((\pi-2x))s\big)}{\sinh \big(\pi s\big)}\Bigg|<1$$
for all $x\in (0,\pi)$ and $s\in (0,+\infty),$
we have 
\begin{equation*}
\begin{split}
\bigg|2\pi\mathbf i b^2\log S_b\Big(\frac{x}{\pi b}+\frac{b}{2}\Big)-L(x)\bigg|<& \int_0^\infty \Bigg(\frac{1}{s^2}-\frac{ \pi b^2 }{ s\sinh (\pi b^2 s)}\bigg)ds \\
= &  \pi b^2 \int_0^\infty \Bigg(\frac{1}{u^2}-\frac{1}{u\sinh u}\bigg)du
\end{split}
\end{equation*}
for all $x\in (0,\pi),$ where the equality comes from the change of variable $u=\pi b^2 s.$

Finally, let
$$B=\pi\int_0^\infty \Bigg(\frac{1}{u^2}-\frac{1}{u\sinh u}\bigg)du.$$
As $u\to 0^+,$ we have
$$\frac{1}{u^2}-\frac{1}{u\sinh u} = \frac{1}{6} + O(u^2);$$
and as $u\to +\infty,$ we have
$$\frac{1}{u^2}-\frac{1}{u\sinh u} =  O(u^{-2}).$$
Hence the improper integral converges at both $0$ and $+\infty,$ and $B<+\infty.$ 
This completes the proof.
\end{proof}


\section{Extended quantum $6j$-symbols}\label{extq6j}

In this section, we let $r\geqslant 3$ be an odd integer, and let $q=e^{\frac{2\pi \mathbf i}{r}}.$ For $n\in \mathbb Z_{\geqslant 0},$ we let 
$$\{n\}=q^{n}-q^{-n}=2\mathbf i \sin \frac{2n\pi}{r},$$
and let 
$$\{n\}!=\prod_{k=1}^n\{k\}\quad\text{and}\quad \{0\}!=0.$$

\begin{lemma} \label{factorial} Let $b=\sqrt{\frac{2}{r}}\in (0,1),$ and let $S_b$ be the double sine function defined in (\ref{eq:def-S}). 
\begin{enumerate}[(1)]
\item For $0\leqslant n \leqslant r-2,$
\begin{equation}
\{n\}!=\mathbf i^n b^{-1}S_b\big((n+1)b\big).
\end{equation}
\item For $\frac{r-1}{2}\leqslant n \leqslant r-2,$
\begin{equation}  \label{move}
\{n\}!=\frac{(-1)^{n+\frac{r-1}{2}}2\mathbf i^n b^{-1}}{S_b\big((r-n)b\big)}.
\end{equation}
\end{enumerate}
\end{lemma}

\begin{proof} Letting  $z=bk$ in (\ref{FE1}), we have
$$2\sin \Big(\frac{2k\pi}{r}\Big)=2\sin\big(\pi k  b^2\big)=\frac{S_b\big((k+1)b\big)}{S_b(kb)},$$
from which and Lemma \ref{Sb=b} below, (1) follows.

To see (2), by (1), we have
\begin{equation}\label{3.3}
    \{n\}!=\mathbf i^n b^{-1}\frac{S_b\Big(\big(n+1\big)b\Big)} {S_b\Big(\big(n+1-\frac{r}{2}\big)b\Big)}S_b\bigg(\Big(n+1-\frac{r}{2}\Big)b\bigg);
\end{equation}
by (\ref{FE1}) with $z=\big(n+1-\frac{r}{2}\big)b=(n+1)b-b^{-1},$
we have 
\begin{equation}\label{3.4}
\frac{S_b\Big(\big(n+1\big)b\Big)}{S_b\Big(\big(n+1-\frac{r}{2}\big)b\Big)}=2\sin\bigg(\Big(n+1-\frac{r}{2}\Big)\pi\bigg)=(-1)^{n+\frac{r-1}{2}}2;
\end{equation}
and by the functional equation (see e.g. \cite[A.16]{TesVar})
\begin{equation}\label{prod}
S_b(z)S_b(Q-z)=1,
\end{equation}
we have 
\begin{equation}\label{3.5}
S_b\bigg(\Big(n+1-\frac{r}{2}\Big)b\bigg)=\frac{1}{S_b\big((r-n)b\big)}.
\end{equation}
Putting (\ref{3.3}), (\ref{3.4}) and (\ref{3.5}) together, we have (2).
\end{proof}

We consider (\ref{move}) because the function $S_b$ has poles on the right of $Q,$ and to avoid these poles we need to move the variables there into $(0,Q).$

\begin{lemma}\label{Sb=b}
For $b\in (0,1),$ $S_b(b)=b.$
\end{lemma}

\begin{proof} By (\ref{FE1}), we have
$$\frac{S_b(z+b)}{S_b(z+b^{-1})}=\frac{\sin(\pi z b)}{\sin(\pi z b^{-1})};$$
and letting $z\to 0,$ we have 
\begin{equation}\label{3.7}
    \frac{S_b(b)}{S_b(b^{-1})}=b^2.
\end{equation}
By (\ref{prod}), we have 
\begin{equation}\label{3.8}
 {S_b(b)}{S_b(b^{-1})}=1.
\end{equation}
Putting (\ref{3.7}) and (\ref{3.8}) together, we have that $S_b(b)=\pm b.$

It now suffices to show that $S_b(b)>0.$ To see this, by \cite[A.7, A.13]{TesVar}, we have $\overline{S_b(z)}=S_b(\overline z).$ As a consequence, $S_b(z)$ is real for all $z\in (0,Q).$ Since $S_b$ has no zeros in $(0,Q)$ and $S_b\big(\frac{Q}{2}\big)=1>0,$ we have $S_b(z)>0$ for all $z\in (0,Q).$ In particular, $S_b(b)>0.$
\end{proof}

As a direct consequence of Lemma \ref{factorial}, we have

\begin{proposition}\label{6jqd} Let 
$Q_4=r-2.$
Then the quantum $6j$-symbol of an $r$-hyperideal  six-tuple $(a_1,\dots,a_6)$  can be computed as 
$$\bigg|
\begin{matrix}
        a_1 & a_2 & a_3 \\
        a_4 & a_5 & a_6 
      \end{matrix} \bigg|=\frac{(-1)^{\frac{r+1}{2}}2}{r}\sin\frac{2\pi}{r}\sum_{k=\max\{T_1,\dots,T_4\}}^{\min\{Q_1,\dots,Q_4\}}\frac{\prod_{i=1}^4\prod_{j=1}^4\sqrt{S_b\big((Q_j-T_i+1)b\big)}}{\prod_{i=1}^4S_b\big((k-T_i+1)b\big)\prod_{j=1}^4S_b\big((Q_j-k+1)b\big)}.$$
\end{proposition}

\begin{definition} A six-tuple $(a_1,\dots,a_6)$ of real numbers in $\mathbb [0,r-2]$ is called \emph{$r$-hyperideal} if for $\{i,j,k\}=\{1,2,3\},$ $\{1,5,6\},$ $\{2,4,6\}$ and $\{3,4,5\},$
\begin{enumerate}[(1)]
\item $0\leqslant a_i+a_j-a_k\leqslant r-2,$ and
\item $r-2\leqslant a_i+a_j+a_k\leqslant 2(r-2).$
\end{enumerate}
\end{definition}

\begin{definition}\label{eq6j}
The \emph{extended quantum $6j$-symbol} of an $r$-hyperideal six-tuple $(a_1,\dots,a_6)\in [0,r-2]^6$ is defined as
$$\bigg\|
\begin{matrix}
        a_1 & a_2 & a_3 \\
        a_4 & a_5 & a_6 
      \end{matrix} \bigg\|\doteq\frac{(-1)^{\frac{r+1}{2}}2}{r}\sin\frac{2\pi}{r}\sum_k\frac{\prod_{i=1}^4\prod_{j=1}^4\sqrt{S_b\big((Q_j-T_i+1)b\big)}}{\prod_{i=1}^4S_b\big((k-T_i+1)b\big)\prod_{j=1}^4S_b\big((Q_j-k+1)b\big)},$$
      where $Q_4=r-2,$ and the summation is over all the integers $k$ in between $\max\{T_1,T_2,T_3,T_4\}$ and $\min\{Q_1,Q_2,Q_3,Q_4\}.$ As a convention, if there is no integer in between, we define 
      $\footnotesize{\bigg\|
\begin{matrix}
        a_1 & a_2 & a_3 \\
        a_4 & a_5 & a_6 
      \end{matrix}} \bigg\|\doteq 0.$
\end{definition}

The extended quantum $6j$-symbols satisfy the following basic properties. 

\begin{proposition}
    \begin{enumerate}[(1)]
        \item For any $r$-admissible six-tuple $(a_1,\dots,a_6)$ of the hyperideal type,
$$\bigg\|\begin{matrix} 
      a_1 & a_2 & a_3 \\
      a_4 & a_5 & a_6
   \end{matrix}\bigg\|=\bigg|\begin{matrix} 
      a_1 & a_2 & a_3 \\
      a_4 & a_5 & a_6
   \end{matrix}\bigg|.$$

        \item  For any $r$-hyperideal six-tuple $(a_1,\dots,a_6)$ of real numbers,
        \begin{align*}
\bigg\|\begin{matrix} 
      a_1 & a_2 & a_3 \\
      a_4 & a_5 & a_6
   \end{matrix}\bigg\|=\bigg\|\begin{matrix} 
      a_2 & a_1 & a_3 \\
      a_5 & a_4 & a_6
   \end{matrix}\bigg\|=\bigg\|\begin{matrix} 
      a_1 & a_3 & a_2 \\
      a_4 & a_6 & a_5
   \end{matrix}\bigg\|=\bigg\|\begin{matrix} 
      a_1 & a_5 & a_6 \\
      a_4 & a_2 & a_3
   \end{matrix}\bigg\|.
\end{align*}
\end{enumerate}
\end{proposition}

\begin{proof} (1) follows directly from Proposition \ref{6jqd}, and (2) follows directly from the definition.
\end{proof}

The following Proposition \ref{key}, whose proof is postponed to Section \ref{pfKey}, will play the central role in the proof of Theorem \ref{main1} and Theorem \ref{main2}. 

\begin{proposition}\label{key}
There exists a $\delta_2>0$ such that for any $r$-hyperideal six-tuple $(a_1,\dots,a_6)$ of real numbers,
      $$  \Bigg|\bigg\|\begin{matrix} 
      a_1 & a_2 & a_3 \\
      a_4 & a_5 & a_6
 \end{matrix}\bigg\|-\bigg\|\begin{matrix} 
      r-2-a_1 & a_2 & a_3 \\
      a_4 & a_5 & a_6
   \end{matrix}\bigg\|\Bigg| < O\Big(e^{\frac{r}{2\pi}(v_8-\delta_2)}\Big),$$
   where $v_8$ is the volume of the regular ideal octahedron. Moreover, the same result holds if we change $a_k\mapsto r-2-a_k$ for any $k\in\{1,\dots,6\}.$
\end{proposition}

In the rest of this paper, we will intensively use the following change of variables 
\begin{equation*}
\begin{split}
   & \alpha_k=\frac{2\pi (a_k+1)}{r}=\pi b^2a_k+\pi b^2\text{ for }k\in\{1,\dots,6\},\quad \xi=\frac{2\pi k}{r}=\pi b^2k,\\
&\tau_i=\frac{2\pi T_i}{r}+\frac{3\pi}{r}=\pi b^2 T_i+\frac{3\pi b^2}{2}\text{ for }i\in\{1,2,3,4\}, \text{and}\\
&\eta_j=\frac{2\pi Q_j}{r}+\frac{4\pi}{r}=\pi b^2 Q_j+2\pi b^2\text{ for  }j\in\{1,2,3,4\}.
\end{split}
\end{equation*}
Then
\begin{equation*}
\begin{split}
    &\tau_1=\frac{\alpha_1+\alpha_2+\alpha_3}{2},   \quad \tau_2=\frac{\alpha_1+\alpha_5+\alpha_6}{2}, \quad \tau_3=\frac{\alpha_2+\alpha_4+\alpha_6}{2}, \quad \tau_4=\frac{\alpha_3+\alpha_4+\alpha_5}{2},\\
&\eta_1=\frac{\alpha_1+\alpha_2+\alpha_4+\alpha_5}{2}, \ \ \eta_2=\frac{\alpha_1+\alpha_3+\alpha_4+\alpha_6}{2}, 
\ \ \eta_3=\frac{\alpha_2+\alpha_3+\alpha_5+\alpha_6}{2},\quad\text{and}\quad \eta_4=2\pi.   
\end{split}
\end{equation*}

In these new variables, we have
\begin{equation}\label{3.9}
\bigg\|
\begin{matrix}
        a_1 & a_2 & a_3 \\
        a_4 & a_5 & a_6 
      \end{matrix} \bigg\|=e^{\big(\frac{r+1}{2}\big)\pi\mathbf i}\frac{2}{r}\sin\frac{2\pi}{r}\sum_{k}e^{\frac{r}{4\pi\mathbf i}\big(U_r(\alpha_1,\dots,\alpha_6,\xi)\big)},
\end{equation}
where  the summation is over the integers $k$ in the interval  $[\max\{T_1,\dots,T_4\},\min\{Q_1,\dots, Q_4\}],$ and 
\begin{equation}\label{Ur}
\begin{split}
U_r(\alpha_1,\dots,\alpha_6,\xi)= &\frac{1}{2}\sum_{i=1}^4\sum_{j=1}^42\pi\mathbf i b^2 \log S_b\bigg(\frac{\eta_j-\tau_i}{\pi b}+\frac{b}{2}\bigg)\\
&-\sum_{i=1}^42\pi\mathbf i b^2 \log S_b\bigg(\frac{\xi-\tau_i+2\pi b^2}{\pi b}+\frac{b}{2}\bigg)-\sum_{j=1}^42\pi\mathbf i b^2 \log S_b\bigg(\frac{\eta_j-\xi-\frac{3\pi b^2}{2}}{\pi b}+\frac{b}{2}\bigg).
\end{split}    
\end{equation}
In light of Proposition \ref{EST}, we define
\begin{align}
&U(\alpha_1,\dots,\alpha_6,\xi)=\frac{1}{2}\sum_{i=1}^4\sum_{j=1}^4L(\eta_j-\tau_i)-\sum_{i=1}^4L(\xi-\tau_i)-\sum_{j=1}^4L(\eta_j-\xi),\label{U}\\
&\kappa(\alpha_1,\dots,\alpha_6,\xi)=16\pi^2+14\pi\sum_{i=1}^4\tau_i-56\pi\xi-8\pi\mathbf i\sum_{i=1}^4\log \Big(1-e^{2\mathbf i (\xi-\tau_i)}\Big)+6\pi\mathbf i\sum_{{j}=1}^4\log \Big(1-e^{2\mathbf i (\eta_j-\xi)}\Big),\label{kappa}\\
&\nu_r(\alpha_1,\dots,\alpha_6,\xi)=\bigg(U_r(\alpha_1,\dots,\alpha_6,\xi)-U(\alpha_1,\dots,\alpha_6,\xi)-\frac{\kappa(\alpha_1,\dots,\alpha_6,\xi)}{r}\bigg)r^2.\label{nu}
\end{align}
Then we have 
\begin{equation}\label{6jexp}
\bigg\|
\begin{matrix}
        a_1 & a_2 & a_3 \\
        a_4 & a_5 & a_6 
      \end{matrix} \bigg\|=e^{\big(\frac{r+1}{2}\big)\pi\mathbf i}\frac{2}{r}\sin\frac{2\pi}{r}\sum_{k}e^{\frac{r}{4\pi\mathbf i}\big(U(\alpha_1,\dots,\alpha_6,\xi)+\frac{\kappa(\alpha_1,\dots,\alpha_6,\xi)}{r}+\frac{\nu_r(\alpha_1,\dots,\alpha_6,\xi)}{r^2}\big)},
\end{equation}
where the summation is over the integers $k$ in the interval  $[\max\{T_1,\dots,T_4\},\min\{Q_1,\dots, Q_4\}]$ as in (\ref{3.9}).

The new variables also satisfy $\alpha_k\in[\frac{2\pi}{r},2\pi-\frac{2\pi}{r}]$ for each $k\in \{1,\dots,6\},$ and for $\{i,j,k\}=\{1,2,3\},$ $\{1,5,6\},$ $\{2,4,6\}$ and $\{3,4,5\},$ $\frac{2\pi}{r}\leqslant \alpha_i+\alpha_j-\alpha_k\leqslant 2\pi-\frac{2\pi}{r},$ and $2\pi+\frac{2\pi}{r}\leqslant \alpha_i+\alpha_j+\alpha_k\leqslant 4\pi-\frac{2\pi}{r}.$ These conditions  suggest us to make the following

\begin{definition} A six-tuple $(\alpha_1,\dots,\alpha_6)\in \mathbb (0,2\pi)^6$ is of the \emph{hyperideal type} if for $\{i,j,k\}=\{1,2,3\},$ $\{1,5,6\},$ $\{2,4,6\}$ and $\{3,4,5\},$
\begin{enumerate}[(1)]
\item $0< \alpha_i+\alpha_j-\alpha_k<2\pi,$ and
\item $2\pi<\alpha_i+\alpha_j+\alpha_k<4\pi.$
\end{enumerate}
\end{definition}
In the rest of this paper, we let 
$$\mathrm B=\Big\{ (\alpha_1,\dots,\alpha_6,\xi)\in (0,2\pi)^6\times \mathbb R  \ \Big|\ (\alpha_1,\dots,\alpha_6) \text{ is of the hyperideal type, and} \max\{\tau_i\}<\xi<\min\{\eta_j\}\Big\};$$
and for $\epsilon>0,$ let 
$$\mathrm B_\epsilon =\Big\{ (\alpha_1,\dots,\alpha_6,\xi)\in \mathrm B  \ \Big|\ d\big((\alpha_1,\dots,\alpha_6) ,\partial \mathrm B\big)>\epsilon\Big\}.$$
By (\ref{LL}), for $(\alpha_1,\dots,\alpha_6,\xi)$ in the region $\mathrm{B},$ we have
\begin{equation}\label{UV} 
U(\alpha_1,\dots,\alpha_6,\xi)=2\mathbf i V(\alpha_1,\dots,\alpha_6,\xi),
\end{equation}
where $V:\mathrm{B}\to\mathbb R$ is given by 
\begin{equation}\label{V}
V(\alpha_1,\dots,\alpha_6,\xi)\doteq-\frac{1}{2}\sum_{i=1}^4\sum_{j=1}^4\Lambda(\eta_j-\tau_i)+\sum_{i=1}^4\Lambda(\xi-\tau_i)+\sum_{j=1}^4\Lambda(\eta_j-\xi).
\end{equation}
A result of Costantino\,\cite{C1} shows that for each $\boldsymbol\alpha=(\alpha_1,\dots,\alpha_6)$ of the hyperideal type, there exists a unique $\xi_{\boldsymbol\alpha}$ so that $(\boldsymbol\alpha,\xi_{\boldsymbol\alpha})\in\mathrm{B}$ and 
\begin{equation}\label{criticalequation}
\frac{\partial V(\boldsymbol\alpha,\xi)}{\partial \xi}\Big|_{\xi=\xi_{\boldsymbol\alpha}}=0.
\end{equation}
Indeed, he proves that for each fixed $\boldsymbol\alpha,$ $V$ is strictly concave down in $\xi,$ i.e.,
\begin{equation}\label{convexU}
 \frac{\partial^2 V(\boldsymbol \alpha, \xi)}{\partial\xi^2} <0,
\end{equation}
and with derivatives $\pm\infty$ at the boundary points of the interval of $\xi,$ hence there is a unique critical point $\xi_{\boldsymbol\alpha}$ at which $V$ achieves the absolute maximum.
Moreover, by using the Murakami-Yano formula\,\cite{MY,U} he shows that 
\begin{equation}\label{VV}
V(\boldsymbol\alpha,\xi_{\boldsymbol\alpha})=\mathrm{Vol}(\Delta_{\boldsymbol \alpha}),
\end{equation}
where $\Delta_{\boldsymbol \alpha}$ is the truncated hyperideal tetrahedron with dihedral angles $|\alpha_1-\pi|,\dots, |\alpha_6-\pi|.$ 
As a consequence, for any $(\boldsymbol \alpha ,\xi)\in \mathrm B,$ we have
\begin{equation}\label{320}
V(\boldsymbol \alpha,\xi)\leqslant V(\boldsymbol \alpha,\xi_{\boldsymbol\alpha})=\mathrm{Vol}(\Delta_{\boldsymbol \alpha})\leqslant v_8=V\Big(\boldsymbol \pi, \frac{7\pi}{4}\Big),    
\end{equation}
where $v_8$ is the volume of the regular ideal octahedron; 
and the equality holds if and only if $(\boldsymbol \alpha,\xi)=\big(\boldsymbol \pi,\frac{7\pi}{4}\big).$ This, in particular, implies that 
\begin{equation}\label{boundaryv8}
\max_{(\boldsymbol \alpha,\xi)\in \mathrm B\setminus \mathrm B_\epsilon} V(\boldsymbol \alpha,\xi)< v_8 
\end{equation}
for any sufficiently small $\epsilon>0;$
and for each $\boldsymbol \alpha$ of the hyperideal type,
\begin{equation}\label{nbhdvol}
\max_{\xi\in (\max\{\tau_i\},\min\{\eta_j\})\setminus [\xi_{\boldsymbol \alpha}-\epsilon, \xi_{\boldsymbol \alpha}+\epsilon]}V(\boldsymbol \alpha,\xi) < \mathrm{Vol}(\Delta_{\boldsymbol \alpha})
\end{equation}
for any sufficiently small $\epsilon>0.$

The function $U$  extends holomorphically to 
$$\mathrm B_{\mathbb C}=\Big\{ (\alpha_1,\dots,\alpha_6,\xi)\in \mathbb C^7\ \Big|\  (\mathrm{Re}\alpha_1,\dots,\mathrm{Re}\alpha_6,\mathrm\xi)\in \mathrm B   \Big\}.$$
For $\boldsymbol\alpha=(\alpha_1,\dots,\alpha_6)\in\mathbb C^6$ such that $(\mathrm{Re}\alpha_1,\dots,\mathrm{Re}\alpha_6)$ is of the hyperideal type, we let $\xi_{\boldsymbol \alpha}$ be such that
\begin{equation}\label{xia}
\frac{\partial U(\boldsymbol\alpha,\xi)}{\partial \xi}\Big|_{\xi=\xi_{\boldsymbol \alpha}}=0.
\end{equation}
Following the idea of \cite{MY}, see also \cite{BY}, it is proved that there is at most one $\xi_{\boldsymbol \alpha}$ such that $(\boldsymbol\alpha,\xi_{\boldsymbol \alpha})\in \mathrm{B_{\mathbb C}};$ and it is not clear that   $(\boldsymbol\alpha,\xi_{\boldsymbol \alpha})$ is always in $\mathrm{B_{\mathbb C}}.$


\section{Computation of the Reshetikhin-Turaev invariants}\label{comp}

\subsection{The relative case}

Let $L_{\text{FSL}}=L_1\cup\dots\cup L_n$  be a fundamental shadow link in $M_c=\#^{c+1}(S^2\times S^1),$ and let $U_{c+1}\subset S^3$ be the disjoint union of $c+1$ unknots with the $0$-framings by doing surgery along which we get $M_c.$ Let $(M,L)$ be the pair obtained from $(M_c,L_{\text{FSL}})$ by doing a change-of-pair operation $T_{(L_\text{FSL};L)}$ as introduced in Section \ref{COP}, i.e., $M=(M_c)_{L_\text{FSL}}$ and $L=\cup_{i=1}^n L^*_i$ and $L^*_i$ is the framed unknot in $M_c\setminus L_{\text{FSL}}$ with the core curve isotopic to the meridian of the tubular neighborhood of $L_i.$ Suppose $L^*_i$ has the framing $q_i$ and $L_i$ has the framing  $p_i$  for each $i\in \{1,\dots,n\}.$
Let $\boldsymbol b^{(r)}$ be a sequence of colorings of $L,$ where $b_i^{(r)} \in \{0,2,\dots, r-3\}$ is an even integer for each $i.$ Assume that for each $i,$ either $b_i^{(r)}>\frac{r-1}{2}$ for all $r$ or $b_i^{(r)}<\frac{r-1}{2}$ for all $r.$ In the former case we let $\mu_i=1$ and in the latter case we let $\mu_i=-1,$ and we let
$$\theta_i^{(r)}=\mu_i\bigg(\frac{4\pi (b_i^{(r)}+1)}{r}-2\pi\bigg)=\bigg|\frac{4\pi (b_i^{(r)}+1)}{r}-2\pi\bigg|.$$
 Let $\boldsymbol \theta^{(r)}=(\theta_1^{(r)},\dots,\theta_n^{(r)}),$ and let $M^{(r)}=M_{L_{\theta^{(r)}}}$ be the hyperbolic cone manifold consisting of $M$ and a hyperbolic cone metric on $M$  with singular locus $L$ and cone angles $\boldsymbol \theta^{(r)},$ if it exists, and let $\mathrm {Vol}(M^{(r)})$ and $\mathrm {CS}(M^{(r)})$ respectively be the volume and the Chern-Simons invariant of $M^{(r)}.$ 
For $i\in \{1,\dots, n\},$ let $u_i$ and $v_i$ respectively be the meridian and the preferred longitude of a tubular neighborhood of $L_i,$ let $m_i=p_iu_i+v_i$ be the meridian of a tubular neighborhood of $L_i^*$ and let $\mathbf m=(m_i)_{i\in \{1,\dots,n\}}.$  Let $\rho_{M^{(r)}}:\pi_1(M\setminus L)\to\mathrm{PSL}(2;\mathbb C)$ be the holonomy representation of the restriction of hyperbolic cone metric $M^{(r)}$ to $M\setminus L,$ and let  $\mathbb{T}_{(M\setminus L,\mathbf m)}([\rho_{M^{(r)}}])$ be the Reidemeister torsion of $M\setminus L$ twisted by the adjoint action of $\rho_{M^{(r)}}.$

Also let 
$$\theta_i=\lim_{r\to \infty}\theta_i^{(r)}$$  for $i\in \{1,\dots,n\},$  
 let  $\boldsymbol \theta=(\theta_1,\dots,\theta_n),$ let $M_{L_\theta}$ be the hyperbolic cone manifold consisting of $M$ and a hyperbolic cone metric on $M$  with singular locus $L$ and cone angles $\boldsymbol \theta,$ if it exists, and let $\mathrm {Vol}(M_{L_\theta})$ and $\mathrm {CS}(M_{L_\theta})$ respectively be the volume and the Chern-Simons invariant of $M_{L_\theta}.$ Then Theorem \ref{main1} can be rephrased  as follows.

\begin{theorem}\label{re1} There exists a $\delta>0$ such that if
$$\mathrm{Vol}(M_{L_{\boldsymbol \theta}})>\mathrm{Vol}(M_c\setminus L_\text{FSL})-\delta,$$
then the following hold:
\begin{enumerate}[(1)]
\item
    As $r$ varies over all positive odd integers, $$\mathrm{RT}_r\big(M,L,\boldsymbol b^{(r)}\big)=\frac{C_r}{\sqrt{\pm\mathbb{T}_{(M\setminus L,\mathbf m)}([\rho_{M^{(r)}}])}}e^{\frac{r}{4\pi}\big(\mathrm{Vol}(M^{(r)})+\mathbf i\mathrm{CS}(M^{(r)})\big)}\Big( 1 + O \Big( \frac{1}{r} \Big) \Big),$$
where $C_r$
is a quantity of norm $1,$ which is independent of the geometric structure on $M.$

\item 
As $r$ varies over all positive odd integers, 
$$\lim_{r\to\infty} \frac{4\pi }{r}\log \mathrm{RT}_r\big(M,L,\boldsymbol b^{(r)}\big)=\mathrm{Vol}(M_{L_{\boldsymbol\theta}})+\mathbf i\mathrm{CS}(M_{L_{\boldsymbol\theta}})\mod \mathbf i\pi^2\mathbb Z.$$ 
\end{enumerate}
\end{theorem}

Clearly, (2) is a direct consequence of (1); and the equivalence  between (1) and the original statement of Conjecture \ref{CAsymp} can be seen from Remark \ref{rmk}.

In the rest of this section, we will drop the superscript and write $\boldsymbol b=\boldsymbol b^{(r)}$ for the simplicity of the notations, and keep in mind that it is a quantity depending on $r.$ 
 Then the $r$-th relative Reshetikhin-Turaev invariant of $(M,L)$ with coloring  $\boldsymbol b$ can be computed as
\begin{equation*}\label{rt}
\mathrm{RT}_r\big(M,L,\boldsymbol b\big)=A_{\boldsymbol b,r}\sum_{\boldsymbol a}h_{\boldsymbol b,r}(\boldsymbol a),
\end{equation*}
where
$$A_{\boldsymbol b,r}=\bigg(\frac{2\sin\frac{2\pi}{r}}{\sqrt{r}}\bigg)^{n-c}e^{-\sigma(U_{c+1}\cup L)\big(-\frac{3}{r}-\frac{r+1}{4}\big)\pi\mathbf i}\prod_{i=1}^nq^{\frac{q_ib_i(b_i+2)}{2}},$$
the sum is over all multi-even integers $\boldsymbol a=(a_i)_{i\in \{1,\dots,n\}}$ in $\{0,2,\dots,r-3\},$ 
$$h_{\boldsymbol b,r}(\boldsymbol a)=\prod_{i=1}^n(-1)^{\frac{\iota_ia_i}{2}}q^{\big(p_i+\frac{\iota_i}{2}\big)\frac{a_i(a_i+2)}{2}}[(a_i+1)(b_i+1)]\prod_{s=1}^c\bigg|\begin{matrix}
        a_{s_1} & a_{s_2} & a_{s_3} \\
        a_{s_4} & a_{s_5} & a_{s_6} 
      \end{matrix} \bigg|,$$
and $a_{s_1},\dots,a_{s_6}$ are the colors of the edges of the building block $\Delta_s$ inherited  from the colors on $L_{\text{FSL}},$ with the convention that $\footnotesize{\bigg|\begin{matrix}
        a_{s_1} & a_{s_2} & a_{s_3} \\
        a_{s_4} & a_{s_5} & a_{s_6} 
      \end{matrix} \bigg|}=0$ if the six-tuple $(a_{s_1},\dots,a_{s_6})$ is not $r$-admissible. As $a_i$ is even for each $i,$   $(-1)^{(p-\iota_i+1)a_i}=1,$ multiplying which to $h_{\boldsymbol b,r}(\boldsymbol a)$ we have 
$$h_{\boldsymbol b,r}(\boldsymbol a)=\prod_{i=1}^n(-1)^{(p_i-\frac{\iota_i}{2})a_i}q^{\big(p_i+\frac{\iota_i}{2}\big)\frac{a_i(a_i+2)}{2}}(-1)^{a_i}[(a_i+1)(b_i+1)]\prod_{s=1}^c\bigg|\begin{matrix}
        a_{s_1} & a_{s_2} & a_{s_3} \\
        a_{s_4} & a_{s_5} & a_{s_6} 
      \end{matrix} \bigg|.$$      

Now let 
$$\mathrm H_{r}=\Big\{ \boldsymbol a \in [0,r-2]^n\ \Big|\ \text{for each building block }\Delta_s, (a_{s_1},\dots,a_{s_6}) \text{ is $r$-hyperideal}  \Big\}.$$ Then by Proposition \ref{bdky} and the constant $\delta_1$ therein, we have
\begin{equation*}\label{rt2}
\mathrm{RT}_r\big(M,L,\boldsymbol b\big)=A_{\boldsymbol b,r}\sum_{\boldsymbol a\in \mathrm H_{r}\cap (2\mathbb Z)^n}h_{\boldsymbol b,r}(\boldsymbol a)+O\Big(e^{\frac{r}{4\pi}\big(\mathrm{Vol}(M_\text{FSL})-\delta_1\big)}\Big).
\end{equation*}
On $\mathrm H_{r},$ we define the function 
$$h_{\boldsymbol b,r}(\boldsymbol a)\doteq \prod_{i=1}^n(-1)^{(p_i-\frac{\iota_i}{2})a_i}q^{\big(p_i+\frac{\iota_i}{2}\big)\frac{a_i(a_i+2)}{2}}(-1)^{a_i}[(a_i+1)(b_i+1)]\prod_{s=1}^c\bigg\|\begin{matrix}
        a_{s_1} & a_{s_2} & a_{s_3} \\
        a_{s_4} & a_{s_5} & a_{s_6} 
      \end{matrix} \bigg\|,$$
where $\footnotesize{\bigg\|\begin{matrix}
        a_{s_1} & a_{s_2} & a_{s_3} \\
        a_{s_4} & a_{s_5} & a_{s_6} 
\      \end{matrix} \bigg\|}$ is the extended quantum $6j$-symbol of $(a_{s_1},\dots,a_{s_6})$ defined in Definition \ref{eq6j}.
Then by Proposition \ref{key} and the constant $\delta_2$ therein, we have 
\begin{equation}\label{rt3}
\mathrm{RT}_r\big(M,L,\boldsymbol b\big)=\frac{A_{\boldsymbol b,r}}{2^{n}}\sum_{\boldsymbol a\in \mathrm H_{r}\cap \mathbb Z^n}h_{\boldsymbol b,r}(\boldsymbol a)+O\Big(e^{\frac{r}{4\pi}\big(\mathrm{Vol}(M_\text{FSL})-\delta_0\big)}\Big)
\end{equation}
for any $\delta_0 < \min\{\delta_1,\delta_2\}.$

\begin{proposition}\label{rt4}
\begin{equation*}
\mathrm{RT}_r\big(M,L,\boldsymbol b\big)=B_r\sum_{\boldsymbol a,\boldsymbol k}\bigg(\sum_{\boldsymbol \epsilon}g^{\boldsymbol \epsilon}_{\boldsymbol b,r}(\boldsymbol a,\boldsymbol k)\bigg)+O\Big(e^{\frac{r}{4\pi}\big(\mathrm{Vol}(M_\text{FSL})-\delta_0\big)}\Big),
\end{equation*}
where 
$$B_r=2^{-n}r^{-\frac{n+c}{2}}e^{\big(-\sigma(U_{c+1}\cup L)(-\frac{3}{r}-\frac{r+1}{4})-(\frac{r}{4}+1+\frac{1}{r})\sum_{i=1}^n(p_i+q_i)+(\frac{1}{2}-\frac{1}{2r})\sum_{i=1}^n\iota_i-\frac{r}{2}(n-c)-\frac{n}{2}\big)\pi\mathbf i},$$
$\boldsymbol \epsilon=(\epsilon_i)_{i\in \{1,\dots,n\}}\in\{-1,1\}^n$ runs over all multi-signs, $\boldsymbol a=(a_i)_{i\in \{1,\dots,n\}}$ runs over multi-integers in $\mathrm H_{r},$ and $\mathbf k=(k_s)_{s\in \{1,\dots,c\}}$ runs over all multi-integers with each $k_s$ lying in between $\max\{T_{s_1},T_{s_2},T_{s_3},T_{s_4}\}$ and $\min\{Q_{s_1},Q_{s_2},Q_{s_3},Q_{s_4}\},$
with 
$$g^{\boldsymbol \epsilon}_{\boldsymbol b,r}(\boldsymbol a,\boldsymbol  k)=e^{\frac{r}{4\pi\mathbf i}\mathcal W_{\boldsymbol \beta,r}^{\boldsymbol\epsilon}\big(\frac{2\pi (\boldsymbol a+\boldsymbol 1)}{r},\frac{2\pi\boldsymbol k}{r}\big)},$$
where $\frac{2\pi (\boldsymbol a+\boldsymbol 1)}{r}=\big(\frac{2\pi (a_i+1)}{r}\big)_{i\in \{1,\dots,n\}}$ and   $\frac{2\pi\boldsymbol k}{r}=\big(\frac{2\pi k_s}{r}\big)_{s\in \{1,\dots, c\}},$ and with $\boldsymbol \beta=(\beta_i)_{i\in\{1,\dots,n\}}=\big(\frac{2\pi (b_i+1)}{r}\big)_{i\in \{1,\dots,n\}},$ $\boldsymbol{\alpha}=(\alpha_i)_{i\in\{1,\dots,n\}},$ $\boldsymbol\xi=(\xi_s)_{s\in\{1,\dots,c\}},$
and
\begin{equation*}
\begin{split}
\mathcal W_{\boldsymbol \beta,r}^{\boldsymbol \epsilon}(\boldsymbol{\alpha},\boldsymbol{\xi})=&-\sum_{i=1}^nq_i(\beta_i-\pi)^2-\sum_{i=1}^np_i(\alpha_i-\pi)^2-\sum_{i=1}^n2\epsilon_i(\alpha_i-\pi)(\beta_i-\pi)\\
&-\sum_{i=1}^n\frac{\iota_i}{2}(\alpha_i-\pi)^2+\sum_{s=1}^c U_r(\alpha_{s_1},\dots,\alpha_{s_6},\xi_s)+\Big(\sum_{i=1}^n\frac{\iota_i}{2}\Big)\pi^2.
\end{split}
\end{equation*}
\end{proposition}

\begin{proof}
For an integer $a$ and an even integer $b,$ let $\alpha=\frac{2\pi (a+1)}{r}$ and $\beta=\frac{2\pi(b+1)}{r}.$ Then we have
\begin{equation*}
\begin{split}
    (-1)^aq^{(a+1)(b+1)}=&\,q^{(a+1-\frac{r}{2})(b+1-\frac{r}{2})+(a+1)r+\frac{br}{2}-\frac{r^2}{4}}=e^{\big(-\frac{r}{2}\big)\pi\mathbf i}e^{\frac{r}{4\pi\mathbf i}\big(-2(\alpha-\pi)(\beta-\pi)\big)},\\
      (-1)^aq^{-(a+1)(b+1)}=&\,q^{-(a+1-\frac{r}{2})(b+1-\frac{r}{2})-(a+1)r-\frac{br}{2}+\frac{r^2}{4}}=-e^{\big(-\frac{r}{2}\big)\pi\mathbf i}e^{\frac{r}{4\pi\mathbf i}\big(2(\alpha-\pi)(\beta-\pi)\big)}.
\end{split}    
\end{equation*}
As a consequence, we have 
\begin{equation*}
 \prod_{i=1}^n(-1)^{a_i}[(a_i+1)(b_i+1)]=\frac{(-1)^{-\frac{r}{2}}}{\{1\}^{n}}\sum_{\boldsymbol \epsilon\in\{-1,1\}^n}e^{\frac{r}{4\pi\mathbf i}\sum_{i=1}^n\big(-2\epsilon_i(\alpha_i-\pi)(\beta_i-\pi)\big)}.
\end{equation*}
Then the result follows from (\ref{rt3}), (\ref{3.9}) and the following computations that for any integer $a,$
\begin{equation*}
\begin{split}
    (-1)^aq^{\frac{a(a+2)}{2}}=&\,q^{\frac{1}{2}\big((a+1-\frac{r}{2})^2+2(a+1)r-\frac{r^2}{4}-r-1\big)}=e^{\big(-\frac{r}{4}-1-\frac{1}{r}\big)\pi\mathbf i}e^{\frac{r}{4\pi\mathbf i}\big(-(\alpha-\pi)^2\big)},\\
    (-1)^{-\frac{a}{2}}q^{\frac{a(a+2)}{4}}=&\,q^{\frac{1}{4}\big((a+1-\frac{r}{2})^2-\frac{r^2}{4}+r-1\big)}=e^{\big(\frac{1}{2}-\frac{1}{2r}\big)\pi\mathbf i}e^{\frac{r}{4\pi\mathbf i}\big(-\frac{(\alpha-\pi)^2}{2}+\frac{\pi^2}{2}\big)}.
\end{split}    
\end{equation*}
 \end{proof}

Next we apply  the Poisson Summation Formula. Notice that the summation in Proposition \ref{rt4} is finite, and we need a summation over all the integer points to use the Poisson Summation Formula. For this reason, we consider the following regions and a bump function.

Let $\beta_i=\frac{2\pi(b_i+1)}{r}$ and $\alpha_i=\frac{2\pi(a_i+1)}{r}$ for $i\in\{1,\dots,n\},$ $\xi_s=\frac{2\pi k_s}{r}$ for $s\in\{1,\dots,c\},$ $\tau_{s_i}=\frac{2\pi T_{s_i}}{r}+\frac{3\pi}{r}$ for $i\in\{1,2,3,4\},$ and $\eta_{s_j}=\frac{2\pi Q_{s_j}}{r}+\frac{4\pi}{r}$ for $j\in\{1,2,3,4\}.$ We let 
$$\mathrm D=\Big\{ (\boldsymbol \alpha,\boldsymbol \xi) \ \Big|\ (\alpha_{s_1},\dots,\alpha_{s_6},\xi_s)\in\mathrm B\text{ for each }s\in\{1,\dots,c\}\Big\},$$
and for a sufficiently small $\epsilon >0,$ let 
$$\mathrm D_{\epsilon}=\Big\{(\boldsymbol{\alpha},\boldsymbol{\xi})\in\mathrm D\ \Big|\ d\big((\boldsymbol{\alpha},\boldsymbol{\xi}), \partial\mathrm D\big)>\epsilon \Big\},$$
where $d$ is the Euclidean distance on $\mathbb R^{n+c}.$
Let $\psi:\mathbb R^{n+c}\to\mathbb R$ be a $C^{\infty}$-smooth bump function supported on $(\mathrm D, \mathrm D_{\epsilon}),$ i.e.,  \begin{equation*}
\left \{\begin{array}{rl}
\psi(\boldsymbol{\alpha},\boldsymbol{\xi})=1, & (\boldsymbol{\alpha},\boldsymbol{\xi})\in \overline{\mathrm D_{\epsilon}}\\
0<\psi(\boldsymbol{\alpha},\boldsymbol{\xi})<1, &  (\boldsymbol{\alpha},\boldsymbol{\xi})\in \mathrm D\setminus \overline{\mathrm D_{\epsilon}}\\
\psi(\boldsymbol{\alpha},\boldsymbol{\xi})=0, & (\boldsymbol{\alpha},\boldsymbol{\xi})\notin \mathrm D,\\
\end{array}\right.
\end{equation*}
 let 
$$f^{\boldsymbol \epsilon}_{\boldsymbol b,r}(\boldsymbol a,\boldsymbol k)=\psi\Big(\frac{2\pi (\boldsymbol a+\boldsymbol 1)}{r},\frac{2\pi \boldsymbol k}{r}\Big)g^{\boldsymbol \epsilon}_{\boldsymbol b,r}(\boldsymbol a,\boldsymbol k),$$
and let $$f_{\boldsymbol b,r}=\sum_{\boldsymbol \epsilon\in\{-1,1\}^{n}} f_{\boldsymbol b,r}^{\boldsymbol \epsilon}.$$ Then
$$\mathrm{RT}_r\big(M,L,\boldsymbol b\big)=B_r\sum_{(\boldsymbol a,\boldsymbol k)\in\mathbb Z^{n+c}}f_{\boldsymbol b,r}\big(\boldsymbol a,\boldsymbol k\big)+\text{error term}+O\Big(e^{\frac{r}{4\pi}\big(\mathrm{Vol}(M_\text{FSL})-\delta_0\big)}\Big),$$
where the error term comes from the use of the bump function. Since $f_{\boldsymbol b,r}$ is $C^{\infty}$-smooth and equals zero out of $\mathrm D,$ it is in the Schwartz space on $\mathbb R^{n+c}.$ Then by the Poisson Summation Formula (see e.g. \cite[Theorem 3.1]{SS}),
$$\sum_{(\boldsymbol a,\boldsymbol k)\in\mathbb Z^{n+c}}f_{\boldsymbol b,r}\big(\boldsymbol a,\boldsymbol k\big)=\sum_{(\boldsymbol m,\boldsymbol n)\in\mathbb Z^{n+c}}\widehat {f_{\boldsymbol b,r}}(\boldsymbol m,\boldsymbol n),$$
where $\boldsymbol m=(m_i)_{i\in\{1,\dots,n\}}\in \mathbb Z^n,$ $\boldsymbol n=(n_s)_{s\in\{1,\dots,c\}}\in\mathbb Z^c,$ and  $\widehat {f_{\boldsymbol b,r}}(\boldsymbol m,\boldsymbol n)$ is the $(\boldsymbol m,\boldsymbol n)$-th Fourier coefficient of $f_{\boldsymbol b,r}$ defined by
\begin{equation*}
\begin{split}
\widehat {f_{\boldsymbol b,r}}(\boldsymbol m,\boldsymbol n)=\int_{\mathbb R^{n+c}}&f_{\boldsymbol b,r}\big(\boldsymbol a,\boldsymbol k\big)e^{\sum_{i=1}^n2\pi \mathbf im_ia_i+\sum_{s=1}^c2\pi \mathbf in_sk_s}d\boldsymbol ad\boldsymbol k,
\end{split}
\end{equation*}
where $d\boldsymbol ad\boldsymbol k=\prod_{i=1}^nda_i\prod_{s=1}^cdk_s.$

\begin{proposition}\label{Poisson}
\begin{equation}\label{P1}
\mathrm{RT}_r\big(M,L,\boldsymbol b\big)=B_r\sum_{(\boldsymbol m,\boldsymbol n)\in\mathbb Z^{n+c}}\widehat{ f_{\boldsymbol b,r}}(\boldsymbol m,\boldsymbol n)+\text{error term}+O\Big(e^{\frac{r}{4\pi}\big(\mathrm{Vol}(M_\text{FSL})-\delta_0\big)}\Big),
\end{equation}
where in terms of the variables $(\boldsymbol \alpha,\boldsymbol \xi),$ the Fourier coefficients can be computed as
$$\widehat{f_{\boldsymbol b,r}}(\boldsymbol m,\boldsymbol n)=\sum_{\boldsymbol \epsilon\in\{-1,1\}^{n}} \widehat{f^{\boldsymbol \epsilon}_{\boldsymbol b,r}}(\boldsymbol m,\boldsymbol n)$$
with
\begin{equation*}
\widehat{f^{\boldsymbol \epsilon}_{\boldsymbol b,r}}(\boldsymbol m,\boldsymbol n)=\Big(\frac{r}{2\pi}\Big)^{n+c}\int_{\mathrm D}\psi(\boldsymbol{\alpha},\boldsymbol{\xi})e^{\frac{r}{4\pi \mathbf i}\big(\mathcal W_{\boldsymbol \beta,r}^{\boldsymbol \epsilon}(\boldsymbol{\alpha},\boldsymbol{\xi})-\sum_{i=1}^n4\pi m_i\alpha_i-\sum_{s=1}^c4\pi n_s\xi_s\big)}d\boldsymbol{\alpha}d\boldsymbol{\xi},
\end{equation*}
where $d\boldsymbol{\alpha}d\boldsymbol{\xi}=\prod_{i=1}^nd\alpha_i\prod_{s=1}^cd\xi_s,$ and $\mathcal W_{\boldsymbol \beta,r}^{\boldsymbol \epsilon}(\boldsymbol{\alpha},\boldsymbol{\xi})$ is as defined in Proposition \ref{rt4}.
\end{proposition}

We will estimate the leading Fourier coefficients, the non-leading Fourier coefficients and the error term respectively in Sections \ref{leading}, \ref{ot} and \ref{ee}, and prove Theorem \ref{re1} (hence Theorem \ref{main1}) in Section \ref{pf}.

\subsection{The closed  case}

In the setting of the previous section, 
the ambient manifold $M=(M_c)_{L_\text{FSL}}$ is the closed oriented $3$-manifold obtained from $M_c$ by doing surgery along the framed fundamental shadow link $L_\text{FSL}.$ As a consequence, 
$\mathrm{RT}_r(M)=\mathrm{RT}_r(M,L,\boldsymbol 0).$

Let 
$$\mathrm H_r=\Big\{ \boldsymbol a \in [0,r-2]^n\ \Big|\ \text{for each building block }\Delta_s, (a_{s_1},\dots,a_{s_6}) \text{ is $r$-hyperideal}  \Big\}$$
be as in the previous section. Then we have

\begin{proposition} \label{rt5}
\begin{equation*}
\mathrm{RT}_r(M)=B_r\sum_{\boldsymbol a,\boldsymbol k}\bigg(\sum_{\boldsymbol \epsilon}g^{\boldsymbol \epsilon}_r(\boldsymbol a,\boldsymbol k)\bigg)+O\Big(e^{\frac{r}{4\pi}\big(\mathrm{Vol}(M_\text{FSL})-\delta_0\big)}\Big),
\end{equation*}
where
$$B_r=2^{-n}r^{-\frac{n+c}{2}}e^{\big(-\sigma(U_{c+1}\cup L)(-\frac{3}{r}-\frac{r+1}{4})-(\frac{r}{4}+1+\frac{1}{r})\sum_{i=1}^np_i+(\frac{1}{2}-\frac{1}{2r})\sum_{i=1}^n\iota_i-\frac{r}{2}(n-c)-\frac{n}{2}\big)\pi\mathbf i},$$
$\boldsymbol \epsilon=(\epsilon_i)_{i\in\{1,\dots,n\}}\in\{-1,1\}^n$ runs over all multi-signs, $\boldsymbol a=(a_i)_{i\in \{1,\dots,n\}}$ runs over multi-integers in $\mathrm H_r,$ and $\boldsymbol k=(k_s)_{s\in \{1,\dots,c\}}$ runs over all multi-integers with each $k_s$ lying in between $\max\{T_{s_1},T_{s_2},T_{s_3},T_{s_4}\}$ and $\min\{Q_{s_1},Q_{s_2},Q_{s_3},Q_{s_4}\},$
with 
$$g^{\boldsymbol \epsilon}_{r}(\boldsymbol a,\boldsymbol k)=e^{\sum_{i=1}^n\epsilon_i(\alpha_i-\pi)\mathbf i}\cdot e^{\frac{r}{4\pi\mathbf i}\mathcal W_{r}^{\boldsymbol\epsilon}\big(\frac{2\pi (\boldsymbol a+\boldsymbol 1)}{r},\frac{2\pi \boldsymbol k}{r}\big)},$$
where $\frac{2\pi (\boldsymbol a+\boldsymbol 1)}{r}=\big(\frac{2\pi (a_i+1)}{r}\big)_{i\in \{1,\dots,n\}},$  $\frac{2\pi\boldsymbol k}{r}=\big(\frac{2\pi k_s}{r}\big)_{s\in \{1,\dots, c\}},$ and with  $\boldsymbol \alpha=(\alpha_i)_{i\in\{1,\dots,n\}},$ $\boldsymbol \xi=(\xi_s)_{s\in\{1,\dots,c\}},$ and
\begin{equation*}
\begin{split}
\mathcal W_r^{\boldsymbol \epsilon}(\boldsymbol{\alpha},\boldsymbol{\xi})=&-\sum_{i=1}^np_i(\alpha_i-\pi)^2+\sum_{i=1}^n2\epsilon_i\pi (\alpha_i-\pi)\\
&-\sum_{i=1}^n\frac{\iota_i}{2}(\alpha_i-\pi)^2+\sum_{s=1}^c U_r(\alpha_{s_1},\dots,\alpha_{s_6},\xi_s)+\Big(\sum_{i=1}^n\frac{\iota_i}{2}\Big)\pi^2.
\end{split}
\end{equation*}
\end{proposition}

\begin{proof}
Since the coloring  $\boldsymbol b=\boldsymbol 0,$ we have $\beta_i=\frac{2\pi}{r}$ for each $i\in\{1,\dots, n\}.$ Then the result follows directly from Proposition \ref{rt4}.
\end{proof}

To use the Poisson Summation Formula, we let 
$$\mathrm D=\Big\{ (\boldsymbol \alpha,\boldsymbol \xi) \ \Big|\ (\alpha_{s_1},\dots,\alpha_{s_6},\xi_s)\in\mathrm B\text{ for each }s\in\{1,\dots,c\}\Big\},$$
for a sufficiently small $\epsilon >0$ let 
$$\mathrm D_{\epsilon}=\Big\{(\boldsymbol{\alpha},\boldsymbol{\xi})\in\mathrm {D}\ \Big|\ d\big((\boldsymbol{\alpha},\boldsymbol{\xi}), \partial\mathrm {D}\big)>\epsilon \Big\},$$
and let 
 $\psi:\mathbb R^{n+c}\to\mathbb R$ be a $C^{\infty}$-smooth bump function supported on $(\mathrm{D}, \mathrm D_\epsilon)$ as in the previous section.
Let 
$$f^{\boldsymbol \epsilon}_r(\boldsymbol a,\boldsymbol k)=\psi\Big(\frac{2\pi (\boldsymbol a+\boldsymbol 1)}{r},\frac{2\pi \boldsymbol k}{r}\Big)g^{\boldsymbol \epsilon}_r(\boldsymbol a,\boldsymbol k),$$
and let $f_r=\sum_{\boldsymbol \epsilon\in\{-1,1\}^n} f_r^{\boldsymbol \epsilon}.$ Then
$$\mathrm{RT}_r(M)=B_r\sum_{(\boldsymbol a,\boldsymbol k)\in\mathbb Z^{n+c}}f_r\big(\boldsymbol a,\boldsymbol k\big)+\text{error term}+O\Big(e^{\frac{r}{4\pi}\big(\mathrm{Vol}(M_\text{FSL})-\delta_0\big)}\Big),$$
where the error term comes from the use of the bump function; and by using the Poisson Summation Formula as in Proposition \ref{Poisson}, we have 
\begin{proposition}\label{Poisson2}
\begin{equation}\label{P2}
    \mathrm{RT}_r(M)=B_r\sum_{(\boldsymbol m,\boldsymbol n)\in\mathbb Z^{n+c}}\widehat{ f_r}(\boldsymbol m,\boldsymbol n)+\text{error term}+O\Big(e^{\frac{r}{4\pi}\big(\mathrm{Vol}(M_\text{FSL})-\delta_0\big)}\Big),
\end{equation}
where in terms of the variables $(\boldsymbol \alpha,\boldsymbol \xi),$ the Fourier coefficients can be computed as
$$\widehat{f_r}(\boldsymbol m,\boldsymbol n)=\sum_{\boldsymbol \epsilon\in\{-1,1\}^n} \widehat{f^{\boldsymbol \epsilon}_r}(\boldsymbol m,\boldsymbol n)$$
with
\begin{equation*}
\widehat{f^{\boldsymbol \epsilon}_r}(\boldsymbol m,\boldsymbol n)=\Big(\frac{r}{2\pi}\Big)^{n+c}\int_{\mathrm{D}}\psi(\boldsymbol{\alpha},\boldsymbol{\xi})e^{\sum_{i=1}^n\epsilon_i(\alpha_i-\pi)\mathbf i}e^{\frac{r}{4\pi \mathbf i}\big(\mathcal W_r^{\boldsymbol \epsilon}(\boldsymbol{\alpha},\boldsymbol{\xi})-\sum_{i=1}^n4\pi m_i\alpha_i-\sum_{s=1}^c4\pi n_s\xi_s\big)}d\boldsymbol{\alpha}_Id\boldsymbol{\xi},
\end{equation*}
where $d\boldsymbol{\alpha}d\boldsymbol{\xi}=\prod_{i=1}^nd\alpha_i\prod_{s=1}^cd\xi_s,$ and $\mathcal W_r^{\boldsymbol \epsilon}(\boldsymbol{\alpha},\boldsymbol{\xi})$ is as defined in Proposition \ref{rt5}.
\end{proposition}

We will estimate the leading Fourier coefficients, the non-leading Fourier coefficients and the error term respectively in Sections \ref{leading}, \ref{ot} and \ref{ee}, and prove Theorem \ref{main2} in Section \ref{pf}.


\section{Asymptotics}\label{pfMain}

The goal of this section is to prove Theorem \ref{main1} (by proving Theorem \ref{re1}) and Theorem \ref{main2}. The argument follows almost verbatim that in \cite{WY1,LMSWY}.  We will point out the main difference and give more details there.  The main tool we use is still the following  Saddle Point Approximations.

\begin{proposition}\cite[Proposition 5.1]{WY1}\label{saddle}
Let $D_{\boldsymbol z}$ be a region in $\mathbb C^n$ and let $D_{\boldsymbol a}$ be a region in $\mathbb R^k.$ Let $f(\boldsymbol z,\boldsymbol a)$ and $g(\boldsymbol z,\boldsymbol a)$ be complex valued functions on $D_{\boldsymbol z}\times D_{\boldsymbol a}$  which are holomorphic in $\boldsymbol z$ and smooth in $\boldsymbol a.$ For each positive integer $r,$ let $f_r(\boldsymbol z,\boldsymbol a)$ be a complex valued function on $D_{\boldsymbol z}\times D_{\boldsymbol a}$ holomorphic in $\boldsymbol z$ and smooth in $\boldsymbol a.$
For a fixed $\boldsymbol a\in D_{\boldsymbol a},$ let $f^{\boldsymbol a},$ $g^{\boldsymbol a}$ and $f_r^{\boldsymbol a}$ be the holomorphic functions  on $D_{\boldsymbol z}$ defined by
$f^{\boldsymbol a}(\boldsymbol z)=f(\boldsymbol z,\boldsymbol a),$ $g^{\boldsymbol a}(\boldsymbol z)=g(\boldsymbol z,\boldsymbol a)$ and $f_r^{\boldsymbol a}(\boldsymbol z)=f_r(\boldsymbol z,\boldsymbol a).$ Suppose $\{\boldsymbol a_r\}$ is a convergent sequence in $D_{\boldsymbol a}$ with $\lim_r\boldsymbol a_r=\boldsymbol a_0,$ $f_r^{\boldsymbol a_r}$ is of the form
$$ f_r^{\boldsymbol a_r}(\boldsymbol z) = f^{\boldsymbol a_r}(\boldsymbol z) + \frac{\upsilon_r(\boldsymbol z,\boldsymbol a_r)}{r^2},$$
$\{S_r\}$ is a sequence of embedded real $n$-dimensional closed disks in $D_{\boldsymbol z}$ sharing the same boundary and converging to an embedded $n$-dimensional disk $S_0,$ and $\boldsymbol c_r$ is a point on $S_r$ such that $\{\boldsymbol c_r\}$ is convergent  in $D_{\boldsymbol z}$ with $\lim_r\boldsymbol c_r=\boldsymbol c_0.$ If for each $r$
\begin{enumerate}[(i)]
\item $\boldsymbol c_r$ is a critical point of $f^{\boldsymbol a_r}$ in $D_{\boldsymbol z},$
\item $\mathrm{Re}f^{\boldsymbol a_r}(\boldsymbol c_r) > \mathrm{Re}f^{\boldsymbol a_r}(\boldsymbol z)$ for all $\boldsymbol z \in S_r\setminus \{\boldsymbol c_r\},$
\item the domain $\{\boldsymbol z\in D_{\boldsymbol z}\ |\ \mathrm{Re} f^{\boldsymbol a_r}(\boldsymbol z) < \mathrm{Re} f^{\boldsymbol a_r}(\boldsymbol c_r)\}$ deformation retracts to $S_r\setminus\{\boldsymbol c_r\},$
\item  the Hessian matrix $\mathrm{Hess}(f^{\boldsymbol a_0})$ of $f^{\boldsymbol a_0}$ at $\boldsymbol c_0$ is non-singular,
\item $|g^{\boldsymbol a_r}(\boldsymbol c_r)|$ is bounded from below by a positive constant independent of $r,$ and
\item $|\upsilon_r(\boldsymbol z, \boldsymbol a_r)|$ is bounded from above by a constant independent of $r$ on $D_{\boldsymbol z},$ 
\end{enumerate}
then
\begin{equation*}
\begin{split}
 \int_{S_r} g^{\boldsymbol a_r}(\boldsymbol z) e^{rf_r^{\boldsymbol a_r}(\boldsymbol z)} d\boldsymbol z= \Big(\frac{2\pi}{r}\Big)^{\frac{n}{2}}\frac{g^{\boldsymbol a_r}(\boldsymbol c_r)}{\sqrt{\det\big(-\mathrm{Hess}(f^{\boldsymbol a_r})(\boldsymbol c_r)\big)}} e^{rf^{\boldsymbol a_r}(\boldsymbol c_r)} \Big( 1 + O \Big( \frac{1}{r} \Big) \Big).
 \end{split}
 \end{equation*}
\end{proposition}

\begin{proposition}\cite[Proposition 6.1]{WY}\label{saddle2}
Let $D$ be a region in $\mathbb C^n$ and let $f(\boldsymbol z)$ and $g(\boldsymbol z)$ be holomorphic functions on $D$ independent of $r.$ Let $f_r(\boldsymbol z)$ be a holomorphic function of the form
$$ f_r(\boldsymbol z) = f(\boldsymbol z) + \frac{\upsilon_r(\boldsymbol z)}{r^2}.$$
Let $S$ be an embedded real $n$-dimensional closed disk in $D$ and let $\boldsymbol c$ be a point on $S.$ If
\begin{enumerate}[(i)]
\item $\boldsymbol c$ is a critical point of $f$ in $D,$
\item $\mathrm{Re}(f)(\boldsymbol c) > \mathrm{Re}(f)(\boldsymbol z)$ for all $\boldsymbol z \in S\setminus \{\boldsymbol c\},$
\item the domain $\{\boldsymbol z\in D\ |\ \mathrm{Re}{f}(\boldsymbol z)<\mathrm{Re}{f}(\boldsymbol c)\}$ deformation retracts to $S\setminus\{(\boldsymbol c)\},$
\item the Hessian matrix $\mathrm{Hess}(f)(\boldsymbol c)$ of $f$ at $\boldsymbol c$ is non-singular,
\item $g(\boldsymbol c) \neq 0,$ and
\item $|\upsilon_r(\boldsymbol z)|$ is bounded from above by a constant independent of $r$ in $D,$
\end{enumerate}
then
\begin{equation*}
\begin{split}
 \int_S g(\boldsymbol z) &e^{rf_r(\boldsymbol z)} d\boldsymbol z= \Big(\frac{2\pi}{r}\Big)^{\frac{n}{2}}\frac{g(\boldsymbol c)}{\sqrt{\det\big(-\mathrm{Hess}(f)(\boldsymbol c)\big)}} e^{rf(\boldsymbol c)} \Big( 1 + O \Big( \frac{1}{r} \Big) \Big).
 \end{split}
 \end{equation*}
\end{proposition}


\subsection{Properties of relevant functions}

In light of Proposition \ref{EST}, we consider the functions $\mathcal W_{\boldsymbol \beta}^{\boldsymbol \epsilon}$ and $\mathcal W^{\boldsymbol \epsilon}$ on $\mathrm D,$ respectively defined by
\begin{equation*}
\begin{split}
\mathcal W_{\boldsymbol \beta}^{\boldsymbol \epsilon}(\boldsymbol{\alpha},\boldsymbol{\xi})=&-\sum_{i=1}^nq_i(\beta_i-\pi)^2-\sum_{i=1}^np_i(\alpha_i-\pi)^2-\sum_{i=1}^n2\epsilon_i(\alpha_i-\pi)(\beta_i-\pi)\\
&-\sum_{i=1}^n\frac{\iota_i}{2}(\alpha_i-\pi)^2+\sum_{s=1}^c U(\alpha_{s_1},\dots,\alpha_{s_6},\xi_s)+\Big(\sum_{i=1}^n\frac{\iota_i}{2}\Big)\pi^2
\end{split}
\end{equation*}
and
\begin{equation*}
\begin{split}
\mathcal W^{\boldsymbol \epsilon}(\boldsymbol{\alpha},\boldsymbol{\xi})=&-\sum_{i=1}^np_i(\alpha_i-\pi)^2+\sum_{i=1}^n2\epsilon_i\pi (\alpha_i-\pi)\\
&-\sum_{i=1}^n\frac{\iota_i}{2}(\alpha_i-\pi)^2+\sum_{s=1}^c U(\alpha_{s_1},\dots,\alpha_{s_6},\xi_s)+\Big(\sum_{i=1}^n\frac{\iota_i}{2}\Big)\pi^2,
\end{split}
\end{equation*}
where the function $U$ is as defined in (\ref{U}); and by abuse of notations consider the functions $\kappa$ and $\nu_b$ on $D,$ respectively defined by
$$\kappa(\boldsymbol \alpha,\boldsymbol \xi)=\sum_{s=1}^c\kappa(\alpha_{s_1},\dots,\alpha_{s_6},\xi_s)$$
and
$$\nu_r(\boldsymbol \alpha,\boldsymbol \xi)=\sum_{s=1}^c\nu_r(\alpha_{s_1},\dots,\alpha_{s_6},\xi_s),$$
where $\kappa$ and $\nu_r$ on the right hand side are as defined in  (\ref{kappa}) and (\ref{nu}). Notice that all the functions above  extend analytically to 
$$\mathrm D_{\mathbb C}=\Big\{(\boldsymbol{\alpha},\boldsymbol{\xi})\in\mathbb C^{n+c}\ \Big|\ (\mathrm{Re}\boldsymbol{\alpha},\mathrm{Re}\boldsymbol\xi)\in \mathrm{D}\Big\},$$ 
where $\mathrm{Re}\boldsymbol{\alpha}=(\mathrm{Re}\alpha_i)_{i\in \{1,\dots,n\}}$ and $\mathrm{Re}\boldsymbol\xi=(\mathrm{Re}\xi_s)_{s\in\{1,\dots,c\}}.$

For $\epsilon>0,$ we denote by $B_\epsilon(\boldsymbol z_0)$ the $\epsilon$-neighborhood of $\boldsymbol z_0=\big(\pi,\dots,\pi,\frac{7\pi}{4},\dots, \frac{7\pi}{4}\big)$ in $\mathbb C^{n+c}$ in the $L^\infty$ distance,  that is 
$$B_\epsilon(\boldsymbol z_0)=\Big\{(\boldsymbol{\alpha},\boldsymbol{\xi})\in \mathbb C^{n+c}\ \Big|\ d_{L^\infty}\big((\boldsymbol{\alpha},\boldsymbol{\xi}),\boldsymbol z_0\big)<\epsilon\Big\},$$
where $d_{L^\infty}$ is the $L^\infty$-distance on $\mathbb C^{n+c}.$

\begin{proposition}\label{convexity} 
There exists an $\epsilon'>0$  independent of $M_c,$ $L_\text{FSL}$ and $L$ such that for any $\boldsymbol \beta\in [0,2\pi]^n$ and $\boldsymbol \epsilon\in\{-1,1\}^n,$ $\mathrm{Im}\mathcal W_{\boldsymbol \beta}^{\epsilon}$ is strictly concave down in $\{\mathrm{Re}(\alpha_i)\}_{i\in\{1,\dots,n\}}$ and in $\{\mathrm{Re}(\xi_s)\}_{s\in\{1,\dots,c\}},$ and is strictly concave up in $\{\mathrm{Im}(\alpha_i)\}_{i\in \{1,\dots,n\}}$ and in  $\{\mathrm{Im}(\xi_s)\}_{s\in\{1,\dots,c\}}$  on $B_{\epsilon'}(\boldsymbol z_0).$
\end{proposition}

\begin{proof} In the special case $(\alpha_1,\dots,\alpha_6,\xi)=(\pi,\dots,\pi,\frac{7\pi}{4}),$ by a direct computation, the Hessian matrix of  $\mathrm{Im}U$ is strictly negative definite in $\mathrm{Re}\alpha_i$'s and in $\mathrm{Re}\xi.$ Then by continuity, there exists an $\epsilon'>0$ such that if $d_{L^\infty}\big((\alpha_1,\dots,\alpha_6,\xi),(\pi,\dots,\pi, \frac{7\pi}{4})\big)<\epsilon',$ the Hessian matrix of $\mathrm{Im}U$ is strictly negative definite, i.e., $\mathrm{Im}U$ is strictly concave down in $\mathrm{Re}\alpha_i$'s and in $\mathrm{Re}\xi.$ This implies that $\sum_{s=1}^c\mathrm{Im}U(\boldsymbol \alpha_s,\xi_s)$ is strictly  concave down in $\{\mathrm{Re}(\alpha_i)\}_{i\in\{1,\dots,n\}}$ and in $\{\mathrm{Re}(\xi_s)\}_{s\in\{1,\dots,c\}}$ on $B_{\epsilon'}(\boldsymbol z_0).$ As all the other terms in $\mathrm{Im}\mathcal W_{\boldsymbol \beta}^{\boldsymbol \epsilon}$ are linear functions of $\{\mathrm{Re}(\alpha_i)\}_{i\in\{1,\dots,n\}}$ and $\{\mathrm{Re}(\xi_s)\}_{s\in\{1,\dots,c\}},$ $\mathrm{Im}\mathcal W_{\boldsymbol \beta}^{\boldsymbol \epsilon}$  is  strictly  concave down in $\{\mathrm{Re}(\alpha_i)\}_{i\in\{1,\dots,n\}}$ and in $\{\mathrm{Re}(\xi_s)\}_{s\in\{1,\dots,c\}}$ on $B_{\epsilon'}(\boldsymbol z_0);$ and by the harmonicity, it is strictly concave up in $\{\mathrm{Im}(\alpha_i)\}_{i\in \{1,\dots,n\}}$ and in  $\{\mathrm{Im}(\xi_s)\}_{s\in\{1,\dots,c\}}$  on $B_{\epsilon'}(\boldsymbol z_0).$
\end{proof}

For $\boldsymbol \beta\in (0,2\pi)^n,$ $\boldsymbol \epsilon = (\epsilon_i)_{i\in\{1,\dots,n\}}\in \{-1,1\}^n$ and  $i\in \{1,\dots, n\},$ let 
\begin{equation}\label{alpha1}
\alpha^{\boldsymbol \epsilon}_i=\pi+\frac{\epsilon_i\mu_i}{2} \mathbf i \mathrm H(u_i)
\end{equation}
with $\mathrm H(u_i)$ the logarithmic holonomy of $u_i$ of the hyperbolic cone manifold $M_{L_{\boldsymbol\theta}},$
and let 
\begin{equation}\label{alpha2}
    \boldsymbol \alpha^{\boldsymbol \epsilon}=\big(\alpha^{\boldsymbol \epsilon}_i\big)_{i\in \{1,\dots,n\}}.
\end{equation}
For $s\in\{1,\dots,c\},$ 
let $\xi^{\boldsymbol \epsilon}_s=\xi_{(\alpha^{\boldsymbol \epsilon}_{s_1},\dots,\alpha^{\boldsymbol \epsilon}_{s_6})}$ be as defined in (\ref{xia}), and let $$\boldsymbol \xi^{\boldsymbol \epsilon}=\big(\xi^{\boldsymbol \epsilon}_s\big)_{s\in\{1,\dots,c\}}.$$
Keep in mind that both $\boldsymbol \alpha^{\boldsymbol \epsilon}$ and $\boldsymbol \xi^{\boldsymbol \epsilon}$ depend on $\boldsymbol \beta.$

Also, for $\boldsymbol \sigma\in \{0,1\}^n,$ $\boldsymbol \epsilon = (\epsilon_i)_{i\in\{1,\dots,n\}}\in \{-1,1\}^n$ and  $i\in \{1,\dots, n\},$ let 
\begin{equation}\label{alpha3}
\alpha^{\boldsymbol
\sigma,\boldsymbol \epsilon}_i=\pi-\frac{(-1)^{\sigma_i}\epsilon_i}{2} \mathbf i \mathrm H(u_i)
\end{equation}
with $\mathrm H(u_i)$ the logarithmic holonomy of $u_i$ of the hyperbolic metric on $M,$
and let 
\begin{equation}\label{alpha4}
    \boldsymbol \alpha^{\boldsymbol\sigma,\boldsymbol \epsilon}=\big(\alpha^{\boldsymbol\sigma,\boldsymbol \epsilon}_i\big)_{i\in \{1,\dots,n\}}.
\end{equation}
For $s\in\{1,\dots,c\},$ 
let $\xi^{\boldsymbol\sigma,\boldsymbol \epsilon}_s=\xi_{(\alpha^{\boldsymbol\sigma,\boldsymbol \epsilon}_{s_1},\dots,\alpha^{\boldsymbol\sigma,\boldsymbol \epsilon}_{s_6})}$ be as defined in (\ref{xia}), and let $$\boldsymbol \xi^{\boldsymbol\sigma,\boldsymbol \epsilon}=\big(\xi^{\boldsymbol\sigma,\boldsymbol \epsilon}_s\big)_{s\in\{1,\dots,c\}}.$$

\begin{proposition}\label{crit} There exists an $\epsilon''>0$ independent of $M_c,$ $L_\text{FSL}$ and $L$ and depending only on $\epsilon'$ in Proposition \ref{convexity} such that:
\begin{enumerate}[(1)]
    \item If the $L^\infty$-distance 
 $$d_{L^\infty}(\boldsymbol \alpha^{\boldsymbol \epsilon},\boldsymbol \pi)<\epsilon'',$$ then
$$\boldsymbol z^{\boldsymbol \epsilon}=\big(\boldsymbol \alpha^{\boldsymbol \epsilon},\boldsymbol \xi^{\boldsymbol \epsilon}\big)\in B_{\epsilon'}(\boldsymbol z_0),$$
which is a critical point of $\mathcal W_{\boldsymbol \beta}^{\boldsymbol \epsilon}$ 
with critical value $$\mathcal W_{\boldsymbol \beta}^{\boldsymbol \epsilon}(\boldsymbol z^{\boldsymbol \epsilon})=\mathbf i\Big(\mathrm{Vol}(M_{L_{\boldsymbol\theta}})+\mathbf i\mathrm{CS}(M_{L_{\boldsymbol\theta}})\Big)$$
for a suitable choice of the value of the Chern-Simons invariant $\mathrm{CS}(M_{L_{\boldsymbol\theta}}).$

\item If the $L^\infty$-distance 
 $$d_{L^\infty}(\boldsymbol \alpha^{\boldsymbol\sigma,\boldsymbol \epsilon},\boldsymbol \pi)<\epsilon'',$$ then
$$\boldsymbol z^{\boldsymbol\sigma,\boldsymbol \epsilon}=\big(\boldsymbol \alpha^{\boldsymbol\sigma,\boldsymbol \epsilon},\boldsymbol \xi^{\boldsymbol\sigma,\boldsymbol \epsilon}\big)\in B_{\epsilon'}(\boldsymbol z_0),$$
which is a critical point of $\mathcal W^{\boldsymbol \epsilon}-4\pi\sum_{i=1}^n\sigma_i\epsilon_i\alpha_i$ 
with critical value $$\mathcal W^{\boldsymbol \epsilon}(\boldsymbol z^{\boldsymbol\sigma,\boldsymbol\epsilon})-4\pi\sum_{i=1}^n\sigma_i\epsilon_i\alpha^{\boldsymbol\sigma,\boldsymbol\epsilon}_i=\mathbf i\Big(\mathrm{Vol}(M)+\mathbf i\mathrm{CS}(M)\Big) -4\pi^2\sum_{i=1}^n\sigma_i\epsilon_i$$
for a suitable choice of the value of the Chern-Simons invariant $\mathrm{CS}(M).$
\end{enumerate}
\end{proposition}

\begin{proof}
In the special case that  $(\alpha_1,\dots,\alpha_6)=(\pi,\dots,\pi),$ 
by solving equation (\ref{xia}) directly,  we have $\xi_{(\pi,\dots,\pi)}=\frac{7\pi}{4},$ and $(\pi,\dots,\pi,\frac{7\pi}{4})\in B_{\epsilon'}(\boldsymbol z_0).$ Then by continuity, there exists an $\epsilon''>0$ such that if $d_{L^\infty}\big((\alpha_1,\dots,\alpha_6),(\pi,\dots,\pi)\big)<\epsilon'',$ then $d_{L^\infty}\big((\alpha_1,\dots,\alpha_6,\xi_{(\alpha_1,\dots,\alpha_6)}),(\pi,\dots,\pi,\frac{7\pi}{4})\big)<\epsilon'.$ As a consequence, if $d_{L^\infty}(\boldsymbol \alpha^{\boldsymbol \epsilon}, \boldsymbol \pi)<\epsilon'',$ then $(\boldsymbol \alpha^{\boldsymbol \epsilon}, \boldsymbol \xi^{\boldsymbol \epsilon})\in B_{\epsilon'}(\boldsymbol z_0)$ as desired in (1); and if $d_{L^\infty}(\boldsymbol \alpha^{\boldsymbol\sigma,\boldsymbol \epsilon}, \boldsymbol \pi)<\epsilon'',$ then $(\boldsymbol \alpha^{\boldsymbol\sigma,\boldsymbol \epsilon}, \boldsymbol \xi^{\boldsymbol\sigma,\boldsymbol \epsilon})\in B_{\epsilon'}(\boldsymbol z_0)$ as desired in (2).

The rest of the proof of (1) follows verbatim that of \cite[Propositions 4.1 and 5.2]{WY1}; and we include a sketch of the proof of (2), which is by essentially the same argument, for the readers' convenience. As proved in \cite[Proposition 4.1]{WY1}, the function 
$$\mathcal U(\boldsymbol\alpha)\doteq-\sum_{i=1}^n\frac{\iota_i}{2}(\alpha_i-\pi)^2+\sum_{s=1}^c U_r(\alpha_{s_1},\dots,\alpha_{s_6},\xi_{(\alpha_{s_1},\dots,\alpha_{s_6})})+\Big(\sum_{i=1}^n\frac{\iota_i}{2}\Big)\pi^2$$
coincides with the Neumann-Zagier potential function $\Phi$ recalled in Section \ref{CCS}; and by (\ref{char}), we have   
\begin{equation}\label{char2}
\frac{\partial \mathcal U(\boldsymbol \alpha)}{\partial \mathrm{H}(u_i)}=\frac{\mathrm H(v_i)}{2}.\\
\end{equation} 
Then by the chain rule, one verifies that for each $s\in\{1,\dots,c\},$
$$\frac{\partial}{\partial \xi_s}\Big|_{\boldsymbol z^{\boldsymbol\sigma,\boldsymbol\epsilon}}\Big(\mathcal W^{\boldsymbol \epsilon}(\boldsymbol\alpha,\boldsymbol\xi)-4\pi\sum_{i=1}^n\sigma_i\epsilon_i\alpha_i\Big)=\frac{\partial}{\partial \xi_s}\Big|_{\xi_s^{\boldsymbol\sigma,\boldsymbol\epsilon}} U(\alpha_{s_1}^{\boldsymbol\sigma,\boldsymbol\epsilon},\dots,\alpha_{s_6}^{\boldsymbol\sigma,\boldsymbol\epsilon},\xi_s)=0;$$
and by (\ref{char2}) and the chain rule, for each $i\in\{1,\dots,n\},$
\begin{equation*}
\begin{split}
    \frac{\partial}{\partial \alpha_i}\Big|_{\boldsymbol z^{\boldsymbol\sigma,\boldsymbol\epsilon}}\Big(\mathcal W^{\boldsymbol \epsilon}(\boldsymbol\alpha,\boldsymbol\xi)-4\pi\sum_{i=1}^n\sigma_i\epsilon_i\alpha_i\Big)=& (-1)^{\sigma_i}\epsilon_i\mathbf i \Big(p_i\mathrm H(u_i)+\mathrm H(v_i)-(-1)^{\sigma_i}\big(2-4\sigma_i\big)\pi\mathbf i\Big)\\
=& (-1)^{\sigma_i}\epsilon_i\mathbf i \Big(p_i\mathrm H(u_i)+\mathrm H(v_i)-2\pi\mathbf i\Big)=0,
\end{split}    
\end{equation*}
where the last equality comes from Thurston's Hyperbolic Gluing Equations. Hence $\boldsymbol z^{\boldsymbol\sigma,\boldsymbol \epsilon}$ is a critical point of $\mathcal W^{\boldsymbol \epsilon}-4\pi\sum_{i=1}^n\sigma_i\epsilon_i\alpha_i.$ 

For each $i\in I,$ let $\gamma_i=-u_i$ be the curve on the boundary of a tubular neighborhood of $L^*_i$ that is isotopic to the core of $L^*_i$ and with the orientation so that $(p_iu_i+v_i)\cdot \gamma_i=1.$ 
Then 
$$\mathrm H(v_i)=2\pi\mathbf i - p_i\mathrm H(u_i)\quad\text{and}\quad\mathrm H(\gamma_i)=-\mathrm H(u_i);$$
and plugging in (\ref{alpha3}), we have
\begin{equation*}
-\frac{\mathrm H(u_i)\mathrm H(v_i)}{4}+\frac{2\pi\mathbf i \mathrm H(\gamma_i)}{4}
    =- p_i(\alpha^{\boldsymbol\sigma,\boldsymbol\epsilon}_i-\pi)^2+ 2(-1)^{\sigma_1}\epsilon_i\pi(\alpha^{\boldsymbol\sigma,\boldsymbol\epsilon}_i-\pi).
\end{equation*}
If $\sigma_1=0,$ then we have 
$$- p_i(\alpha^{\boldsymbol\sigma,\boldsymbol\epsilon}_i-\pi)^2+ 2\epsilon_i\pi(\alpha^{\boldsymbol\sigma,\boldsymbol\epsilon}_i-\pi)-4\pi\sigma_i\epsilon_i\alpha_i^{\boldsymbol\sigma,\boldsymbol\epsilon}=-\frac{\mathrm H(u_i)\mathrm H(v_i)}{4}+\frac{2\pi\mathbf i \mathrm H(\gamma_i)}{4};$$
and if  $\sigma_i=1,$ then we have     
$$- p_i(\alpha^{\boldsymbol\sigma,\boldsymbol\epsilon}_i-\pi)^2+ 2\epsilon_i\pi(\alpha^{\boldsymbol\sigma,\boldsymbol\epsilon}_i-\pi)-4\pi\sigma_i\epsilon_i\alpha_i^{\boldsymbol\sigma,\boldsymbol\epsilon}=-\frac{\mathrm H(u_i)\mathrm H(v_i)}{4}+\frac{2\pi\mathbf i \mathrm H(\gamma_i)}{4}-4\epsilon_i\pi^2.$$
Then by (\ref{Yo}), we have 
$$\mathcal W^{\boldsymbol \epsilon}(\boldsymbol z^{\boldsymbol\sigma,\boldsymbol\epsilon})-4\pi\sum_{i=1}^n\sigma_i\epsilon_i\alpha^{\boldsymbol\sigma,\boldsymbol\epsilon}_i=\mathbf i\Big(\mathrm{Vol}(M)+\mathbf i\mathrm{CS}(M)\Big) \mod \pi^2\mathbb Z;$$
and with the choice of the value of the Chern-Simons invariant so that 
$$\mathcal W^{\boldsymbol \epsilon}(\boldsymbol z^{\boldsymbol 0,\boldsymbol\epsilon})=\mathbf i\Big(\mathrm{Vol}(M)+\mathbf i\mathrm{CS}(M)\Big),$$
we have 
$$\mathcal W^{\boldsymbol \epsilon}(\boldsymbol z^{\boldsymbol\sigma,\boldsymbol\epsilon})-4\pi\sum_{i=1}^n\sigma_i\epsilon_i\alpha^{\boldsymbol\sigma,\boldsymbol\epsilon}_i = \mathbf i\Big(\mathrm{Vol}(M)+\mathbf i\mathrm{CS}(M)\Big)-4\pi^2\sum_{i=1}^n\sigma_i\epsilon_i.$$
This completes the proof of (2).
\end{proof}

\begin{proposition}\label{Hessian}
\begin{enumerate}[(1)]
    \item For each $\boldsymbol\epsilon\in\{-1,1\}^n,$ the determinant of the Hessian matrix of $\mathcal W_{\boldsymbol \beta}^{\boldsymbol \epsilon}$ at the critical point $\boldsymbol z^{\epsilon}$ equals
$$\frac{\det\big(-\mathrm{Hess}\mathcal W_{\boldsymbol \beta}^{\boldsymbol\epsilon}(\boldsymbol z^{\boldsymbol \epsilon})\big)}{\exp\big({\frac{\kappa(\boldsymbol z^{\boldsymbol \epsilon})}{2\pi \mathbf i}}\big)}= \pm{2^{n+c}}\cdot \mathbb{T}_{(M\setminus L,\mathbf m)}([\rho_{M^{(r)}}]).$$
As a consequence, $\boldsymbol z^{\boldsymbol \epsilon}$ is a nondegenerate critical point of $\mathcal W_{\boldsymbol \beta}^{\boldsymbol\epsilon}.$

   \item For each $\boldsymbol\sigma\in\{0,1\}^n$ and $\boldsymbol\epsilon\in\{-1,1\}^n,$ the determinant of the Hessian matrix of $\mathcal W^{\boldsymbol \epsilon}$ at the critical point $\boldsymbol z^{\boldsymbol\sigma,\boldsymbol\epsilon}$ of $\mathcal W^{\boldsymbol \epsilon}-4\pi\sum_{i=1}^n\sigma_i\epsilon_i\alpha_i$ equals
$$\frac{\det\big(-\mathrm{Hess}\mathcal W^{\boldsymbol \epsilon}(\boldsymbol z^{\boldsymbol\sigma,\boldsymbol \epsilon})\big)}{\exp\big({\frac{\kappa(\boldsymbol z^{\boldsymbol\sigma, \boldsymbol \epsilon})}{2\pi \mathbf i}}\big)}= \pm{2^{n+c}}\cdot \mathbb{T}_{(M\setminus L,\mathbf m)}([\rho_{M\setminus L}]),$$
where $\rho_{M\setminus L}$ is the holonomy representation of the restriction of the hyperbolic metric on $M$ to $M\setminus L.$
As a consequence, $\boldsymbol z^{\boldsymbol \sigma,\boldsymbol \epsilon}$ is a nondegenerate critical point of $\mathcal W^{\boldsymbol \epsilon}-4\pi\sum_{i=1}^n\sigma_i\epsilon_i\alpha_i$ 
\end{enumerate}
\end{proposition}

\begin{proof}
The proof follows from the same algebraic computation as in the proof of \cite[Formula (5.11)]{LMSWY}.
\end{proof}


\subsection{Asymptotics of the leading Fourier coefficients}\label{leading}

\begin{proposition}\label{critical}
Let $\epsilon'$ and $\epsilon''$ respectively be the constants in Proposition  \ref{convexity} and Proposition \ref{crit}, and let 
$$B_{\epsilon',\epsilon''}(\boldsymbol z_0)=B_{\min\{\epsilon',\epsilon''\}}(\boldsymbol \pi)\times B_{\epsilon'}\Big(\frac{7\boldsymbol\pi}{4}\Big).$$
\begin{enumerate}[(1)]
\item  Let $\boldsymbol \alpha^{\boldsymbol \epsilon}$ be as defined in (\ref{alpha1}) and (\ref{alpha2}).  If $d_{L^\infty}(\boldsymbol\alpha^{\boldsymbol \epsilon},\boldsymbol\pi)<\min\{\epsilon',\epsilon''\}$ and 
$$\mathrm{Vol}(M_{L_{\boldsymbol \theta}})>\max_{\mathrm{D}\setminus  B_{\epsilon',\epsilon''}(\boldsymbol z_0)} \mathrm{Im}\mathcal W_{\boldsymbol \beta}^{\boldsymbol \epsilon},$$
then
$$\widehat{f^{\boldsymbol \epsilon}_{\boldsymbol b,r}}(\boldsymbol 0,\boldsymbol 0)= \frac{e^{\frac{(n+c)\pi\mathbf i }{4}}r^{\frac{n+c}{2}}}{\sqrt{\pm\mathbb{T}_{(M\setminus L,\mathbf m)}([\rho_{M_{L_{\boldsymbol \theta}}}])}}e^{\frac{r}{4\pi}\big(\mathrm{Vol}(M_{L_{\boldsymbol\theta}})+\mathbf i\mathrm{CS}(M_{L_{\boldsymbol\theta}})\big)}\Big( 1 + O \Big( \frac{1}{r} \Big) \Big).$$

\item Let $\boldsymbol \alpha^{\boldsymbol\sigma,\boldsymbol \epsilon}$ be as defined in (\ref{alpha3}) and (\ref{alpha4}), and let $\boldsymbol \sigma\boldsymbol\epsilon\doteq(\sigma_i\epsilon_i)_{i\in\{1,\dots,n\}}.$   If $d_{L^\infty}(\boldsymbol\alpha^{\boldsymbol\sigma,\boldsymbol \epsilon},\boldsymbol\pi)<\min\{\epsilon',\epsilon''\}$ and 
\begin{equation}\label{assump}
\mathrm{Vol}(M)
>\max_{\mathrm{D}\setminus B_{\epsilon',\epsilon''}(\boldsymbol z_0)} \mathrm{Im}\mathcal W^{\boldsymbol \epsilon},
\end{equation}
then
$$\widehat{f^{\boldsymbol \epsilon}_r}(\boldsymbol\sigma\boldsymbol\epsilon,\boldsymbol 0)=e^{\frac{(n+c)\pi\mathbf i }{4}}r^{\frac{n+c}{2}}\frac{(-1)^{\sum_{i=1}^n\sigma_i}e^{(-1)^{\sigma_i}\sum_{i=1}^n \frac{\mathrm{H}(u_i)}{2}}}{\sqrt{\pm\mathbb{T}_{(M\setminus L,\mathbf m)}([\rho_{M\setminus L}])}}e^{\frac{r}{4\pi}\big(\mathrm{Vol}(M)+\mathbf i\mathrm{CS}(M)\big)}\Big( 1 + O \Big( \frac{1}{r} \Big) \Big),$$
where  $\rho_{M\setminus L}$ is the holonomy representation of the restriction of the hyperbolic metric on $M$ to $M\setminus L.$ 
\end{enumerate}
\end{proposition}

\begin{proof} The proof of (1) follows verbatim that of \cite[Proposition 5.5]{WY1}, namely, deform the contour of the integral to satisfy all the conditions of Proposition \ref{saddle}. The proof of (2) is by essentially the same argument, and we include a sketch here for the readers' convenience.

As depicted in Figure \ref{surface},
\begin{figure}[htbp]
\centering
\includegraphics[scale=0.3]{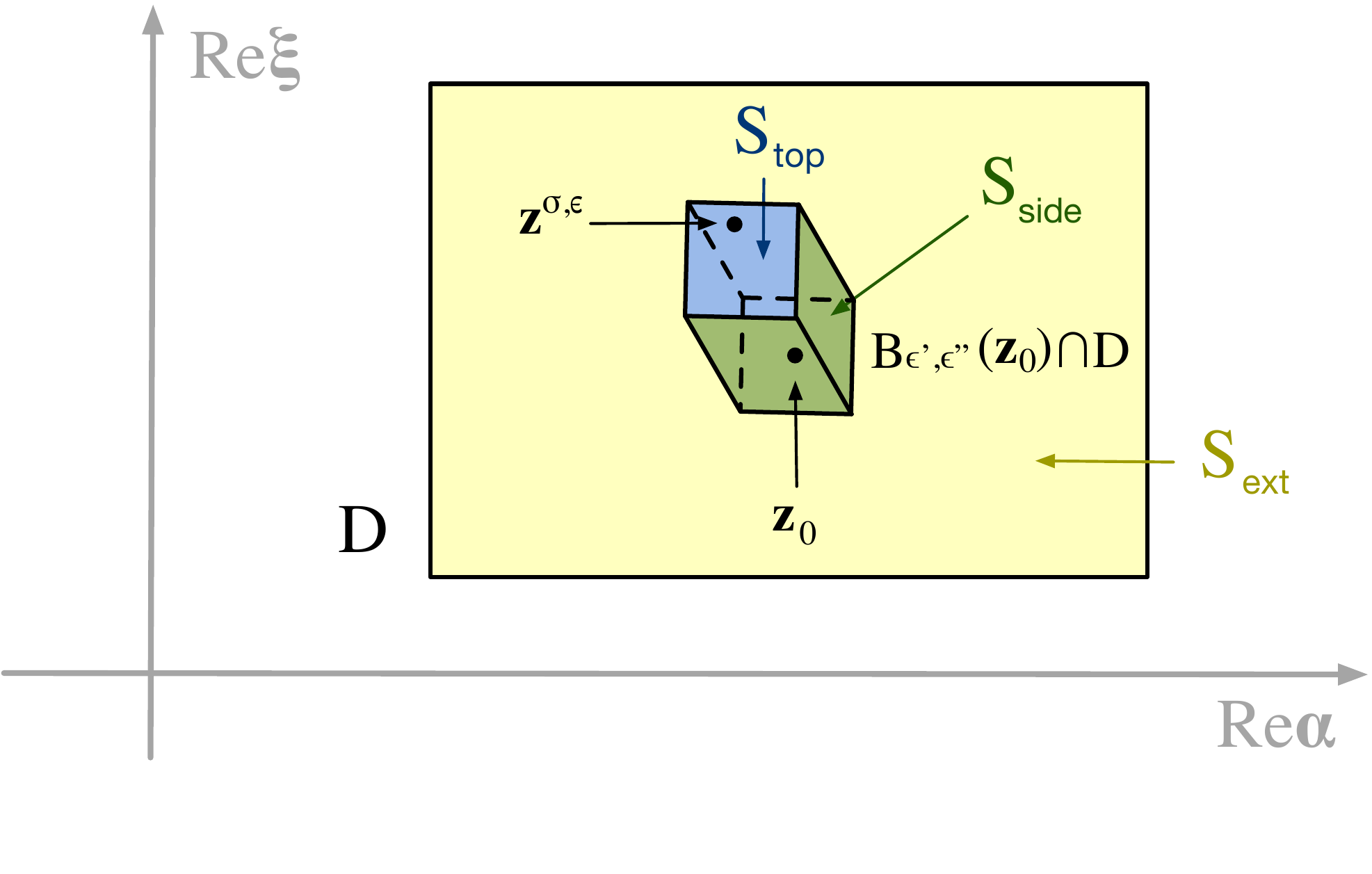}
\caption{The deformed surface $S,$ where $S_{\text{top}}$ is the ``top" of the rectangular cylinder colored in blue, $S_{\text{side}}$ is the ``side" of the rectangular cylinder colored in green, $B_{\epsilon',\epsilon''}(\boldsymbol z_0)\cap \mathrm D$ is the ``bottom" of the rectangular cylinder, and $S_\text{ext}$ is the exterior of $B_{\epsilon',\epsilon''}(\boldsymbol z_0)$ in $\mathrm D$ colored in light yellow.}
\label{surface} 
\end{figure}
let 
\begin{equation*}
\begin{split}
\mathrm S_\text{top}= & \Big\{\boldsymbol z + \mathrm{Im}\boldsymbol z^{\boldsymbol \sigma,\boldsymbol \epsilon}  \ \Big|\  \boldsymbol z \in B_{\epsilon',\epsilon''}(\boldsymbol z_0)\cap \mathrm D \Big\},\\
\mathrm S_\text{side}= & \Big\{\boldsymbol z + t\cdot \mathbf i\mathrm{Im}\boldsymbol z^{\boldsymbol \sigma,\boldsymbol \epsilon}  \ \Big|\  \boldsymbol z \in \partial B_{\epsilon',\epsilon''}(\boldsymbol z_0)\cap \mathrm D, t\in [0,1]\Big\},\\
\mathrm S_\text{ext}=&\mathrm D\setminus B_{\epsilon',\epsilon''}(\boldsymbol z_0);    
\end{split}    
\end{equation*}
and let 
$$\mathrm S=\mathrm S_\text{top}\cup \mathrm S_\text{side}\cup \mathrm S_\text{ext}.$$
Then all the conditions of Proposition \ref{saddle2}
are satisfied with $D=\mathrm D_{\mathbb C},$ $f(\boldsymbol z)= \frac{\mathcal W^{\boldsymbol\epsilon}(\boldsymbol z)-4\pi\sum_{i=1}^n\sigma_i\epsilon_i\alpha_i}{4\pi\mathbf i},$ $g(\boldsymbol z)= e^{\sum_{i=1}^n\epsilon_i(\alpha_i-\pi)\mathbf i + \frac{\kappa(\boldsymbol z)}{4\pi\mathbf i}},$ $f_r(\boldsymbol z)= \frac{\mathcal W_r^{\boldsymbol\epsilon}(\boldsymbol z)-4\pi\sum_{i=1}^n\sigma_i\epsilon_i\alpha_i}{4\pi\mathbf i},$ 
$\upsilon_r(\boldsymbol z)=\frac{\nu_r(\boldsymbol z)}{4\pi \mathbf i},$  $S=\mathrm S$  and $\boldsymbol c=\boldsymbol z^{\boldsymbol\sigma,\boldsymbol \epsilon}.$ Indeed, (i) follows from Proposition \ref{crit} (2), (ii) and (iii) follow from Proposition \ref{convexity} and assumption (\ref{assump}), (iv) follows from Proposition \ref{Hessian}, (v) follows from that $g$ is an exponential function, and (vi) follows from Proposition \ref{EST} (1). See the proof of \cite[Proposition 5.5]{WY1} for more details of the argument. As a consequence, we have 
\begin{equation*}
\begin{split}
&\widehat{f^{\boldsymbol \epsilon}_r}(\boldsymbol\sigma\boldsymbol\epsilon,\boldsymbol 0)\\
=& 
\frac{\Big(\frac{r}{2\pi}\Big)^{\frac{n+c}{2}}e^{\sum_{i=1}^n\epsilon_i(\alpha_i^{\boldsymbol\sigma,\boldsymbol\epsilon}-\pi)\mathbf i + \frac{\kappa(\boldsymbol z^{\boldsymbol \sigma,\boldsymbol\epsilon})}{4\pi\mathbf i}}}{\sqrt{\det\bigg(-\mathrm{Hess}\Big(\frac{\mathcal W^{\boldsymbol\epsilon}(\boldsymbol z)-4\pi\sum_{i=1}^n\sigma_i\epsilon_i\alpha_i}{4\pi\mathbf i}\Big)\Big|_{\boldsymbol z=\boldsymbol z^{\boldsymbol \sigma,\boldsymbol \epsilon}}\bigg)}} e^{\frac{r}{4\pi\mathbf i}\big(\mathcal W^{\boldsymbol\epsilon}(\boldsymbol z^{\boldsymbol \sigma,\boldsymbol\epsilon})-4\pi\sum_{i=1}^n\sigma_i\epsilon_i\alpha_i^{\boldsymbol\sigma,\boldsymbol\epsilon}\big)}\Big( 1 + O \Big( \frac{1}{r} \Big) \Big)\\
=&e^{-\frac{(n+c)\pi\mathbf i }{4}}r^{\frac{n+c}{2}}\frac{(-1)^{\sum_{i=1}^n\sigma_i}e^{(-1)^{\sigma_i}\sum_{i=1}^n \frac{\mathrm{H}(u_i)}{2}}}{\sqrt{\pm\mathbb{T}_{(M\setminus L,\mathbf m)}([\rho_{M\setminus L}])}}e^{\frac{r}{4\pi}\big(\mathrm{Vol}(M)+\mathbf i\mathrm{CS}(M)\big)}\Big( 1 + O \Big( \frac{1}{r} \Big) \Big),    
\end{split}    
\end{equation*}
where the last equality comes from Proposition \ref{crit} (2) and Proposition \ref{Hessian} (2).
\end{proof}


\subsection{Estimate of the other Fourier coefficients} \label{ot}

\begin{proposition}\label{other} 
There exists an $\epsilon'''>0$  such that: 
\begin{enumerate}[(1)]
\item   If $d_{L^\infty}(\boldsymbol z^{\boldsymbol \epsilon},\boldsymbol z_0)<\epsilon'''$ and 
$$\mathrm{Vol}(M_{L_{\boldsymbol \theta}})>\max_{\mathrm{D}\setminus  B_{\epsilon'''}(\boldsymbol z_0)} \mathrm{Im}\mathcal W_{\boldsymbol \beta}^{\boldsymbol \epsilon},$$
then there exists a $\delta'>0$ such that for any $(\boldsymbol m,\boldsymbol n)\neq (\boldsymbol 0,\boldsymbol 0)$ and $k>n+c,$
$$\Big|\widehat{f^{\boldsymbol \epsilon}_{\boldsymbol b,r}}(\boldsymbol m,\boldsymbol n)\Big| <\frac{O\Big(e^{\frac{r}{4\pi}\big(\mathrm{Vol}(M_{L_{\boldsymbol\theta}})-\delta'\big)}\Big)}{\Big(\sum_{i=1}^nm_i^2+\sum_{s=1}^cn_s^2\Big)^k}.$$

\item  
If $d_{L^\infty}(\boldsymbol z^{\boldsymbol\sigma,\boldsymbol \epsilon}, \boldsymbol z_0)<\epsilon'''$ and 
\begin{equation}\label{assump2}
\mathrm{Vol}(M)
>\max_{\mathrm{D}\setminus B_{\epsilon'''}(\boldsymbol z_0)} \mathrm{Im}\mathcal W^{\boldsymbol \epsilon},
\end{equation}
then  there exists a $\delta'>0$ such that for any $(\boldsymbol m,\boldsymbol n)\neq (\boldsymbol \sigma\boldsymbol\epsilon,\boldsymbol 0)$ and $k>n+c$
$$\Big|\widehat{f^{\boldsymbol \epsilon}_r}(\boldsymbol m,\boldsymbol n)\Big| <\frac{O\Big(e^{\frac{r}{4\pi}\big(\mathrm{Vol}(M)-\delta'\big)}\Big)}{\Big(\sum_{i=1}^nm_i^2+\sum_{s=1}^cn_s^2\Big)^k}.$$
\end{enumerate}
\end{proposition}

\begin{proof} We will prove (2) below; and that of (1) follows verbatim. 
\\

By the Integration by Parts as on \cite[Page 703]{O1}, one has
$$\widehat{f^{\boldsymbol \epsilon}_r}(\boldsymbol m,\boldsymbol n)=\frac{r^{n+c-k}}{\Big(\sum_{i=1}^nm_i^2+\sum_{s=1}^cn_s^2\Big)^k}\int_\mathrm D \Psi(\boldsymbol \alpha,\boldsymbol\xi)e^{\frac{r}{4\pi \mathbf i}\big(\mathcal W_r^{\boldsymbol \epsilon}(\boldsymbol{\alpha},\boldsymbol{\xi})-\sum_{i=1}^n4\pi m_i\alpha_i-\sum_{s=1}^c4\pi n_s\xi_s\big)}d\boldsymbol{\alpha}d\boldsymbol{\xi},$$
where $\Psi(\boldsymbol\alpha,\boldsymbol\xi)$ is a polynomial of the partial derivatives of order at most $2k$ of $\psi(\boldsymbol\alpha,\boldsymbol \xi)e^{\sum_{i=1}^n\epsilon_i(\alpha_i-\pi)\mathbf i}$ and $\mathcal W_r^{\boldsymbol \epsilon};$ and it  suffices to bound each $\mathcal W_r^{\boldsymbol \epsilon}(\boldsymbol z)-4\pi\sum_{i=1}^nm_i\alpha_i-4\pi\sum_{s=1}^cn_s\xi_s.$
\medskip

For the simplicity of notations, for  $\boldsymbol \epsilon\in\{-1,1\}^n$ and $(\boldsymbol m, \boldsymbol n)\in \mathbb Z^{n+c},$ we let $\mathcal W^{\boldsymbol \epsilon}_{\boldsymbol m,\boldsymbol n}$ be the function given by 
$$\mathcal W^{\boldsymbol \epsilon}_{\boldsymbol m,\boldsymbol n}(\boldsymbol z)\doteq \mathcal W^{\boldsymbol \epsilon}(\boldsymbol z)-4\pi\sum_{i=1}^nm_i\alpha_i-4\pi\sum_{s=1}^cn_s\xi_s.$$

We first observe that for all $\boldsymbol \epsilon\in\{-1,1\}^n$ and $(\boldsymbol m, \boldsymbol n)\in \mathbb Z^{n+c},$ the functions $\mathrm{Im}\mathcal W^{\boldsymbol \epsilon}_{\boldsymbol m,\boldsymbol n}$ share the same second derivatives, as they only differ by a linear function. As a consequence, there is an $N>0$ such that for any $\boldsymbol \epsilon\in\{-1,1\}^n,$ $(\boldsymbol m, \boldsymbol n)\in \mathbb Z^{n+c},$
unit vector $\mathbf u\in \mathbb C^{n+c}$ and $\boldsymbol z\in B_{\epsilon',\epsilon''}(\boldsymbol z_0),$  
the second directional derivative 
$$\big| D_\mathbf u D_\mathbf u \mathrm{Im}\mathcal W^{\boldsymbol \epsilon}_{\boldsymbol m,\boldsymbol n}(\boldsymbol z) \big|< N.$$
As a consequence, for any $\boldsymbol z_1, \boldsymbol z_2\in  B_{\epsilon',\epsilon''}(\boldsymbol z_0)$ and for the unit vector $\mathbf u=\frac{\boldsymbol z_1-\boldsymbol z_2}{|\boldsymbol z_1-\boldsymbol z_2|},$  the directional derivative 
\begin{equation}\label{dfd}
\big|D_\mathbf u \mathrm{Im}\mathcal W^{\boldsymbol \epsilon}_{\boldsymbol m,\boldsymbol n}(\boldsymbol z_1)-D_\mathbf u \mathrm{Im}\mathcal W^{\boldsymbol \epsilon}_{\boldsymbol m,\boldsymbol n}(\boldsymbol z_2)\big|< N|\boldsymbol z_1-\boldsymbol z_2|\leqslant N\sqrt{2n+2c}\cdot d_{L^\infty}(\boldsymbol z_1,\boldsymbol z_2).
\end{equation}

Next, let 
$$\epsilon'''=\min\Big\{\epsilon',\epsilon'', \frac{\pi}{4N(n+c)}\Big\}.$$
Then as long as $d_{L^\infty}(\boldsymbol z^{\boldsymbol\sigma,\boldsymbol \epsilon},\boldsymbol z_0)<\epsilon''',$ we have $\boldsymbol z^{\boldsymbol\sigma,\boldsymbol \epsilon}\in B_{\epsilon',\epsilon''}(\boldsymbol z_0);$ and for any $\boldsymbol z\in B_{\epsilon'''}(\boldsymbol z_0),$ 
\begin{equation}\label{d<}
    d_{L^\infty}(\boldsymbol z,\boldsymbol z^{\boldsymbol \sigma,\boldsymbol \epsilon})<2\epsilon ''' \leqslant \frac{\pi}{2N(n+c)}.
\end{equation}
By Proposition \ref{crit} (2), the total derivative 
$$D\mathrm{Im}\mathcal W^{\boldsymbol \epsilon}_{\boldsymbol \sigma\boldsymbol \epsilon,\boldsymbol 0}(\boldsymbol z^{\boldsymbol \sigma,\boldsymbol\epsilon})=\mathbf 0;$$
and by (\ref{dfd}) and (\ref{d<}), we have 
\begin{equation}\label{Dest}
    \big|D_{\mathbf u}\mathrm{Im}\mathcal W^{\boldsymbol \epsilon}_{\boldsymbol \sigma\boldsymbol \epsilon,\boldsymbol 0}(\boldsymbol z)\big|<\frac{\pi}{\sqrt{2n+2c}}
\end{equation}
for any $\boldsymbol z\in B_{\epsilon'''}(\boldsymbol z_0)$ and $\mathbf u=\frac{\boldsymbol z-\boldsymbol z^{\boldsymbol \sigma,\boldsymbol \epsilon}}{|\boldsymbol z-\boldsymbol z^{\boldsymbol \sigma,\boldsymbol \epsilon}|},$ the unit vector along the direction of the line segment connecting $\boldsymbol z$ and $\boldsymbol z^{\boldsymbol \sigma,\boldsymbol \epsilon}.$

Now, since $(\boldsymbol m,\boldsymbol n)\neq (\boldsymbol \sigma\boldsymbol\epsilon,\boldsymbol 0),$ at least one of  the following four cases happens:
\begin{enumerate}[(1)]
\item $m_i>\max\{0,\epsilon_i\}$ for some $i\in\{1,\dots,n\}.$

\item $m_i<\min\{0,\epsilon_i\}$ for some $i\in\{1,\dots,n\}.$

\item $n_s>0$ for some $s\in\{1,\dots,c\}.$

\item $n_s<0$ for some $s\in\{1,\dots,c\}.$
\end{enumerate}

In case (1), we let $\mathbf  u_i\in\mathbb C^{n+c}$ be the vector with the $i$-th component $1$ and all the other components $0;$ let
\begin{equation*}
\begin{split}
\mathrm S_\text{top}= &\Big\{\boldsymbol z + \epsilon'''\mathbf i\mathbf u_i  \ \Big|\  \boldsymbol z \in B_{\epsilon'''}(\boldsymbol z_0)\cap \mathrm D \Big\},\\
\mathrm S_\text{side}= & \Big\{\boldsymbol z + t\cdot \mathbf i\mathbf u_i  \ \Big|\  \boldsymbol z \in \partial B_{\epsilon'''}(\boldsymbol z_0)\cap \mathrm D, t\in [0,\epsilon''']\Big\},\\
\mathrm S_\text{ext}= &\mathrm D\setminus B_{\epsilon'''}(\boldsymbol z_0);
\end{split}
\end{equation*}
and let 
$$\mathrm S=\mathrm S_\text{top}\cup \mathrm S_\text{side}\cup \mathrm S_\text{ext}.$$

For $\boldsymbol z\in \mathrm{S}_\text{top},$ as $\mathrm{S}_\text{top}\subset B_{\epsilon'''}(\boldsymbol z_0),$ by the Mean Value Theorem we have 
\begin{equation*}
\begin{split}
    \mathrm{Im}\mathcal W_{\boldsymbol \sigma\boldsymbol \epsilon,\boldsymbol 0}^{\boldsymbol \epsilon}(\boldsymbol z)\leqslant & \mathrm{Im}\mathcal W_{\boldsymbol \sigma\boldsymbol \epsilon,\boldsymbol 0}^{\boldsymbol \epsilon}(\boldsymbol z^{\boldsymbol \sigma,\boldsymbol \epsilon})+ \big|D_{\mathbf u}\mathrm{Im}\mathcal W_{\boldsymbol \sigma\boldsymbol \epsilon,\boldsymbol 0}^{\boldsymbol\epsilon}(\boldsymbol z')\big|\cdot  |\boldsymbol z - \boldsymbol z^{\boldsymbol\sigma,\boldsymbol \epsilon}|\\
< & \mathrm{Im}\mathcal W_{\boldsymbol \sigma\boldsymbol \epsilon,\boldsymbol 0}^{\boldsymbol \epsilon}(\boldsymbol z^{\boldsymbol \sigma,\boldsymbol \epsilon})+2\pi\epsilon''',
  \end{split}
\end{equation*}
where $\boldsymbol z'$ is some point on the line segment from $\boldsymbol z$ to $\boldsymbol z^{\boldsymbol \sigma,\boldsymbol \epsilon}$ and $\mathbf u$ is the unit vector along the direction of  this line segment, and the last inequality comes from (\ref{Dest}) and that $|\boldsymbol z - \boldsymbol z^{\boldsymbol\sigma,\boldsymbol\epsilon}|$ is less than the diameter of $B_{\epsilon'''}(\boldsymbol z_0),$ which equals $2\epsilon'''\sqrt{2n+2c}.$
As a consequence, 
\begin{equation}\label{top}
\begin{split}
   \mathrm{Im}\mathcal W_{\boldsymbol m,\boldsymbol n}^{\boldsymbol\epsilon}(\boldsymbol z)
= & \mathrm{Im}\mathcal W^{\boldsymbol\epsilon}(\boldsymbol z)-4\pi m_i\epsilon''' \\
= & \mathrm{Im}\mathcal W_{\boldsymbol \sigma\boldsymbol \epsilon,\boldsymbol 0}^{\boldsymbol\epsilon}(\boldsymbol z)-4\pi m_i\epsilon''' + 4\pi\sigma_i\epsilon_i\epsilon'''\\
< & \mathrm{Im}\mathcal W_{\boldsymbol \sigma\boldsymbol \epsilon,\boldsymbol 0}^{\boldsymbol \epsilon}(\boldsymbol z^{\boldsymbol \sigma,\boldsymbol \epsilon}) +2\pi\epsilon''' - 4\pi m_i\epsilon'''+ 4\pi\sigma_i\epsilon_i\epsilon'''\\
= & \mathrm{Im}\mathcal W_{\boldsymbol \sigma\boldsymbol \epsilon,\boldsymbol 0}^{\boldsymbol \epsilon}(\boldsymbol z^{\boldsymbol \sigma,\boldsymbol \epsilon}) -2\pi \epsilon''' - 4\pi (m_i-\sigma_i\epsilon_i-1)\epsilon'''\\
\leqslant  & \mathrm{Im}\mathcal W_{\boldsymbol \sigma\boldsymbol \epsilon,\boldsymbol 0}^{\boldsymbol \epsilon}(\boldsymbol z^{\boldsymbol \sigma,\boldsymbol \epsilon})- 2\pi\epsilon''',
\end{split}
\end{equation}
where the last equality comes from that $m_i\geqslant \epsilon_i\sigma_i+ 1.$

For $\boldsymbol z\in \mathrm{S}_\text{side},$ let $$\delta''=\mathrm{Vol}(M)
-\max_{\mathrm{D}\setminus B_{\epsilon'''}(\boldsymbol z_0)} \mathrm{Im}\mathcal W^{\boldsymbol \epsilon}>0.$$
Then by assumption (\ref{assump2}) and Proposition \ref{crit} (2), for $\boldsymbol z'\in \partial B_{\epsilon'''}(\boldsymbol z_0)\cap \mathrm D,$ we have
\begin{equation}\label{sidebottom}
 \mathrm{Im}\mathcal W_{\boldsymbol m,\boldsymbol n}^{\boldsymbol\epsilon}(\boldsymbol z') =   \mathrm{Im}\mathcal W^{\boldsymbol\epsilon}(\boldsymbol z')
\leqslant  \mathrm{Im}\mathcal W_{\boldsymbol \sigma\boldsymbol\epsilon,\boldsymbol 0}^{\boldsymbol \epsilon}(\boldsymbol z^{\boldsymbol \sigma,\boldsymbol \epsilon})-\delta''.
\end{equation}
By Proposition \ref{convexity}, we have that $\mathrm{Im}\mathcal W_{\boldsymbol m,\boldsymbol n}^{\boldsymbol \epsilon}(\boldsymbol z'+t\cdot \mathbf i\mathbf u_i)$ is strictly concave up in $t$ for any $\boldsymbol z'\in \partial B_{\epsilon'''}(\boldsymbol z_0);$ and by (\ref{top}) and  (\ref{sidebottom}), we have
\begin{equation}\label{side}
    \mathrm{Im}\mathcal W_{\boldsymbol m,\boldsymbol n}^{\boldsymbol\epsilon}(\boldsymbol z) < \mathrm{Im}\mathcal W_{\boldsymbol \sigma\boldsymbol \epsilon,\boldsymbol 0}^{\boldsymbol \epsilon}(\boldsymbol z^{\boldsymbol \sigma,\boldsymbol \epsilon}) - \min\{2\pi\epsilon''',\delta''\}.
\end{equation}

For $\boldsymbol z\in \mathrm{S}_\text{ext},$ by the assumption (\ref{assump2}) and Proposition \ref{crit} (2), we have 
\begin{equation}\label{ext}
\begin{split}
   \mathrm{Im}\mathcal W_{\boldsymbol m,\boldsymbol n}^{\boldsymbol\epsilon}(\boldsymbol z)
=  \mathrm{Im}\mathcal W^{\boldsymbol\epsilon}(\boldsymbol z)
\leqslant &   \mathrm{Im}\mathcal W_{\boldsymbol \sigma\boldsymbol\epsilon,\boldsymbol 0}^{\boldsymbol \epsilon}(\boldsymbol z^{\boldsymbol \sigma,\boldsymbol \epsilon})-\delta''.
\end{split}
\end{equation}

Finally, let
$$\delta'=\min\{2\pi \epsilon''', \delta''\};$$ 
and notice that, as both $\epsilon'''$ and $\delta''$ are independent of $(\boldsymbol m,\boldsymbol n),$   $\delta'$ is so as well.
Putting (\ref{top}), (\ref{side}) and (\ref{ext}) together,  we have 
\begin{equation*}
 \mathrm{Im}\mathcal W_{\boldsymbol m,\boldsymbol n}^{\boldsymbol\epsilon}(\boldsymbol z) <\mathrm{Im}\mathcal W_{\boldsymbol \sigma\boldsymbol\epsilon,\boldsymbol 0}^{\boldsymbol \epsilon}(\boldsymbol z^{\boldsymbol \sigma,\boldsymbol \epsilon})-\delta' =\mathrm{Vol}(M)-\delta'
\end{equation*}
on $\mathrm S,$ where the equality comes from Proposition \ref{crit} (2). Then the result follows from Proposition \ref{EST} (1).

In case (2), we deform $B_{\epsilon'''}(\boldsymbol z_0)\cap\mathrm D$ along the $-\mathbf i\mathbf u_i$ direction by $\epsilon''';$ and in case (3) and case (4), we deform $B_{\epsilon'''}(\boldsymbol z_0)\cap\mathrm D$  respectively along the $\pm\mathbf i \mathbf u_{n+s}$ directions by $\epsilon'''.$ Then the rest of the proof follows verbatim that of case (1).
\end{proof}


\subsection{Estimate of the error term}\label{ee}

\begin{proposition}\label{errer} 
\begin{enumerate}[(1)]
\item If 
$$\mathrm{Vol}(M_{L_{\boldsymbol \theta}})>\max_{\partial\mathrm{D}} \mathrm{Im}\mathcal W_{\boldsymbol \beta}^{\boldsymbol \epsilon},$$
then for any 
$$\delta'< \frac{1}{2}\Big(\mathrm{Vol}(M_{L_{\boldsymbol \theta}})-\max_{\partial\mathrm{D}} \mathrm{Im}\mathcal W_{\boldsymbol \beta}^{\boldsymbol \epsilon}\Big),$$
there exists an $\epsilon''''>0$ such that the error term in (\ref{P1}) is less than
$O\Big(e^{\frac{r}{4\pi}\big(\mathrm{Vol}(M_{L_{\boldsymbol\theta}})-\delta'\big)}\Big)$
on $\mathrm D\setminus \mathrm D_{\epsilon''''}.$

\item If 
\begin{equation*}
\mathrm{Vol}(M)
>\max_{\partial \mathrm{D}} \mathrm{Im}\mathcal W^{\boldsymbol \epsilon},
\end{equation*}
then for any 
$$\delta'< \frac{1}{2}\Big(\mathrm{Vol}(M)-\max_{\partial\mathrm{D}} \mathrm{Im}\mathcal W^{\boldsymbol \epsilon}\Big),$$
there exists an $\epsilon''''>0$ such that the error term in (\ref{P2}) is less than
$O\Big(e^{\frac{r}{4\pi}\big(\mathrm{Vol}(M)-\delta'\big)}\Big)$
on $\mathrm D\setminus \mathrm D_{\epsilon''''}.$
\end{enumerate}
\end{proposition}

\begin{proof}
By Proposition \ref{EST} (2), the relation (\ref{LL}), and the continuity of the Lobachevsky function $\Lambda,$ we can choose a sufficiently small $\epsilon''''>0$ so that if $\big(\frac{2\pi (\boldsymbol a+\boldsymbol 1)}{r},\frac{2\pi \boldsymbol k}{r}\big)\in \mathrm D\setminus\mathrm D_{\epsilon''''},$ then for (1),
$$\Big|g^{\boldsymbol\epsilon}_{\boldsymbol b,r}(\boldsymbol a, \boldsymbol  k)\Big|<O\Big(e^{\frac{r}{4\pi}\big(\max_{\partial\mathrm{D}} \mathrm{Im}\mathcal W_{\boldsymbol \beta}^{\boldsymbol \epsilon}+\delta'\big)}\Big)<O\Big(e^{\frac{r}{4\pi} \big(\mathrm{Vol}(M_{L_{\boldsymbol\theta}})-\delta'\big)}\Big),$$
and for (2), 
$$\Big|g^{\boldsymbol\epsilon}_r(\boldsymbol a, \boldsymbol  k)\Big|<O\Big(e^{\frac{r}{4\pi}\big(\max_{\partial\mathrm{D}} \mathrm{Im}\mathcal W^{\boldsymbol \epsilon}+\delta'\big)}\Big)<O\Big(e^{\frac{r}{4\pi} \big(\mathrm{Vol}(M)-\delta'\big)}\Big).$$
Let $\psi$ be the bump function supported on $(\mathrm D, \mathrm D_{\epsilon''''}).$ Then the error term in (\ref{P1}) is less than $O\big(e^{\frac{r}{4\pi}(\mathrm{Vol}(M_{L_{\boldsymbol\theta}})-\delta')}\big),$ and the error term in (\ref{P2}) is less than $O\big(e^{\frac{r}{4\pi}(\mathrm{Vol}(M)-\delta')}\big).$
\end{proof}


\subsection{Proof of Theorem \ref{main1} and Theorem \ref{main2}}\label{pf}

\begin{proof}[Proof of Theorem \ref{main1}]
We will prove Theorem \ref{main1} by proving Theorem \ref{re1}.

Let $\delta_1$ and $\delta_2$  respectively be as in Propositions \ref{bdky} and \ref{key}, and let $$\delta\in (0,\min\{\delta_1, \delta_2\})$$
be sufficiently small so that if 
$$\mathrm {Vol}(M_{L_{\boldsymbol \theta}})>\mathrm{Vol}(M_\text{FSL})-\delta,$$
then all the conditions of Proposition \ref{convexity} and (1) of Propositions \ref{crit},  \ref{critical}, \ref{other} and \ref{errer} are satisfied.
The existence of such $\delta$ is guaranteed by Thurston's Hyperbolic Dehn-filling Theorem (see e.g. \cite[Section 15.4.5]{M}).
Then by Propositions \ref{convexity}, \ref{crit}, \ref{Hessian}, \ref{critical}, \ref{other} and \ref{errer}, 
\begin{equation*}
\begin{split}
\mathrm{RT}_r\big(M,L,\boldsymbol b\big)=& B_r\sum_{\boldsymbol\epsilon\in\{-1,1\}^n}\widehat{ f^{\boldsymbol \epsilon}_{\boldsymbol b,r}}(\boldsymbol 0,\boldsymbol 0)+O\Big(e^{\frac{r}{4\pi}(\mathrm{Vol}(M_\text{FSL})-\delta)}\Big)\\
=&\frac{C_r}{\sqrt{\pm\mathbb{T}_{(M\setminus L,\mathbf m)}([\rho_{M^{(r)}}])}}e^{\frac{r}{4\pi}\big(\mathrm{Vol}(M^{(r)})+\mathbf i\mathrm{CS}(M^{(r)})\big)}\Big( 1 + O \Big( \frac{1}{r} \Big) \Big),
\end{split}
\end{equation*}
where $$C_r=e^{\big(\frac{r+1}{4}\sigma(U_{c+1}\cup L)-(\frac{r}{4}+1)\sum_{i=1}^n(p_i+q_i)+\frac{1}{2}\sum_{i=1}^n\iota_i-\frac{r}{2}(n-c)-\frac{3n}{4}-\frac{c}{4}\big)\pi\mathbf i}$$
is a complex number of norm $1$ independent of the geometry on $(M,L).$  
\end{proof}

\begin{proof}[Proof of Theorem \ref{main2}]
Let $\delta_1$ and $\delta_2$  respectively be as in Propositions \ref{bdky} and \ref{key}, and let $$\delta\in (0,\min\{\delta_1,\delta_2\})$$
be sufficiently small so that  if 
$$\mathrm {Vol}(M)>\mathrm{Vol}(M_\text{FSL})-\delta,$$
then all the conditions of Proposition  \ref{convexity} and (2) of Propositions \ref{crit},  \ref{critical}, \ref{other} and \ref{errer} are satisfied. Again, the existence of  $\delta$ is guaranteed by Thurston's Hyperbolic Dehn-filling Theorem (see e.g. \cite[Section 15.4.5]{M}). Then by Propositions \ref{convexity}, \ref{crit}, \ref{Hessian}, \ref{critical}, \ref{other} and \ref{errer}, 
\begin{equation*}
\begin{split}
\mathrm{RT}_r\big(M)=& B_r\sum_{\boldsymbol\sigma\in\{0,1\}^n,\boldsymbol\epsilon\in\{-1,1\}^n}\widehat{ f^{\boldsymbol \epsilon}_{r}}(\boldsymbol \sigma \epsilon, \boldsymbol 0)+O\Big(e^{\frac{r}{4\pi}\big(\mathrm{Vol}(M_\text{FSL})-\delta\big)}\Big)\\
=& C_r \frac{\sum_{\boldsymbol\sigma\in\{0,1\}^n}(-1)^{\sum_{i=1}^n\sigma_i}e^{(-1)^{\sigma_i}\sum_{i=1}^n \frac{\mathrm{H}(u_i)}{2}}}{\sqrt{\pm\mathbb{T}_{(M\setminus L,\mathbf m)}([\rho_{M^{(r)}}])}}e^{\frac{r}{4\pi}\big(\mathrm{Vol}(M^{(r)})+\mathbf i\mathrm{CS}(M^{(r)})\big)}\Big( 1 + O \Big( \frac{1}{r} \Big) \Big)\\
=& C_r \frac{\prod_{i=1}^n 2\sinh \frac{\mathrm{H}(u_i)}{2}}{\sqrt{\pm\mathbb{T}_{(M\setminus L,\mathbf m)}([\rho_{M^{(r)}}])}}e^{\frac{r}{4\pi}\big(\mathrm{Vol}(M^{(r)})+\mathbf i\mathrm{CS}(M^{(r)})\big)}\Big( 1 + O \Big( \frac{1}{r} \Big) \Big)\\
= & \frac{C_r}{\sqrt{\pm\mathrm{Tor}(M)}}e^{\frac{r}{4\pi}\big(\mathrm{Vol}(M)+\mathbf i\mathrm{CS}(M)\big)}\Big(1+O\Big(\frac{1}{r}\Big)\Big),
\end{split}
\end{equation*}
where 
  
$$C_r=e^{\big(\frac{r+1}{4}\sigma(U_{c+1}\cup L)-(\frac{r}{4}+1)\sum_{i=1}^np_i+\frac{1}{2}\sum_{i=1}^n\iota_i-\frac{r}{2}(n-c)-\frac{3n}{4}-\frac{c}{4}\big)\pi\mathbf i}$$
is a complex number of norm $1$ independent of the geometry on $M,$ and the last equality comes from the Surgery Formula of the adjoint twisted Reidemeister torsion (see \cite[Theorem 4.1]{P}).
\end{proof}


\section{Weak reflection symmetry of extended quantum $6j$-symbols}\label{pfKey}

Proposition \ref{key} is a consequence of Proposition \ref{key'} (2) and Theorem \ref{wrs} below. 

To state them, we recall that a six-tuple $(a_1,\dots,a_6)\in\mathbb [0,r-2]^6$ is $r$-hyperideal if for $\{i,j,k\}=\{1,2,3\},$ $\{1,5,6\},$ $\{2,4,6\}$ and $\{3,4,5\},$
\begin{enumerate}[(1)]
\item $0\leqslant a_i+a_j-a_k\leqslant r-2,$ and
\item $r-2\leqslant a_i+a_j+a_k\leqslant 2(r-2);$
\end{enumerate}
and a six-tuple $(\alpha_1,\dots,\alpha_6)\in \mathbb (0,2\pi)^6$ is of the hyperideal type if for $\{i,j,k\}=\{1,2,3\},$ $\{1,5,6\},$ $\{2,4,6\}$ and $\{3,4,5\},$
\begin{enumerate}[(1)]
\item $0< \alpha_i+\alpha_j-\alpha_k<2\pi,$ and
\item $2\pi<\alpha_i+\alpha_j+\alpha_k<4\pi.$
\end{enumerate}

For $(a_1,\dots,a_6)\in [0,r-2]^6$ and each $k\in \{1,\dots,6\},$ we call 
$$\alpha_k=\frac{2\pi (a_k+1)}{r}$$
the corresponding $\alpha_k.$
Then $(a_1,\dots,a_6)$ being $r$-hyperideal implies that the corresponding  $(\alpha_1,\dots,\alpha_6)$ is of the hyperideal type. In the rest of this paper, we let 
$$\mathrm H=\Big\{ (\alpha_1,\dots,\alpha_6)\in (0,2\pi)^6 \ \Big|\ (\alpha_1,\dots,\alpha_6) \text{ is of the hyperideal type}\Big\};$$
and for $\epsilon_0>0$ sufficiently small, let
$$\mathrm H_{\epsilon_0}=\Big\{ (\alpha_1,\dots,\alpha_6)\in \mathrm H \ \Big|\ d\big((\alpha_1,\dots,\alpha_6),\partial \mathrm H\big)>\epsilon_0\Big\}.$$

For $\boldsymbol\alpha=(\alpha_1,\dots,\alpha_6)\in\mathrm H$ and each $k\in\{1,\dots,6\},$  let 
$$\theta_k=|\alpha_k-\pi|.$$
Then by \cite{BB}, $(\theta_1,\dots,\theta_6)$ form the set of dihedral angles of a truncated hyperideal tetrahedron $\Delta_{\boldsymbol\alpha};$ and we denote the hyperbolic volume of $\Delta_{\boldsymbol\alpha}$ by $\mathrm{Vol}(\Delta_{\boldsymbol\alpha}).$

\begin{proposition}\label{key'}
\begin{enumerate}[(1)]
    \item For $\epsilon_0>0$ sufficiently small and any $r$-hyperideal six-tuple $(a_1,\dots,a_6)$ with the corresponding $(\alpha_1,\dots,\alpha_6)\in \mathrm H\setminus \mathrm H_{\epsilon_0},$ 
      $$  \Bigg|\bigg\|\begin{matrix} 
      a_1 & a_2 & a_3 \\
      a_4 & a_5 & a_6
 \end{matrix}\bigg\|\Bigg| < O\Big(e^{\frac{r}{2\pi}\big(\mathrm{Vol}(\Delta_{\boldsymbol \alpha})+\delta\big)}\Big)$$
 for any $\delta>0.$

    \item For $\epsilon_0>0$ sufficiently small, there exists a $\delta_3>0$ depending only on $\epsilon$ such that for any $r$-hyperideal six-tuple $(a_1,\dots,a_6)$ with the corresponding $(\alpha_1,\dots,\alpha_6)\in \mathrm H\setminus \mathrm H_{\epsilon_0},$ 
      $$  \Bigg|\bigg\|\begin{matrix} 
      a_1 & a_2 & a_3 \\
      a_4 & a_5 & a_6
 \end{matrix}\bigg\|\Bigg| < O\Big(e^{\frac{r}{2\pi}\big(v_8-\delta_3\big)}\Big),$$
 where $v_8$ is the volume of the regular ideal hyperbolic octahedron.
\end{enumerate}
\end{proposition}

Proposition \ref{key'} will be proved in Section \ref{routine}.

\begin{theorem}[Weak Reflection Symmetry]\label{wrs}
For $\epsilon_0>0$ sufficiently small, there exists a $\delta_4>0$ depending only on $\epsilon_0$ such that for any $r$-hyperideal six-tuple $(a_1,\dots,a_6)$ with the corresponding $\boldsymbol\alpha=(\alpha_1,\dots,\alpha_6)\in \mathrm H_{\epsilon_0},$  we have
      $$  \Bigg|\bigg\|\begin{matrix} 
      a_1 & a_2 & a_3 \\
      a_4 & a_5 & a_6
 \end{matrix}\bigg\|-\bigg\|\begin{matrix} 
      r-2-a_1 & a_2 & a_3 \\
      a_4 & a_5 & a_6
   \end{matrix}\bigg\|\Bigg| < O\Big(e^{\frac{r}{2\pi}\big(\mathrm{Vol}
(\Delta_{\boldsymbol\alpha})-\delta_4\big)}\Big).$$
Moreover, the same result holds if we change $a_k\mapsto r-2-a_k$ for any $k\in\{1,\dots,6\}.$
\end{theorem}

 Theorem \ref{wrs} is a consequence of Proposition \ref{error6j}, Proposition \ref{other6j} and Proposition \ref{leading6j} below.

 To state them, we let $\epsilon>0$ be sufficiently small so that for each $\boldsymbol \alpha \in \mathrm H_{\epsilon_0},$ the interval $[\xi_{\boldsymbol \alpha}-3\epsilon,\xi_{\boldsymbol \alpha}+3\epsilon]$ entirely lies in $(\max\{\tau_i\},\min\{\eta_j\}).$ Recall that 
 \begin{equation*}
\begin{split}
    &\tau_1=\frac{\alpha_1+\alpha_2+\alpha_3}{2},   \quad \tau_2=\frac{\alpha_1+\alpha_5+\alpha_6}{2}, \quad \tau_3=\frac{\alpha_2+\alpha_4+\alpha_6}{2}, \quad \tau_4=\frac{\alpha_3+\alpha_4+\alpha_5}{2},\\
&\eta_1=\frac{\alpha_1+\alpha_2+\alpha_4+\alpha_5}{2}, \ \ \eta_2=\frac{\alpha_1+\alpha_3+\alpha_4+\alpha_6}{2}, 
\ \ \eta_3=\frac{\alpha_2+\alpha_3+\alpha_5+\alpha_6}{2}\quad\text{and}\quad \eta_4=2\pi;  
\end{split}
\end{equation*} and keep in mind that $\tau_i$'s and $\eta_j$'s are quantities depending on $\boldsymbol \alpha.$ For $k\in\mathbb Z,$ we call 
$$\xi=\frac{2k\pi}{r}$$
the corresponding $\xi.$

 \begin{proposition}\label{error6j}
There exists a $\delta_5>0$ depending only on $\epsilon$ such that for any $r$-hyperideal six-tuple $(a_1,\dots,a_6)$ with the corresponding $(\alpha_1,\dots,\alpha_6)\in\mathrm H_{\epsilon_0}$ and integer $k$ with the corresponding $\xi\in (\max\{\tau_i\},\min\{\eta_j\})\setminus [\xi_{\boldsymbol \alpha}-\epsilon,\xi_{\boldsymbol \alpha}+\epsilon],$  we have 
$$\Big|e^{\frac{r}{4\pi\mathbf i}U_r(\boldsymbol a, k)}\Big|<O\Big(e^{\frac{r}{2\pi} \big(\mathrm{Vol}(\Delta_{\boldsymbol\alpha})-\delta_5\big)}\Big).$$
\end{proposition}

Proposition \ref{error6j} will be proved in Section \ref{routine}.
\\

For the $r$-hyperideal six-tuple $(a_1,\dots,a_6)$ with the corresponding $(\alpha_1,\dots,\alpha_6)\in\mathrm H_{\epsilon_0}$ and integer $k$ with the corresponding $\xi\in  [\xi_{\boldsymbol \alpha}-\epsilon,\xi_{\boldsymbol \alpha}+\epsilon],$  we will use the Poisson Summation Formula. Let $\psi:\mathbb R\to\mathbb R$ be a $C^\infty$-smooth bump function supported on $([\xi_{\boldsymbol \alpha}-2\epsilon,\xi_{\boldsymbol \alpha}+2\epsilon],[\xi_{\boldsymbol \alpha}-\epsilon,\xi_{\boldsymbol \alpha}+\epsilon]),$ i.e., 
\begin{equation*}
\left \{\begin{array}{rl}
\psi({\xi})=1, & \xi\in  [\xi_{\boldsymbol \alpha}-\epsilon,\xi_{\boldsymbol \alpha}+\epsilon]\\
0<\psi({\xi})<1, &  {\xi}\in [\xi_{\boldsymbol \alpha}-2\epsilon,\xi_{\boldsymbol \alpha}+2\epsilon]\setminus  [\xi_{\boldsymbol \alpha}-\epsilon,\xi_{\boldsymbol \alpha}+\epsilon]\\
\psi({\xi})=0, & {\xi}\notin [\xi_{\boldsymbol \alpha}-2\epsilon,\xi_{\boldsymbol \alpha}+2\epsilon].\\
\end{array}\right.
\end{equation*}
Then by (\ref{3.9}) and Proposition \ref{error6j},  we have 
\begin{equation}\label{6.1}
\bigg\|
\begin{matrix}
        a_1 & a_2 & a_3 \\
        a_4 & a_5 & a_6 
      \end{matrix} \bigg\|=(-1)^{\frac{r+1}{2}}\frac{2}{r}\sin\frac{2\pi}{r}\sum_{k\in \mathbb Z}\psi(\xi)e^{\frac{r}{4\pi\mathbf i}\big(U_r(\alpha_1,\dots,\alpha_6,\xi)\big)}+O\Big(e^{\frac{r}{2\pi} \big(\mathrm{Vol}(\Delta_{\boldsymbol\alpha})-\delta_5\big)}\Big) ;
\end{equation}
and by the Poisson Summation Formula, we have \begin{equation}\label{6.2}
\bigg\|
\begin{matrix}
        a_1 & a_2 & a_3 \\
        a_4 & a_5 & a_6 
      \end{matrix} \bigg\|=\frac{(-1)^{\frac{r+1}{2}}}{\pi}\sin\frac{2\pi}{r}\sum_{n\in \mathbb Z}\mathrm I_{\boldsymbol a}(n)+O\Big(e^{\frac{r}{2\pi} \big(\mathrm{Vol}(\Delta_{\boldsymbol\alpha})-\delta_5\big)}\Big),
\end{equation}
where the $n$-th Fourier coefficient 
\begin{equation}\label{Fn}
\mathrm I_{\boldsymbol a}(n)= \int_\mathbb R\psi(\xi) \exp\bigg(\frac{r}{4\pi\mathbf i}\Big(U_r(\alpha_1,\dots,\alpha_6,\xi)-4\pi n\xi\Big)\bigg) d\xi.
\end{equation}

\begin{proposition}\label{other6j}
There exists a $\delta_6>0$ such that for any $n\neq 0,$ 
$$\big|\mathrm I_{\boldsymbol a}(n)\big| < \frac{1}{k^2}\cdot  O\Big(e^{\frac{r}{2\pi} \big(\mathrm{Vol}(\Delta_{\boldsymbol\alpha})-\delta_6\big)}\Big).$$
\end{proposition}

Proposition \ref{other6j} will be proved in Section \ref{routine}.

\begin{proposition}\label{leading6j}
For $\boldsymbol a=(a_1,\dots, a_6),$ let $\boldsymbol a'$ be obtained from $\boldsymbol a$ by changing $a_k\mapsto r-2-a_k$ for some $k\in\{1,\dots,6\}.$ Then there exists a $\delta_7>0$ such that 
$$\big|\mathrm I_{\boldsymbol a}(0)- \mathrm I_{\boldsymbol a'}(0)\big| < O\Big(e^{\frac{r}{2\pi} \big(\mathrm{Vol}(\Delta_{\boldsymbol\alpha})-\delta_7\big)}\Big).$$
\end{proposition}

Proposition \ref{leading6j} will be proved in Section \ref{challenging}.


\subsection{The routine part}\label{routine}

The goal of this section is to prove Propositions \ref{key'},   \ref{error6j} and \ref{other6j}

\begin{proof}[Proof of Proposition \ref{key'}]
For (1),  by Proposition \ref{EST} (2), the relation (\ref{LL}), (\ref{VV}) and (\ref{320}), for any $\delta>0,$ we have
$$\Big|e^{\frac{r}{4\pi\mathbf i}U_r(\boldsymbol \alpha, \xi)}\Big|<O\Big(e^{\frac{r}{2\pi}\big(V(\boldsymbol \alpha,\xi_{\boldsymbol\alpha})+\delta\big)}\Big)=O\Big(e^{\frac{r}{2\pi} \big(\mathrm{Vol}(\Delta_{\boldsymbol\alpha})+\delta\big)}\Big).$$
Therefore, by (\ref{3.9}), we have
      $$  \Bigg|\bigg\|\begin{matrix} 
      a_1 & a_2 & a_3 \\
      a_4 & a_5 & a_6
 \end{matrix}\bigg\|\Bigg|\leqslant\frac{2}{r}\sin\frac{2\pi}{r} \sum_k  \Big|e^{\frac{r}{4\pi\mathbf i}U_r(\boldsymbol \alpha, \xi)}\Big|< O\Big(e^{\frac{r}{2\pi}\big(\mathrm{Vol}(\Delta_{\boldsymbol\alpha})+\delta\big)}\Big)$$
  for any $\delta_3<\delta,$
where the sum is over all the integers in between $\max\{T_i\}$ and $\min\{Q_j\}.$
\\

For (2), by (\ref{boundaryv8}), for $\epsilon>0$ sufficiently small, we have 
$$\max_{(\boldsymbol \alpha,\xi)\in \mathrm B\setminus \mathrm B_\epsilon} V(\boldsymbol \alpha,\xi)< v_8.$$  
Let 
$$\delta =\frac{1}{2}\Big(v_8-\max_{(\boldsymbol \alpha,\xi)\in \mathrm B\setminus \mathrm B_\epsilon} V(\boldsymbol \alpha,\xi)\Big)>0.$$
 If $\frac{2\pi (\boldsymbol a+\boldsymbol 1)}{r} \in \mathrm H\setminus\mathrm H_{\epsilon},$ then $\big(\frac{2\pi (\boldsymbol a+\boldsymbol 1)}{r}, \frac{2\pi k}{r}\big)\in \mathrm B\setminus \mathrm B_\epsilon$ for any $k\in [\max\{T_i\},\min\{Q_j\}]\cap\mathbb Z;$ and by (1), we have
$$\Big|e^{\frac{r}{4\pi\mathbf i}U_r(\boldsymbol \alpha, \xi)}\Big|<O\Big(e^{\frac{r}{2\pi}\big(\max_{(\boldsymbol\alpha,\xi)\in\mathrm B\setminus\mathrm{B}_\epsilon} V(\boldsymbol \alpha,\xi)+\delta\big)}\Big)<O\Big(e^{\frac{r}{2\pi} \big(v_8-\delta\big)}\Big).$$
Therefore, by (\ref{3.9}), we have
      $$  \Bigg|\bigg\|\begin{matrix} 
      a_1 & a_2 & a_3 \\
      a_4 & a_5 & a_6
 \end{matrix}\bigg\|\Bigg|\leqslant\frac{2}{r}\sin\frac{2\pi}{r} \sum_k  \Big|e^{\frac{r}{4\pi\mathbf i}U_r(\boldsymbol \alpha, \xi)}\Big|< O\Big(e^{\frac{r}{2\pi}\big(v_8-\delta_3\big)}\Big)$$
  for any $\delta_3<\delta,$
where the sum is over all the integers in between $\max\{T_i\}$ and $\min\{Q_j\}.$
\end{proof}

\begin{proof}[Proof of Proposition \ref{error6j}]
By (\ref{VV}) and (\ref{nbhdvol}), let 
$$\delta' = \frac{1}{2}\min_{\boldsymbol \alpha \in \overline{\mathrm H_{\epsilon_0}}}\Big( \mathrm{Vol}(\Delta_{\boldsymbol \alpha})-\max_{\xi\in (\max\{\tau_i\},\min\{\eta_j\})\setminus [\xi_{\boldsymbol \alpha}-\epsilon, \xi_{\boldsymbol \alpha}+\epsilon]}V(\boldsymbol \alpha,\xi) \Big)>0.$$
Then by Proposition \ref{EST} (2) and the relation (\ref{LL}), for any $\delta_5<\delta',$ we have
$$\Big|e^{\frac{r}{4\pi\mathbf i}U_r(\boldsymbol a, k)}\Big|<O\Big(e^{\frac{r}{2\pi}\big(\max_{\xi\in (\max\{\tau_i\},\min\{\eta_j\})\setminus [\xi_{\boldsymbol \alpha}-\epsilon, \xi_{\boldsymbol \alpha}+\epsilon]}V(\boldsymbol \alpha,\xi)+\delta_5\big)}\Big)=O\Big(e^{\frac{r}{2\pi} \big(\mathrm{Vol}(\Delta_{\boldsymbol\alpha})-\delta_5\big)}\Big)$$
for each integer $k$ so that  the corresponding $\xi\in (\max\{\tau_i\},\min\{\eta_j\})\setminus [\xi_{\boldsymbol \alpha}-\epsilon,\xi_{\boldsymbol \alpha}+\epsilon].$
\end{proof}

\begin{proof}[Proof of Proposition \ref{other6j}]
By the Integration by Parts as on \cite[Page 703]{O1}, one has
$$\mathrm I_{\boldsymbol a}(n)=\frac{1}{n^2}\int_\mathbb R \Psi(\xi)e^{\frac{r}{4\pi \mathbf i}\big(U_r(\boldsymbol{\alpha},{\xi})-4\pi n\xi\big)}d{\xi},$$
where $\Psi(\xi)$ is a polynomial of the partial derivatives of at most the second order of  $\psi( \xi)$ and $U_r,$ and it  suffices to bound each $U_r(\boldsymbol \alpha,\xi)-4\pi n\xi.$

By (\ref{criticalequation}) and (\ref{convexU}), there exists an $\epsilon'>0$ sufficiently small such that for any $\boldsymbol \alpha \in \overline{\mathrm H_{\epsilon_0}}$ and for any $\xi\in\mathbb C$ with $\mathrm{Re}\xi \in [\xi_{\boldsymbol \alpha}-\epsilon,\xi_{\boldsymbol \alpha}+\epsilon]$ and $|\mathrm{Im}\xi|<\epsilon',$ we have 
\begin{equation}\label{<2pi}
\bigg|\frac{\partial \mathrm{Im}U(\boldsymbol \alpha,\xi)}{\partial \mathrm{Im}\xi}\bigg|<2\pi,
\end{equation}
and
\begin{equation}\label{<0}
\frac{\partial^2 \mathrm{Im}U(\boldsymbol \alpha,\xi)}{\partial \mathrm{Re}\xi^2}<0.
\end{equation}

Now, in the case that $n>0,$ we let 
\begin{equation*}
\begin{split}
\mathrm C_\text{top}= &\Big\{\xi + \epsilon'\mathbf i \ \Big|\  \xi \in [\xi_{\boldsymbol \alpha}-\epsilon,\xi_{\boldsymbol \alpha}+\epsilon]\Big\},\\
\mathrm C_\text{side}= & \Big\{\xi' + t \mathbf i\mathbf   \ \Big|\ \xi'=\xi_{\boldsymbol \alpha}\pm \epsilon, t\in [0,\epsilon']\Big\},\\
\mathrm C_\text{ext}= &\mathbb R \setminus [\xi_{\boldsymbol \alpha}-\epsilon,\xi_{\boldsymbol \alpha}+\epsilon],
\end{split}
\end{equation*}
and let 
$$\mathrm C=\mathrm C_\text{top}\cup\mathrm C_\text{side}\cup\mathrm C_\text{ext}.$$

For $\xi\in \mathrm C_\text{top},$ by the Mean Value Theorem, (\ref{<2pi}) and (\ref{320}), we have
\begin{equation*}
\begin{split}
\mathrm{Im}  U(\boldsymbol \alpha, \xi)\leqslant &\mathrm{Im}  U(\boldsymbol \alpha, \mathrm{Re}\xi)+\bigg|\frac{\partial \mathrm{Im}U(\boldsymbol \alpha,\xi)}{\partial \mathrm{Im}\xi}\Big|_{\xi=\xi'}\bigg|\cdot  |\mathrm{Im}\xi'|\\
<&\mathrm{Im}  U(\boldsymbol \alpha, \xi_{\boldsymbol\alpha})+ 2\pi \epsilon',  
\end{split}    
\end{equation*}
where $\xi'$ is some point on the line segment connecting $\xi$ and $\mathrm{Re}\xi.$ Then as $n\geqslant 1,$
\begin{equation}\label{top2}
\begin{split}
\mathrm{Im} \Big( U(\boldsymbol \alpha, \xi)-4\pi n\xi\Big)<&\mathrm{Im}  U(\boldsymbol \alpha, \xi_{\boldsymbol\alpha})+ 2\pi \epsilon'-4\pi n\epsilon'\\
\leqslant & \mathrm{Im}  U(\boldsymbol \alpha, \xi_{\boldsymbol\alpha})- 2\pi \epsilon'\\
=& 2\mathrm{Vol}(\Delta_{\boldsymbol\alpha})-2\pi \epsilon'. 
\end{split}    
\end{equation}

For $\xi \in \mathrm{C}_\text{side},$ by (\ref{UV}), (\ref{VV}) and (\ref{nbhdvol}), let $$\delta'=\frac{1}{2}\min_{\boldsymbol\alpha\in \overline {\mathrm H_{\epsilon_0}}}\Big(2\mathrm{Vol}(\Delta_{\boldsymbol\alpha})
-\max_{\xi=\xi_{\boldsymbol\alpha}\pm \epsilon} \mathrm{Im} U(\boldsymbol \alpha,\xi)\Big)>0.$$
Then by (\ref{nbhdvol}),  for $\xi'=\xi_{\boldsymbol\alpha}\pm \epsilon,$ we have
\begin{equation}\label{sidebottom2}
 \mathrm{Im}\Big(U(\boldsymbol\alpha,\xi')-4\pi n \xi'\big) =  \mathrm{Im}U(\boldsymbol\alpha,\xi')
< 2\mathrm{Vol}(\Delta_{\boldsymbol\alpha})-2\delta'.
\end{equation}
By (\ref{convexU}) and the harmonicity of $\mathrm{Im}U,$ we have that $\mathrm{Im}U(\boldsymbol \alpha, \xi'+t  \mathbf i)$ is strictly concave up in $t;$ and by (\ref{top2}) and  (\ref{sidebottom2}), we have
\begin{equation}\label{side2}
   \mathrm{Im} \Big( U(\boldsymbol \alpha, \xi)-4\pi n\xi\Big) < 2\mathrm{Vol}(\Delta_{\boldsymbol\alpha})-2\min\{\pi\epsilon',\delta'\}.
\end{equation}

 For $\xi\in \mathrm C_\text{ext},$ by (\ref{nbhdvol}),  we have
\begin{equation}\label{ext2}
 \mathrm{Im}\Big(U(\boldsymbol\alpha,\xi)-4\pi n \xi\big) =  \mathrm{Im}U(\boldsymbol\alpha,\xi)
< 2\mathrm{Vol}(\Delta_{\boldsymbol\alpha})-2\delta'.
\end{equation}

Finally, let
$$\delta_6=\min\{\pi \epsilon', \delta'\};$$ 
and notice that, as both $\epsilon'$ and $\delta'$ are independent of $n,$ $\delta_6$ is so as well.
Putting (\ref{top2}), (\ref{side2}) and (\ref{ext2}) together,  we have 
\begin{equation*}
  \mathrm{Im}\Big(U(\boldsymbol\alpha,\xi)-4\pi n \xi\big)
< 2\mathrm{Vol}(\Delta_{\boldsymbol\alpha})-2\delta_6.
\end{equation*}
on $\mathrm C.$ Then the result follows from Proposition \ref{EST} (1).
\\

In the case that $n<0,$ let 
$$\mathrm C_\text{top}= \Big\{\xi'- \epsilon\mathbf i \ \Big|\  \xi \in [\xi_{\boldsymbol \alpha}-\epsilon,\xi_{\boldsymbol \alpha}+\epsilon]\Big\}\quad\text{and}\quad \mathrm C_\text{side}= \Big\{\xi' - t \mathbf i\mathbf   \ \Big|\ \xi'=\xi_{\boldsymbol \alpha}\pm \epsilon, t\in [0,\epsilon']\Big\},$$
and the rest of the proof follows verbatim that of the previous case. 
\end{proof}


\subsection{The challenging part}\label{challenging}

The goal of  this section is to prove Proposition \ref{leading6j}. The strategy is as follows. We will show that the $0$-th Fourier coefficient of the extended quantum $6j$-symbols and the $b$-$6j$ symbols recalled in Section \ref{b6jsymbols} satisfy the same real Gevrey-$1$ estimate (see Proposition \ref{r6jb6j}), hence share the same reflection symmetry after the truncated Borel-Laplace transforms (see Proposition \ref{refI0}). This will induce the desired weak reflection symmetry for the $0$-th Fourier coefficient of the extended quantum $6j$-symbols. The rest of this section is organized as follows. In Subsection \ref{621}, we collect the basic properties of real Gevrey-$1$ estimates and truncated Borel-Laplace transforms, with the proofs postponed to Appendix \ref{Gevrey-1}. In Subsection \ref{622}, we establish a real Gevrey-$1$ estimate for the double sine function; and in Subsection \ref{623} we use this estimate to derive the corresponding real  Gevrey-$1$ estimates of the  $0$-th Fourier coefficient of the extended quantum $6j$-symbols and the $b$-$6j$ symbols. Finally, we complete the proof of  Proposition \ref{leading6j} in Subsection \ref{624}.

\subsubsection{Real Gevrey-$1$ estimate and truncated Borel-Laplace transform}\label{621}

For $x_0>0,$ let $f$ be a complex valued function  on $(0,x_0).$ Let 
$$\widehat f(x)=\sum_{n=0}^\infty a_nx^n$$
be a formal power series in $x.$  Here we do not assume the continuity of $f,$ nor the convergence of $\widehat f.$ 
For $T>0,$ and let $g$ be a complex valued integrable function on $[0,T].$

\begin{definition}
\begin{enumerate}[(1)]
    \item 
    If there exist constants $C>0$ and $\sigma>0$ such that
\begin{equation}\label{gevrey-1}
\bigg| f(x)-\sum_{n=0}^{N-1}a_nx^n  \bigg|\leqslant C\sigma^NN!|x|^N    
\end{equation}
 for each  $N\geqslant 0$ and $x\in (0,x_0),$ then  we say that the function 
 $f$ satisfies a \emph{$(C,\sigma)$ real Gevrey-$1$ estimate}, and call the formal power series  $\widehat f$ a \emph{Gevrey-$1$ expansion} of $f$

 \item The \emph{Borel transform} $\mathcal B\widehat f$ of the formal power series 
 $\widehat f$ is the following formal power series 
 $$\mathcal B \widehat f (z) = \sum_{n=0}^\infty \frac{a_n}{n!}z^n.$$

 \item   The \emph{truncated Laplace transform} of $g$ truncated at $T$ is the function $\mathcal L_Tg$ defined  on $(0,\infty)$ by
 $$\mathcal L_Tg(x) = \frac{1}{x}\int_0^T e^{-\frac{t}{x}}g(t)dt.$$   
\end{enumerate}  
\end{definition}

\begin{lemma}\label{GBL}
Let $x_0>0,$ and let $f$ be a complex valued function on $(0,x_0).$ Suppose $f$ satisfies a $(C,\sigma)$ real Gevrey - $1$ estimate for some $C>0$ and $\sigma>0.$  Then:
\begin{enumerate}[(1)]
\item For each $n\geqslant 0,$
\begin{equation}\label{anbound}
|a_n|\leqslant C\sigma^nn!,    
\end{equation}
and the Borel transform $\mathcal B\widehat f (z)$ of $\widehat f$ defines a holomorphic function on 
$\big\{ z\in \mathbb C\ \big|\ |z| < \sigma^{-1}\big\}.$

\item For each $T\in (0,\sigma^{-1}),$ there exists a $D>0$ depending only on $C,$ $\sigma$ and $T$ such that for each $x\in (0, \min\{x_0,T\}),$
$$\Big| f(x)-\mathcal L_T\mathcal B \widehat f(x)\Big| \leqslant D\sqrt{\frac{T}{x}}e^{-\frac{T}{x}}.$$
\end{enumerate}
\end{lemma}

\begin{lemma}\label{Exp}
Let $x_0>0,$ and let $f$ be a complex valued function on $(0,x_0).$ Suppose $f$ satisfies a $(C,\sigma)$ real Gevrey - $1$ estimate for some $C>0$ and $\sigma>0.$ Then $e^{f(x)}$ satisfies a $(C_1,\sigma_1)$ real Gevrey - $1$ estimate
$$\bigg|e^{f(x)}-\sum_{n=0}^{N-1}b_nx^n\bigg|\leqslant C_1\sigma_1^NN!
|x|^N$$
 for some $C_1>0$ depending only on $C,$ $\sigma$ and $x_0$ and some $\sigma_1$ depending only on $\sigma.$ Moreover, if 
$$\widehat f(x)=\sum_{n=0}^\infty a_nx^n$$
is a Gevrey - $1$ expansion of $f,$ then for each $n\geqslant 0,$ 
$$b_n=e^{a_0}\sum_{k=0}^n\sum_{\substack{
m_1,\ldots,m_n\geqslant0\\
m_1+\cdots+m_n=k\\
m_1+2m_2+\cdots+nm_n=n}}
\frac{1}{m_1!\cdots m_n!}
\prod_{j=1}^{n}a_j^{m_j}.$$
\end{lemma}

\begin{lemma}\label{Int}
Let $\rho>0$ and $\rho_0\in (0,\rho),$ and let  $g(u)$ be a holomorphic function defined for $u\in \mathbb C$ with $|u|<\rho$ satisfying 
$$\inf_{u\in [-\rho_0,\rho_0]}\mathrm{Re}g(u)\geqslant c\quad\text{and}\quad \sup_{|u|<\rho}|g(u)|\leqslant G$$
for some $c>0$ and $G>0.$ Let $x_0>0,$ 
and let $f(u,x)$ be a function defined for $u$ with $|u|<\rho$ and  $x\in (0,x_0).$ Suppose that $f(u,x)$ is continuous in $u\in [-\rho_0,\rho_0],$ and satisfies a $(C,\sigma)$ real Gevrey-$1$ estimate 
$$\bigg|f(u,x)-\sum_{n=0}^{N-1}a_n(u)x^n\bigg| \leqslant C\sigma^NN!x^N$$
for each $u$ with $|u|<\rho$ and for some $C>0$ and  $\sigma>0;$  and suppose for each $n\geqslant 0$ that $a_n(u)$ is a holomorphic function on the disk 
$|u|<\rho.$
 Then the integral
$$I(x)=\frac{1}{\sqrt {2\pi x}}\int_{-\rho_0}^{\rho_0} e^{-\frac{g(u)u^2}{x}}f(u,x)du$$
satisfies a $(C_2,\sigma_2)$ real Gevrey-$1$ estimate
$$\bigg|I(x)-\sum\limits_{n=0}^{N-1}c_nx^n\bigg|\leqslant C_2\sigma_2^N N!x^N$$
for some $C_2>0$ depending only on $C,$ $c$ and $G,$ and some $\sigma_2>0$ depending only on $\sigma,$ $c,$ $G,$ $\rho_0$ and $\rho,$
and 
\begin{equation*}
 c_n=\frac1{\sqrt2}\sum\limits_{\substack{j,k\geqslant 0\\j+k=n}}
\frac{F_j^{(2k)}(0)}{4^k k!},
\end{equation*}
where for each $j\geqslant 0,$ $F_j$ is a holomorphic function defined in a neighborhood of $0$ such that if $v=u\sqrt{g(u)},$ then 
$$F_j(v)=\frac{a_j(u)}{\big(u\sqrt{g(u)}\big)'}.$$
\end{lemma}

\subsubsection{Gevrey-$1$ estimate of the double sine function}\label{622}

To have a uniform treatment of the extended quantum $6j$-symbols and the $b$-$6j$ symbols, we will consider the  function $s_\hbar$ defined  for $\hbar \in (-1,1)$ and  $z\in \mathbb C$ with $0<\mathrm{Re}z<\pi$ as follows. Let $\Omega$ be the contour going along the real line and passing above the origin. For $\hbar \neq 0,$ we let
$$s_\hbar(z)\doteq \pi \hbar\int_\Omega\frac{\sinh \big((\pi+\pi \hbar-2z)x\big)}{4x\sinh(\pi \hbar x)\sinh(\pi x)}dx;$$
and let 
$$s_0(z)\doteq \int_\Omega\frac{\sinh \big((\pi-2z)x\big)}{4x^2\sinh(\pi x)}dx.$$
Directly from the definition, we have 
$$s_\hbar\Big(z+\frac{\pi\hbar}{2}\Big)=s_{-\hbar}\Big(z-\frac{\pi\hbar}{2}\Big).$$

Let 
\begin{equation*}
S_b(z)=\exp\Bigg(\int_\Omega\frac{\sinh\Big(\big(\frac{Q}{2}-z\big)t\Big)}{4t\sinh(\frac{bt}{2})\sinh(\frac{t}{2b})}dt\Bigg) 
\end{equation*}
be the double sine function defined in (\ref{eq:def-S}), and let
$$L(z)=z^2-\pi z +\frac{\pi^2}{6}-\mathrm{Li}_2\big(e^{2\mathbf iz}\big)$$
be the complexified Lobachevsky function defined in (\ref{eq:Lx}). Then for $\hbar>0$ and $b=\sqrt{\hbar},$
\begin{equation}\label{sh}
2\pi \mathbf i b^2\log S_b\Big(\frac{z}{\pi b}+\frac{b}{2}\Big)=2\mathbf i  s_\hbar\Big(z+\frac{\pi \hbar}{2}\Big)
=2\mathbf i  s_{-\hbar}\Big(z-\frac{\pi \hbar}{2}\Big);
\end{equation}
and 
\begin{equation}\label{s0}
L(z)=2\mathbf i s_0(z).    
\end{equation}

\begin{proposition}\label{shGevrey}
Let $M>0,$ $\Lambda>0$ and $\epsilon \in (0,\frac{\pi}{2}).$  Suppose that 
$$ \epsilon \leqslant \mathrm{Re}z\leqslant \pi-\epsilon, \quad \big|\mathrm{Im}z\big|\leqslant M \quad\text{and}\quad |\lambda| \leqslant \Lambda.$$
Then there exists a $C>0$ depending only on $M$ and $\Lambda,$ and a $\sigma>0$ depending only on $\epsilon,$ such that for every integer $N\geqslant 1$ and real number $\hbar$ with $0<|\hbar|<\frac{\epsilon}{2\pi(\Lambda +1)},$
\begin{equation}\label{shg}
\bigg|s_\hbar\big(z+\lambda\pi\hbar\big)-\frac{1}{2\mathbf i}\sum_{n=0}^{N-1}\frac{\pi^n}{n!} B_n(\lambda)L^{(n)}(z)\hbar^n\bigg|\leqslant C\sigma^NN!|\hbar|^N, 
\end{equation}
where $B_n$ is the $n$-th Bernoulli polynomial characterized by 
$$\frac{te^{\lambda t}}{e^t-1}=\sum_{n=0}^\infty \frac{B_n(\lambda)}{n!}t^n$$
for $t\in \mathbb C$ with $|t|< 2\pi.$
\end{proposition}

\begin{proof}
Let 
 $$A_\Lambda=\max_{|t|=4,\lambda\in [-\Lambda,\Lambda]}\bigg|\frac{te^{\lambda t}}{e^t-1}\bigg|,$$
 and let 
$$C_\Lambda=2A_\Lambda + \frac{2}{1-e^{-2}}.$$

We first show that for all $t\in \mathbb R$ and $N\geqslant 1,$
\begin{equation}\label{611}
\bigg| \frac{te^{\lambda t}}{e^t-1} - \sum_{n=0}^{N-1} \frac{B_n(\lambda)}{n!}t^n  \bigg| \leqslant C_\Lambda e^{\Lambda|t|} \bigg(\frac{|t|}{2}\bigg)^N.    
\end{equation}
To see this, by Cauchy's Estimate, we have 
\begin{equation}\label{BA}
\frac{|B_n(\lambda)|}{n!}\leqslant \frac{A_\Lambda}{4^n}.    
\end{equation}
As a consequence, for $t\in\mathbb C$ with $|t|<2,$ we have
\begin{equation}\label{t<2}
\bigg| \frac{te^{\lambda t}}{e^t-1} - \sum_{n=0}^{N-1} \frac{B_n(\lambda)}{n!}t^n  \bigg| \leqslant A_\Lambda \sum_{n=N}^\infty \bigg(\frac{|t|}{4}\bigg)^n\leqslant 2A_\Lambda e^{\Lambda |t|}\bigg(\frac{|t|}{2}\bigg)^N.    
\end{equation}
Now for  $t\in \mathbb R$ with $|t|\geqslant 2,$ we have 
 $$\bigg|\frac{te^{\lambda t}}{e^t-1}\bigg|\leqslant \frac{|t|e^{\Lambda |t|}}{1-e^{-2}}\leqslant\frac{2e^{\Lambda |t|}}{1-e^{-2}} \bigg(\frac{|t|}{2}\bigg)^N;$$
 and by (\ref{BA}), we have 
 \begin{equation}\label{t>2}
 \bigg|\sum_{n=0}^{N-1}\frac{B_n(\lambda)}{n!}t^n\bigg|\leqslant A_\Lambda \sum_{n=0}^{N-1}2^{-n}\bigg(\frac{|t|}{2}\bigg)^n\leqslant  2A_\Lambda e^{\Lambda |t|}\bigg(\frac{|t|}{2}\bigg)^N.
 \end{equation}
Then (\ref{611}) follows from (\ref{t<2}) and (\ref{t>2}).
\medskip

Next, to simplify the notation, let the remainder 
$$R_N(t,\lambda)= \frac{te^{\lambda t}}{e^t-1} - \sum_{n=0}^{N-1}\frac{B_n(\lambda)}{n!}t^n.$$
Then we observe that for $N\geqslant 3,$
\begin{equation}\label{616}
\begin{split}
&s_\hbar\big(z+\lambda\pi\hbar\big)-\frac{1}{2\mathbf i}\sum_{n=0}^{N-1}\frac{\pi^n}{n!} B_n(\lambda)L^{(n)}(z)\hbar^n\\
=& \frac{1}{8}\int_\Omega \frac{ e^{(\pi -2z)x}R_N(-2\pi\hbar x,\lambda)-e^{-(\pi -2z)x}R_N(2\pi \hbar x,\lambda)}{x^2\sinh(\pi x)}dx\\
=& \frac{1}{8}\int_{\mathbb R}\frac{ e^{(\pi -2z)x}R_N(-2\pi\hbar x,\lambda)-e^{-(\pi -2z)x}R_N(2\pi \hbar x,\lambda)}{x^2\sinh(\pi x)}dx,
\end{split}    
\end{equation}
where the first equality comes from $\sinh\big((\pi-2z)x\big)=\frac{1}{2}\big(e^{(\pi -2z)x}-e^{-(\pi -2z)x}\big),$ (\ref{s0}) and its consequence  that 
$$\frac{L^{(n)}(z)}{2\mathbf i}=\frac{1}{8}\int_\Omega \frac{(-2x)^ne^{(\pi-2z)x}-(2x)^ne^{-(\pi-2z)x}}{x^2\sinh(\pi x)}dx$$
for $n\geqslant 0,$
and the second equality comes from (\ref{611}) that the integrand is of order $O\big(x^{N-3}\big)$ hence $0$ is a removable singularity of the integrand. For $x>0,$ we have
$$\frac{\big|e^{(\pi-2z)x}\big|+\big|e^{-(\pi-2z)x}\big|}{\big|\sinh(\pi x)\big|}\leqslant \frac{4e^{-2\epsilon x}}{1-e^{-2\pi x}}\leqslant 4e^{-2\epsilon x}\bigg(1+\frac{1}{2\pi x}\bigg).$$
Together with (\ref{611}) and (\ref{616}), we have 
\begin{equation*}
\begin{split}
\bigg| s_\hbar\big(z+\lambda\pi\hbar\big)-\frac{1}{2\mathbf i}\sum_{n=0}^{N-1}\frac{\pi^n}{n!} B_n(\lambda)L^{(n)}(z)\hbar^n \bigg|\leqslant & C_\Lambda \pi^N |\hbar|^N \int_0^\infty x^{N-2}e^{2(\pi|\hbar|\Lambda-\epsilon)x}\bigg(1+\frac{1}{2\pi x}\bigg)dx\\
\leqslant & C_\Lambda \pi^N |\hbar|^N \int_0^\infty x^{N-2}e^{-\epsilon x}\bigg(1+\frac{1}{2\pi x}\bigg)dx\\
= & C_\Lambda \pi^N |\hbar|^N \Big((N-2)!\epsilon^{1-N}+\frac{1}{2\pi}(N-3)!\epsilon^{2-N}\Big)\\
\leqslant & C_\Lambda\pi^N\epsilon^{-N}N!|\hbar|^N,
\end{split}
\end{equation*}
where the second inequality comes from that $|\hbar|< \frac{\epsilon}{2\pi(\Lambda+1)}<\frac{\epsilon}{2\pi\Lambda},$ and the last inequality comes from that $\epsilon< \frac{\pi}{2}$ so that $\epsilon<\frac{(N-1)N}{2}$ and $\frac{\epsilon^2}{2\pi}<\frac{(N-2)(N-1)N}{2}.$ 
Then (\ref{shg}) holds for $N\geqslant 3$ with $C=C_\Lambda$ and $\sigma=\pi\epsilon^{-1}.$

For $N=2,$ we have 
\begin{equation*}
\begin{split}
&\bigg| s_\hbar\big(z+\lambda\pi\hbar\big)-\frac{1}{2\mathbf i}\sum_{n=0}^{1}\frac{\pi^n}{n!} B_n(\lambda)L^{(n)}(z)\hbar^n \bigg|\\
\leqslant & \bigg| s_\hbar\big(z+\lambda\pi\hbar\big)-\frac{1}{2\mathbf i}\sum_{n=0}^{2}\frac{\pi^n}{n!} B_n(\lambda)L^{(n)}(z)\hbar^n \bigg| + \frac{\pi^2}{2}\Big|B_2(\lambda)\Big| \big|\cot z\big||\hbar|^2\\
\leqslant & C_\Lambda\pi^3\epsilon^{-3}3!|\hbar|^3+ \frac{\pi^2}{2}\Big|B_2(\lambda)\Big| \big|\cot z\big||\hbar|^2.
\end{split}    
\end{equation*}
Since $|\mathrm{Im}z|\leqslant M,$ there is a $C_M>0$ such that 
$$\big|\cot z\big|\leqslant C_M\epsilon^{-2};$$
and since $B_2(\lambda)$ is bounded on $[-\Lambda,\Lambda]$ and $\frac{|\hbar|}{\epsilon}<\frac{1}{2\pi (\Lambda+1)}<1,$ we have (\ref{shg}) for $N=2,$ by increasing $C=C_\Lambda$ if necessary and with $\sigma=\pi\epsilon^{-1}.$

For $N=1,$ we have 
\begin{equation*}
\begin{split}
&\bigg| s_\hbar\big(z+\lambda\pi\hbar\big)-\frac{1}{2\mathbf i} B_0(\lambda)L(z) \bigg|\\
\leqslant & \bigg| s_\hbar\big(z+\lambda\pi\hbar\big)-\frac{1}{2\mathbf i}\sum_{n=0}^{1}\frac{\pi^n}{n!} B_n(\lambda)L^{(n)}(z)\hbar^n \bigg| + \pi\Big|B_1(\lambda)\Big| \big|\log \big(2\sin z)\big||\hbar|\\
\leqslant & C\pi^2\epsilon^{-2}2!|\hbar|^2 + \pi\Big|B_1(\lambda)\Big| \big|\log \big(2\sin z)\big||\hbar|.
\end{split}    
\end{equation*}
Since $|\mathrm{Im}z|\leqslant M,$ there is a $C_M>0$ such that 
$$\big|\log\big(2\sin z\big)\big|\leqslant C_M\epsilon^{-1};$$
and since $B_1(\lambda)$ is bounded on $[-\Lambda,\Lambda]$ and  $\frac{|\hbar|}{\epsilon}<\frac{1}{2\pi (\Lambda+1)}<1,$ we have (\ref{shg}) for $N=1,$ by increasing $C$ further if necessary and with $\sigma=\pi\epsilon^{-1}.$ 
\end{proof}

\subsubsection{Extended quantum $6j$-symbols and $b$-$6j$symbols}\label{623}

The goal of this subsection is to understand the relationship between the $0$-th Fourier coefficient of the extended quantum $6j$-symbols and the $b$-$6j$ symbols.

For $\boldsymbol \alpha = (\alpha_1,\dots,\alpha_6)\in \mathrm H,$ $\xi\in (\max\{\tau_i\},\min\{\eta_j\}),$ and $h\in (-1,1)$ with $|\hbar|$ sufficiently small, define the function
\begin{equation*}
 V_\hbar (\boldsymbol \alpha,\xi)\doteq \frac{1}{2}\sum_{i=1}^4\sum_{j=1}^4s_\hbar\bigg(\eta_j-\tau_i+\frac{\pi \hbar}{2}\bigg)-\sum_{i=1}^4s_\hbar\bigg(\xi-\tau_i+\frac{5\pi\hbar}{2}\bigg)-\sum_{j=1}^4s_\hbar\big(\eta_j-\xi-\pi\hbar\big). 
\end{equation*}

Let $\epsilon_0>0$ be as in Proposition \ref{key'} and Theorem \ref{wrs}, and let $\boldsymbol a$ by an $r$-hyperideal six-tuple with the corresponding $\boldsymbol\alpha\in\mathrm H_{\epsilon_0}.$
Then by (\ref{Fn}), (\ref{Ur}) and (\ref{sh}), for $\hbar=b^2=\frac{2}{r},$ we have that the $0$-th Fourier coefficient of the extended quantum $6j$-symbol
\begin{equation}\label{I0}
\mathrm I_{\boldsymbol a}(0)= \int_\mathbb R\psi(\xi)\exp\bigg(\frac{V_\hbar(\alpha_1,\dots,\alpha_6,\xi)}{\pi \hbar }\bigg)d\xi.
\end{equation}

For $b\in (0,1),$ let $(a_1,\dots,a_6)\in\mathbb R^6$ be a $b$-admissible six-tuple, and recall that for each $k\in\{1,\dots,6\},$ $\alpha_k=2\pi b a_k-\pi b^2.$
Then by (\ref{6jint}) and (\ref{sh}), for $\hbar=b^2,$ we have that  the $b$-$6j$ symbol
\begin{equation}\label{b6js}
\bigg\{\begin{matrix} a_1 & a_2 & a_3 \\ a_4 & a_5 & a_6 \end{matrix} \bigg\}_b=\frac{1}{\pi\mathbf i b}\int_{\Gamma}\exp\bigg(\frac{V_{-\hbar}(\alpha_1,\dots,\alpha_6,\xi)}{-\pi \hbar} \bigg)d \xi,
\end{equation}
where $\Gamma$ is a vertical straight line that intersects the real axis in the interval  $\big(\max\{\tau_i\}+\frac{3\pi\hbar}{2},\min\{\eta_j\}+2\pi\hbar\big)$. 

Recall that $\epsilon>0$ is sufficiently small so that for each $\boldsymbol \alpha \in \mathrm H_{\epsilon_0},$ the interval $[\xi_{\boldsymbol \alpha}-3\epsilon,\xi_{\boldsymbol \alpha}+3\epsilon]$ entirely lies in $(\max\{\tau_i\},\min\{\eta_j\});$ and also recall that   $\Delta_{\boldsymbol \alpha}$ is the truncated hyperideal tetrahedron with dihedral angles $|\alpha_1-\pi|,\dots, |\alpha_6-\pi|$ with volume $\mathrm{Vol}(\Delta_{\boldsymbol \alpha}).$

 We let 
\begin{equation}\label{Ir}
\mathcal I^r_{\boldsymbol \alpha}( \hbar)=\frac{e^{-\frac{\mathrm{Vol}(\Delta_{\boldsymbol \alpha})}{\pi\hbar}}}{\sqrt{2\pi\hbar}}\int_{\xi_{\boldsymbol\alpha}-\epsilon}^{\xi_{\boldsymbol\alpha}+\epsilon}    \exp\bigg(\frac{V_\hbar(\alpha_1,\dots,\alpha_6,\xi)}{\pi \hbar }\bigg)d\xi;
\end{equation}
and let 
\begin{equation}\label{Ib}
\mathcal I^b_{\boldsymbol \alpha}( \hbar)=\frac{e^{\frac{\mathrm{Vol}(\Delta_{\boldsymbol \alpha})}{\pi\hbar}}}{\mathbf i\sqrt{2\pi\hbar}}\int_{\xi_{\boldsymbol\alpha}-\mathbf i\epsilon}^{\xi_{\boldsymbol\alpha}+\mathbf i\epsilon}    \exp\bigg(\frac{V_{-\hbar}(\alpha_1,\dots,\alpha_6,\xi)}{-\pi \hbar }\bigg)d\xi,
\end{equation}
where the integral is along the vertical line connecting $\xi_{\boldsymbol\alpha}\pm\mathbf i\epsilon.$

\begin{proposition} \label{r6jb6j}
There exist functions $c_n,d_n: \mathrm H_{\epsilon_0}\to\mathbb R$ for all integers $n\geqslant 0,$ and constants  $C>0,$ $\sigma>0$ and $\hbar_0>0$ depending only on $\epsilon,$ such that the following properties hold:  
\begin{enumerate}[(1)]
    \item  For all $\boldsymbol\alpha \in\mathrm H_{\epsilon_0},$ $N\geqslant 0$ and $\hbar\in(0,\hbar_0),$
    $$\bigg|\mathcal I^r_{\boldsymbol \alpha}( \hbar)-\sum_{n=0}^{N-1}c_n(\boldsymbol\alpha)\hbar^n\bigg| \leqslant C\sigma^NN!\hbar^N$$
    and
 $$\bigg|\mathcal I^b_{\boldsymbol \alpha}( \hbar)-\sum_{n=0}^{N-1}d_n(\boldsymbol\alpha)\hbar^n\bigg| \leqslant C\sigma^NN!\hbar^N,$$
i.e., both $\mathcal I^r_{\boldsymbol \alpha}( \hbar)$ and $\mathcal I^b_{\boldsymbol \alpha}( \hbar)$ satisfy a $(C,\sigma)$ Gevrey-$1$ estimate. Moreover, 
\begin{equation*}
    d_n(\boldsymbol\alpha)=(-1)^nc_n(\boldsymbol\alpha)
\end{equation*}
for all $\boldsymbol\alpha\in\mathrm H_{\epsilon_0}$ and $n\geqslant 0.$

\item Consider the formal power series
$$\widehat{\mathcal I^r_{\boldsymbol\alpha}}(\hbar)=\sum_{n=0}^\infty c_n(\boldsymbol\alpha)\hbar^n\quad\text{and}\quad \widehat{\mathcal I^b_{\boldsymbol\alpha}}(\hbar)=\sum_{n=0}^\infty d_n(\boldsymbol\alpha)\hbar^n.$$
Then for all $\boldsymbol \alpha \in\mathrm H_{\epsilon_0},$
both of the Borel transforms $\mathcal B\widehat{\mathcal I^r_{\boldsymbol\alpha}}$ and $\mathcal B\widehat{\mathcal I^b_{\boldsymbol\alpha}}$ have radius of convergence at least $\sigma^{-1}.$ Moreover,
$$\mathcal B\widehat{\mathcal I^r_{\boldsymbol\alpha}}(t)=\mathcal B\widehat{\mathcal I^b_{\boldsymbol\alpha}}(-t)$$
for all $t\in\mathbb C$ with $|t|<\sigma^{-1}.$

\item For each $T\in (0,\sigma^{-1}),$ there exists a $D>0$ depending only on $\epsilon$ and $T$ such that for each $\hbar\in (0, \min\{\hbar_0,T\}),$
$$\Big| \mathcal I^r_{\boldsymbol \alpha}( \hbar)-\mathcal L_T\mathcal B \widehat {\mathcal I^r_{\boldsymbol \alpha}}( \hbar)\Big| \leqslant D\sqrt{\frac{T}{\hbar}}e^{-\frac{T}{\hbar}}.$$

\end{enumerate}    
\end{proposition}

\begin{proof} By Lemma \ref{GBL}, (2) and (3) follow directly from (1).
\medskip

To see (1), for $\boldsymbol\alpha\in\mathrm H_{\epsilon_0},$ let $s=\xi-\xi_{\boldsymbol\alpha}.$ Then by Proposition \ref{shGevrey} with $M=2\epsilon$ and $\Lambda=3,$ there exist $C_0>0$ and $\sigma_0>0$ depending only on $\epsilon$ such that for each $N\geqslant 1,$ $\hbar\in\mathbb R$ with $|\hbar| \in (0, \frac{\epsilon}{8\pi}),$ and $s\in\mathbb C$ with $|s|< 2\epsilon,$
$$\bigg|V_{\hbar}\big(\boldsymbol\alpha, \xi_{\boldsymbol \alpha}+s\big)-\sum_{n=0}^{N-1}v_n(\boldsymbol\alpha,\xi_{\boldsymbol \alpha}+s)\hbar^n\bigg|\leqslant C_0\sigma_0^NN!|\hbar|^N,$$
where for each $n\geqslant 0,$
\begin{equation*}
\begin{split}
&v_n(\boldsymbol \alpha, \xi)\\
=& \frac{\pi ^n}{2\mathbf i n!} \bigg(\frac{B_n\big(\frac{1}{2}\big)}{2}\sum_{i=1}^4\sum_{j=1}^4L^{(n)}(\eta_j-\tau_i)-B_n\Big(\frac{5}{2}\Big)\sum_{i=1}^4L^{(n)}(\xi-\tau_i)-B_n(-1)\sum_{j=1}^4L^{(n)}(\eta_j-\xi) \bigg)
\end{split}
\end{equation*}
is holomorphic in $\xi.$
As a consequence, we have 
\begin{equation*}
\begin{split}
\bigg|\frac{V_{\hbar}\big(\boldsymbol\alpha, \xi_{\boldsymbol \alpha}+s\big)-V_0\big(\boldsymbol\alpha, \xi_{\boldsymbol \alpha}+s\big)}{\pi\hbar}-\frac{1}{\pi}\sum_{n=0}^{N-1}v_{n+1}(\boldsymbol\alpha,\xi_{\boldsymbol \alpha}+s)\hbar^n\bigg|\leqslant & \frac{C_0}{\pi}\sigma_0^{N+1}(N+1)!|\hbar|^N\\
\leqslant & \frac{C_0\sigma_0}{\pi}(2\sigma_0)^NN!|\hbar|^N,
\end{split}
\end{equation*}
where the last inequality comes from that  $N+1\leqslant 2^N$ for $N\geqslant 0.$
Then by Lemma \ref{Exp}, there exist $C_1>0$ and $\sigma_1>0$ depending only on $\epsilon$ such that 
for each $N\geqslant 0,$  $\hbar\in\mathbb R$ with $|\hbar| \in (0, \frac{\epsilon}{8\pi}),$  and $s\in\mathbb C$ with $|s|<2\epsilon,$
\begin{equation}\label{exp}
 \bigg|\exp\bigg({\frac{V_{\hbar}\big(\boldsymbol\alpha, \xi_{\boldsymbol \alpha}+s\big)-V_0\big(\boldsymbol\alpha, \xi_{\boldsymbol \alpha}+s\big)}{\pi\hbar}}\bigg)-\sum_{n=0}^{N-1}b_n(\boldsymbol\alpha,\xi_{\boldsymbol \alpha}+s)\hbar^n\bigg|\leqslant C_1\sigma_1^NN!|\hbar|^N,   
\end{equation}
where for  each $n\geqslant 0,$
\begin{equation*}
b_n(\boldsymbol \alpha, \xi)=e^{\frac{v_1(\boldsymbol\alpha,\xi)}{\pi}}\sum_{k=0}^n\sum_{\substack{
m_1,\ldots,m_n\geqslant0\\
m_1+\cdots+m_n=k\\
m_1+2m_2+\cdots+nm_n=n}}
\frac{1}{m_1!\cdots m_n!}
\prod_{j=1}^{n}\bigg(\frac{v_{j+1}(\boldsymbol\alpha,\xi)}{\pi}\bigg)^{m_j}
\end{equation*}
is holomorphic in $\xi.$

Next, for $\boldsymbol\alpha\in \mathrm H_{\epsilon_0}$ and $s\in\mathbb C$ with $|s|< 2\epsilon,$ we let
$$g_{\boldsymbol\alpha}(s)\doteq \frac{V_0(\boldsymbol\alpha,\xi_{\boldsymbol\alpha})-V_0(\boldsymbol\alpha,\xi_{\boldsymbol\alpha}+s)}{s^2}$$
for $s\neq 0,$ and let 
$$g_{\boldsymbol\alpha}(0)\doteq-\frac{1}{2}\frac{\partial ^2 V_0(\boldsymbol\alpha,\xi)}{\partial \xi^2}\Big|_{\xi=\xi_{\boldsymbol\alpha}}.$$
Since $\frac{\partial  V_0(\boldsymbol\alpha,\xi)}{\partial \xi}\Big|_{\xi=\xi_{\boldsymbol\alpha}}=0,$ $g_{\boldsymbol\alpha}$ is continuous at $0,$ and hence is holomorphic on the disk  $|s|<2\epsilon.$
Then by the compactness of the region $\big\{ (\boldsymbol \alpha, \xi)\in \mathbb R^6\times \mathbb C\ \big|\ \boldsymbol\alpha\in \overline{\mathrm H_{\epsilon_0}},|\xi-\xi_{\boldsymbol\alpha}|\leqslant 2\epsilon\big\},$ there exists a $G>0$ depending only on $\epsilon$ such that 
\begin{equation}\label{supG}
    \sup_{|s|\leqslant 2\epsilon} |g_{\boldsymbol\alpha}(s)|\leqslant \pi G
\end{equation}
for all $\boldsymbol\alpha\in \mathrm H_{\epsilon_0}.$ Also, by Taylor's Formula and that $\frac{\partial  V_0(\boldsymbol\alpha,\xi)}{\partial \xi}\big|_{\xi=\xi_{\boldsymbol\alpha}}=0,$  we have
\begin{equation}\label{Taylor}
g_{\boldsymbol \alpha}(s)=-\int_0^1(1-t)\frac{\partial ^2 V_0(\boldsymbol\alpha,\xi)}{\partial \xi^2}\Big|_{\xi=\xi_{\boldsymbol\alpha}+ts}dt.    
\end{equation}
By the fact that $\mathrm{Re}V_0(\boldsymbol \alpha,\xi)$ is strictly concave down in $\mathrm{Re}\xi$ (see (\ref{convexU}) or \cite{C1}) and is strictly concave up in $\mathrm{Im}\xi$ when sufficiently close to $\xi_{\boldsymbol\alpha}$ (see \cite[Proposition 3.11]{LMSWY2}), and by the compactness of the region $\big\{ (\boldsymbol \alpha, \xi)\in \mathbb R^6\times \mathbb C\ \big|\ \boldsymbol\alpha\in \overline{\mathrm H_{\epsilon_0}},|\xi-\xi_{\boldsymbol\alpha}|\leqslant \epsilon\big\},$  there exists a $c>0$ depending only on $\epsilon$ such that 
$$\sup_{s\in\mathbb R, |s|\leqslant \epsilon}\mathrm{Re}\bigg(\frac{\partial ^2V_0(\boldsymbol\alpha,\xi)}{\partial \xi^2}\Big|_{\xi=\xi_{\boldsymbol\alpha}+ts} \bigg)=\sup_{s\in\mathbb R, |s|\leqslant \epsilon}\frac{\partial ^2 \mathrm{Re}V_0(\boldsymbol\alpha,\xi)}{\partial \mathrm{Re}\xi^2}\Big|_{\xi=\xi_{\boldsymbol\alpha}+ts} \leqslant -2\pi c$$
and 
$$\sup_{s\in\mathbf i\mathbb R, |s|\leqslant \epsilon}\mathrm{Re}\bigg(\frac{\partial ^2 V_0(\boldsymbol\alpha,\xi)}{\partial \xi^2}\Big|_{\xi=\xi_{\boldsymbol\alpha}+ts} \bigg)=-\inf_{s\in\mathbf i\mathbb R, |s|\leqslant \epsilon}\frac{\partial ^2 \mathrm{Re}V_0(\boldsymbol\alpha,\xi)}{\partial \mathrm{Im}\xi^2}\Big|_{\xi=\xi_{\boldsymbol\alpha}+ts} \leqslant -2\pi c.$$
Together with (\ref{Taylor}), we have 
\begin{equation}\label{infc}
\inf_{s\in\mathbb R\cup\mathbf i\mathbb R,|s|\leqslant \epsilon} \mathrm{Re}g_{\boldsymbol\alpha}(s)\geqslant 2\pi c \int_0^1(1-t)dt=\pi c.    
\end{equation}

On the one hand, letting $s=u,$ we have
\begin{equation}\label{Ir'}
\mathcal I^r_{\boldsymbol \alpha}( \hbar)=\frac{1}{\sqrt{2\pi\hbar}}\int_{-\epsilon}^{\epsilon}    e^{-\frac{u^2g_{\boldsymbol\alpha}(u)}{\pi\hbar}}\exp\bigg(\frac{V_\hbar(\boldsymbol\alpha,\xi_{\boldsymbol\alpha}+u)-V_0(\boldsymbol\alpha,\xi_{\boldsymbol\alpha}+u)}{\pi \hbar }\bigg)du.
\end{equation}
Then  by (\ref{exp}), (\ref{supG}) and (\ref{infc}), all the conditions of Lemma \ref{Int} are satisfied with $x=\hbar,$ $\rho=2\epsilon,$ $\rho_0=\epsilon,$ $c=c,$ $G=G,$ $g(u)=\frac{g_{\boldsymbol\alpha}(u)}{\pi},$ $f(u,x)=\exp\big(\frac{V_\hbar(\boldsymbol\alpha,\xi_{\boldsymbol\alpha}+u)-V_0(\boldsymbol\alpha,\xi_{\boldsymbol\alpha}+u)}{\pi \hbar }\big)$ and $I(x)=\mathcal I^r_{\boldsymbol \alpha}( \hbar).$ Then by Lemma \ref{Int},  $\mathcal I^r_{\boldsymbol \alpha}( \hbar)$ satisfies a $(C,\sigma)$ Gevrey-$1$ estimate, where $C>0$ and  $\sigma>0$ depend only on $C_1,$ $\sigma_1,$  $c,$ $G$ and $\epsilon,$ which in turn depend only on $\epsilon,$
and
\begin{equation}\label{cn}
c_n(\boldsymbol\alpha)=\frac1{\sqrt2}\sum\limits_{\substack{j,k\geqslant 0\\j+k=n}}
\frac{F_{r,j}^{(2k)}(0)}{4^k k!},
\end{equation}
where for each $j\geqslant 0,$ $F_{r,j}$ is a holomorphic function defined in a neighborhood of $0$ such that if $v_r=v_r(u)\doteq u\sqrt{\frac{g_{\boldsymbol\alpha}(u)}{\pi}},$ then 
\begin{equation}\label{Frj}
    F_{r,j}(v_r)=\frac{b_j\big(\boldsymbol\alpha,\xi_{\boldsymbol\alpha}+u\big)}{\Big(u\sqrt{\frac{g_{\boldsymbol\alpha}(u)}{\pi}}\Big)'}=\frac{b_j\big(\boldsymbol\alpha,\xi_{\boldsymbol\alpha}+u\big)}{v_r'(u)}.
\end{equation}

On the other hand, letting $s=\mathbf iu,$ we have
\begin{equation}\label{Ib'}
\mathcal I^b_{\boldsymbol \alpha}( \hbar)=\frac{1}{\sqrt{2\pi\hbar}}\int_{-\epsilon}^{\epsilon}    e^{-\frac{u^2g_{\boldsymbol\alpha}(\mathbf iu)}{\pi\hbar}}\exp\bigg(\frac{V_{-\hbar}(\boldsymbol\alpha,\xi_{\boldsymbol\alpha}+\mathbf iu)-V_0(\boldsymbol\alpha,\xi_{\boldsymbol\alpha}+\mathbf iu)}{-\pi \hbar }\bigg)du.
\end{equation}
Then  by (\ref{exp}), (\ref{supG}) and (\ref{infc}), all the conditions of Lemma \ref{Int} are satisfied with $x=\hbar,$ $\rho=2\epsilon,$ $\rho_0=\epsilon,$ $c=c,$ $G=G,$ $g(u)=\frac{g_{\boldsymbol\alpha}(\mathbf i u)}{\pi},$ $f(u,x)=\exp\big(\frac{V_{-\hbar}(\boldsymbol\alpha,\xi_{\boldsymbol\alpha}+\mathbf iu)-V_0(\boldsymbol\alpha,\xi_{\boldsymbol\alpha}+\mathbf iu)}{-\pi \hbar }\big)$ and $I(x)=\mathcal I^b_{\boldsymbol \alpha}( \hbar).$ Then by Lemma \ref{Int},  $\mathcal I^b_{\boldsymbol \alpha}( \hbar)$ satisfies a $(C,\sigma)$ Gevrey-$1$ estimate for the same  $C$ and $\sigma$ as above, and
\begin{equation}\label{dn}
d_n(\boldsymbol\alpha)=\frac1{\sqrt2}\sum\limits_{\substack{j,k\geqslant 0\\j+k=n}}
\frac{F_{b,j}^{(2k)}(0)}{4^k k!},
\end{equation}
where for each $j\geqslant 0,$ $F_{b,j}$ is a holomorphic function defined in a neighborhood of $0$ such that if $v_b=v_b(u)\doteq u\sqrt{\frac{g_{\boldsymbol\alpha}(\mathbf iu)}{\pi}},$ then 
\begin{equation}\label{Fbj}
    F_{b,j}(v_b)=\frac{(-1)^jb_j\big(\boldsymbol\alpha,\xi_{\boldsymbol\alpha}+\mathbf iu\big)}{\Big(u\sqrt{\frac{g_{\boldsymbol\alpha}(\mathbf iu)}{\pi}}\Big)'}=\frac{(-1)^jb_j\big(\boldsymbol\alpha,\xi_{\boldsymbol\alpha}+\mathbf iu\big)}{v_b'(u)}.
\end{equation}

Since $\mathbf iv_b(u)=v_r(\mathbf iu)=\mathbf iu\sqrt{\frac{g_{\boldsymbol\alpha}(\mathbf iu)}{\pi}},$ we have $v_b'(u)=\big(v_r(\mathbf iu)\big)'.$ Together with (\ref{Frj}) and (\ref{Fbj}), we have
$$F_{b,j}\big(v_b(u)\big)=(-1)^j F_{r,j}\big(v_r(\mathbf iu)\big)=(-1)^j F_{r,j}\big(\mathbf iv_b(u)\big).$$
Therefore, for each $k\geqslant 0,$ we have
$$F^{(2k)}_{b,j}(0)=(-1)^{j+k} F^{(2k)}_{r,j}(0);$$ 
and together with (\ref{cn}) and (\ref{dn}), we have 
\begin{equation*}
    d_n(\boldsymbol\alpha)=(-1)^nc_n(\boldsymbol\alpha).
\end{equation*}
\end{proof}

\subsubsection{Proof of Proposition \ref{leading6j}}\label{624}

The proof of Proposition \ref{leading6j} relies on Proposition \ref{refI0} and Lemma \ref{I0I6j} below.

\begin{proposition}\label{refI0}
For $\boldsymbol\alpha=(\alpha_1,\dots,\alpha_6)\in \mathrm H_{\epsilon_0},$ let $\boldsymbol\alpha'\in \mathrm H_{\epsilon_0}$ be obtained from $\boldsymbol\alpha$ by changing $\alpha_k\mapsto 2\pi-\alpha_k$ for some $k\in\{1,\dots,6\}.$ Let $T$ be as in Proposition \ref{r6jb6j} (3). Then 
$$\mathcal L_T\mathcal B\widehat{\mathcal I_{\boldsymbol\alpha}^r}(\hbar)=\mathcal L_T\mathcal B\widehat{\mathcal I_{\boldsymbol\alpha'}^r}(\hbar).$$
\end{proposition}

To prove Proposition \ref{refI0}, we need the following 

\begin{lemma}\label{b6jIb}
Let $\hbar=b^2.$ Then there is a $\delta'>0$ such that for each $b$-admissible $\boldsymbol a=(a_1,\dots,a_6)$ with the corresponding $\boldsymbol\alpha\in\mathrm H_{\epsilon_0},$ the $b$-$6j$ symbol satisfies 
$$\Bigg|\sqrt{\frac{\pi}{2}
}e^{\frac{\mathrm{Vol}(\Delta_{\boldsymbol\alpha})}{\pi \hbar}}\bigg\{\begin{matrix} a_1 & a_2 & a_3 \\ a_4 & a_5 & a_6 \end{matrix} \bigg\}_b- \mathcal I^b_{\boldsymbol\alpha}(\hbar)\Bigg|<O\Big(e^{\frac{r}{2\pi} \big(-\delta'\big)}\Big).$$
\end{lemma}

\begin{proof} Let 
$U_b(\boldsymbol \alpha,\xi)$
be as defined (\ref{Ub}). Then by (\ref{sh}), we have
$$U_b(\boldsymbol\alpha,\xi)=-2\mathbf i V_{-\hbar}(\boldsymbol\alpha,\xi);$$
and by Proposition \ref{EST} (2), we have 
\begin{equation*}
\Bigg|\pi \mathbf i b\bigg\{\begin{matrix} a_1 & a_2 & a_3 \\ a_4 & a_5 & a_6 \end{matrix} \bigg\}_b- \frac{\mathbf i\sqrt{2\pi\hbar}}{e^{\frac{\mathrm{Vol}(\Delta_{\boldsymbol\alpha})}{\pi \hbar}}}\mathcal I^b_{\boldsymbol\alpha}(\hbar)\Bigg|\leqslant \int_{\Gamma\setminus [\xi_{\boldsymbol\alpha}-\mathbf i\epsilon,\xi_{\boldsymbol\alpha}+\mathbf i\epsilon]}\Big|\psi(\xi)e^{\frac{1}{2\pi \mathbf i b^2}U_b(\boldsymbol\alpha,\xi)}\Big|d\xi,
\end{equation*} 
where $\Gamma$ here is chosen to pass $\xi_{\boldsymbol\alpha}.$
Let $\epsilon_1$ be as in the proof of \cite[Theorem 1.5]{LMSWY2}, and let $\delta'=\epsilon_1.$
Then by (II) and (III) and letting $d=\epsilon$ in the proof of \cite[Theorem 1.5]{LMSWY2}, we have 
$$\int_{\Gamma\setminus [\xi_{\boldsymbol\alpha}-\mathbf i\epsilon,\xi_{\boldsymbol\alpha}+\mathbf i\epsilon]}\Big|\psi(\xi)e^{\frac{1}{2\pi \mathbf i b^2}U_b(\boldsymbol\alpha,\xi)}\Big|d\xi<O\Big(e^{\frac{r}{2\pi} \big(-\mathrm{Vol}(\Delta_{\boldsymbol\alpha})-\delta'\big)}\Big),$$
from which the result follows.
\end{proof}

\begin{proof}[Proof of Proposition \ref{b6jIb}]
By Lemma \ref{b6jIb}, $\frac{e^{\frac{\mathrm{Vol}(\Delta_{\boldsymbol\alpha})}{\pi \hbar}}}{\sqrt{2\pi\hbar}}\bigg\{\begin{matrix} a_1 & a_2 & a_3 \\ a_4 & a_5 & a_6 \end{matrix} \bigg\}_b$ and $\mathcal I^b_{\boldsymbol\alpha}(\hbar)$ satisfy the same Gevrey-$1$ estimate, hence share the same Borel transform. If we change $a_k\mapsto Q-a_k,$ then the corresponding $\alpha_k\mapsto 2\pi-\alpha_k.$ Then by Proposition \ref{b6jreflection} and the fact that $\Delta_{\boldsymbol\alpha}=\Delta_{\boldsymbol\alpha'},$ we have 
\begin{equation}\label{634}
d_n(\boldsymbol\alpha)=d_n(\boldsymbol\alpha')   
\end{equation}
for all $n\geqslant 0.$ As a consequence, 
$$\mathcal B\widehat{\mathcal I_{\boldsymbol\alpha}^b}(-\hbar)=\mathcal B\widehat{\mathcal I_{\boldsymbol\alpha'}^b}(-\hbar),$$ and by Proposition \ref{r6jb6j} (2),
$$\mathcal B\widehat{\mathcal I_{\boldsymbol\alpha}^r}(\hbar)=\mathcal B\widehat{\mathcal I_{\boldsymbol\alpha'}^r}(\hbar),$$ 
from which the result follows.
\end{proof}

\begin{remark}
By (\ref{634}) and Proposition \ref{r6jb6j} (1), we also get 
$$c_n(\boldsymbol\alpha)=c_n(\boldsymbol\alpha')$$
for all $n\geqslant 0,$ which gives an affirmative answer to \cite[Conjecture 3]{CM}.
\end{remark}

\begin{lemma}\label{I0I6j}
Let $\hbar=\frac{2}{r}$ and let $\delta_5>0$ be as in Proposition \ref{error6j}. Then for each $r$-hyperideal  $\boldsymbol a$ with the corresponding $\boldsymbol\alpha\in\mathrm H_{\epsilon_0},$ the $0$-th Fourier coefficient of the extended quantum $6j$-symbol satisfies
$$\Big|\mathrm I_{\boldsymbol a}(0)- \sqrt{2\pi\hbar}e^{\frac{\mathrm{Vol}(\Delta_{\boldsymbol\alpha})}{\pi \hbar}}\mathcal I^r_{\boldsymbol\alpha}(\hbar)\Big|<O\Big(e^{\frac{r}{2\pi} \big(\mathrm{Vol}(\Delta_{\boldsymbol\alpha})-\delta_5\big)}\Big).$$
\end{lemma}

\begin{proof}
By Proposition \ref{EST} (2), the relation (\ref{sh}) and the choice of $\delta_5$ in Proposition \ref{error6j}, we have 
\begin{equation*}
\begin{split}
\Big|\mathrm I_{\boldsymbol a}(0)- \sqrt{2\pi\hbar}e^{\frac{\mathrm{Vol}(\Delta_{\boldsymbol\alpha})}{\pi \hbar}}\mathcal I^r_{\boldsymbol\alpha}(\hbar)\Big|=& \bigg|\int_{\mathbb R\setminus [\xi_{\boldsymbol\alpha}-\epsilon,\xi_{\boldsymbol\alpha}+\epsilon]}\psi(\xi)e^{\frac{r}{4\pi \mathbf i}U_r(\boldsymbol\alpha,\xi)}d\xi\bigg|\\
\leqslant &\int_{[\xi_{\boldsymbol\alpha}-2\epsilon,\xi_{\boldsymbol\alpha}+2\epsilon]\setminus [\xi_{\boldsymbol\alpha}-\epsilon,\xi_{\boldsymbol\alpha}+\epsilon]}\Big|e^{\frac{r}{4\pi \mathbf i}U_r(\boldsymbol\alpha,\xi)}\Big|d\xi\\
<&O\Big(e^{\frac{r}{2\pi} \big(\mathrm{Vol}(\Delta_{\boldsymbol\alpha})-\delta_5\big)}\Big).
\end{split} 
\end{equation*}
\end{proof}

\begin{proof}[Proof of Proposition \ref{leading6j}] By Proposition \ref{refI0} we have
\begin{equation*}
\begin{split}
\mathrm I_{\boldsymbol a}(0)- \mathrm I_{\boldsymbol a'}(0) = &  \Big(\mathrm I_{\boldsymbol a}(0)- \sqrt{2\pi\hbar}e^{\frac{\mathrm{Vol}(\Delta_{\boldsymbol\alpha})}{\pi \hbar}}\mathcal I^r_{\boldsymbol\alpha}(\hbar)\Big)+ \sqrt{2\pi\hbar}e^{\frac{\mathrm{Vol}(\Delta_{\boldsymbol\alpha})}{\pi \hbar}}\Big(\mathcal I^r_{\boldsymbol\alpha}(\hbar)-\mathcal L_T\mathcal B\widehat{\mathcal I_{\boldsymbol\alpha}^r}(\hbar)\Big)\\
&-\sqrt{2\pi\hbar}e^{\frac{\mathrm{Vol}(\Delta_{\boldsymbol\alpha})}{\pi \hbar}}\Big(\mathcal I^r_{\boldsymbol\alpha'}(\hbar)-\mathcal L_T\mathcal B\widehat{\mathcal I_{\boldsymbol\alpha'}^r}(\hbar)\Big)-\Big(\mathrm I_{\boldsymbol a'}(0)- \sqrt{2\pi\hbar}e^{\frac{\mathrm{Vol}(\Delta_{\boldsymbol\alpha})}{\pi \hbar}}\mathcal I^r_{\boldsymbol\alpha'}(\hbar)\Big).
\end{split}
\end{equation*}
Let $\delta_5$ be as in Lemma \ref{I0I6j}, let $T$ be as in Proposition \ref{r6jb6j}, and let  
$$\delta_7< \min\{\delta_5, \pi T\}.$$
Then by Lemma \ref{I0I6j} and Proposition \ref{r6jb6j} (3), each of the four summands above has order less than $O\big(e^{\frac{r}{2\pi} (\mathrm{Vol}(\Delta_{\boldsymbol\alpha})-\delta_7)}\big),$ from which the result follows.
\end{proof}

\subsection{Proof of Theorem \ref{wrs} and Proposition \ref{key}}

\begin{proof}[Proof of Theorem \ref{wrs}]
Let $\delta_5,$ $\delta_6$ and $\delta_7$ respectively be as in Proposition \ref{error6j}, Proposition \ref{other6j} and Proposition \ref{leading6j}, and let 
$$\delta_4< \min\{\delta_5,\delta_6,\delta_7\}.$$
Then by (\ref{6.2}), Proposition \ref{error6j}, Proposition \ref{other6j} and Proposition \ref{leading6j}, we have 
\begin{equation*}
 \begin{split}
\Bigg|\bigg\|\begin{matrix} 
      a_1 & a_2 & a_3 \\
      a_4 & a_5 & a_6
 \end{matrix}\bigg\|-&\bigg\|\begin{matrix} 
      r-2-a_1 & a_2 & a_3 \\
      a_4 & a_5 & a_6
   \end{matrix}\bigg\|\Bigg| \\
   <& \sin\frac{2\pi}{r} \Big|\mathrm I_{\boldsymbol a}(0)- \mathrm I_{\boldsymbol a'}(0) \Big|+ O\Big(e^{\frac{r}{2\pi}\big(\mathrm{Vol}
(\Delta_{\boldsymbol\alpha})-\delta_4\big)}\Big)\\
=&O\Big(e^{\frac{r}{2\pi}\big(\mathrm{Vol}
(\Delta_{\boldsymbol\alpha})-\delta_4\big)}\Big).
\end{split}
\end{equation*}
\end{proof}

\begin{proof}[Proof of Proposition \ref{key}]
Let $\delta_3$ and $\delta_4$ respectively be as in Proposition \ref{key'} and Theorem \ref{wrs}, and let 
$$\delta_2=\min\{\delta_3,\delta_4\}.$$
Then the result follows from Proposition \ref{key'} and Theorem \ref{wrs} and the fact that the volume of any truncated hyperideal hyperbolic tetrahedron $\Delta_{\boldsymbol\alpha}$ is at most that of the regular ideal hyperbolic octahedron,  which equals $v_8.$ 
\end{proof}


\appendix
\section{Proof of the technical lemmas}\label{Gevrey-1}

\begin{proof}[Proof of Lemma \ref{GBL}]
To see (1), by (\ref{gevrey-1}) with $N=n+1,$ we have 
$$\bigg|\frac{f(x)-\sum_{k=0}^{n-1}a_kx^k}{x^n}-a_n\bigg|\leqslant C\sigma^{n+1}(n+1)!x;$$
and letting $x\to 0^+,$ we have
$$a_n=\lim_{x\to 0^+}\frac{f(x)-\sum_{k=0}^{n-1}a_kx^k}{x^n}.$$
By (\ref{gevrey-1}) again with $N=n,$ we have
$$\bigg|\frac{f(x)-\sum_{k=0}^{n-1}a_kx^k}{x^n}\bigg|\leqslant C\sigma^nn!;$$
and as a consequence, 
$$|a_n|\leqslant C\sigma^nn!$$
for each $n\geqslant 0.$
Therefore, the power series $\mathcal B\widehat f(x)$ converges for $x\in\mathbb C$ with $|x|<\sigma^{-1},$ defining a holomorphic function there. 
\\

To see (2), we use the identity that for any integer $n\geqslant 0,$
\begin{equation}\label{gamma}
x^n=\frac{1}{x}\int_0^\infty \frac{e^{-\frac{t}{x}}}{n!}t^ndt,    
\end{equation}
which could be obtained using the  Integration by Parts.
Then for any $T\in (0,\sigma^{-1})$ and $N\geqslant 0,$
 we have 
 \begin{equation*}
\begin{split}
\mathcal L_T\mathcal B\widehat f(x) -\sum_{n=0}^{N-1}a_nx^n = & \frac{1}{x}\int_0^Te^{-\frac{t}{x}}\mathcal B\widehat f(t)dt -  \frac{1}{x}\int_0^\infty e^{-\frac{t}{x}} \sum_{n=0}^{N-1} \frac{a_n}{n!}t^n dt\\
= & \frac{1}{x}\int_0^Te^{-\frac{t}{x}}\bigg(\mathcal B\widehat f(t)- \sum_{n=0}^{N-1} \frac{a_n}{n!}t^n\bigg) dt -  \frac{1}{x}\int_T^\infty e^{-\frac{t}{x}}\sum_{n=0}^{N-1} \frac{a_n}{n!}t^n dt.
\end{split}     
 \end{equation*}
To estimate the terms above, for $t<T$ we have 
\begin{equation}\label{t<T}
\bigg|\mathcal B\widehat f(t)- \sum_{n=0}^{N-1} \frac{a_n}{n!}t^n\bigg|\leqslant \sum_{n=N}^\infty \frac{|a_n|}{n!}|t|^n \leqslant \sum_{n=N}^\infty C\sigma^n|t|^n \leqslant \sum_{n=N}^\infty C(\sigma T)^n\bigg(\frac{t}{T}\bigg)^N  \leqslant  \frac{C}{1-\sigma T}\bigg(\frac{t}{T}\bigg)^N,    
\end{equation}
where the second inequality comes from (\ref{anbound}) and the third inequality comes from that $\frac{t}{T}<1;$
and for $t>T$ we have
\begin{equation}\label{t>T}\bigg|\sum_{n=0}^{N-1} \frac{a_n}{n!}t^n\bigg| \leqslant \sum_{n=0}^{N-1}C(\sigma T)^n\bigg(\frac{t}{T}\bigg)^n  \leqslant  \frac{C}{1-\sigma T}\bigg(\frac{t}{T}\bigg)^N,
\end{equation}
where the first inequality comes from (\ref{anbound}) and the second inequality comes from that $\frac{t}{T}>1.$
Putting (\ref{t<T}) and (\ref{t>T}) together, and by (\ref{gamma}), we have
$$\bigg| \mathcal L_T\mathcal B\widehat f(x) -\sum_{n=0}^{N-1}a_nx^n\bigg|\leqslant \frac{C}{1-\sigma T}\frac{1}{x}\int_0^\infty e^{-\frac{t}{x}}\bigg(\frac{t}{T}\bigg)^Ndt=\frac{C}{1-\sigma T} N!\bigg(\frac{x}{T}\bigg)^N;$$
and by (\ref{gevrey-1}), we have
\begin{equation}\label{f-LBf}
\begin{split}
  \bigg|f(x) - \mathcal L_T\mathcal B\widehat f(x)\bigg| \leqslant &   \bigg| f(x) -\sum_{n=0}^{N-1}a_nx^n\bigg| +  \bigg| \mathcal L_T\mathcal B\widehat f(x) -\sum_{n=0}^{N-1}a_nx^n\bigg|\\
  \leqslant & C\sigma^NN!x^N+\frac{C}{1-\sigma T} N!\bigg(\frac{x}{T}\bigg)^N  \leqslant  C\bigg(1+\frac{1}{1-\sigma T}\bigg)N! \bigg(\frac{x}{T}\bigg)^N
\end{split}    
\end{equation}
for all integers $N\geqslant 0.$
Finally, for $x\in (0,\min\{x_0,T\}),$ we let $N=\lfloor \frac{T}{x}\rfloor$ so that $\frac{T}{x}-1< N \leqslant  \frac{T}{x}.$ Then by the Stirling's Formula that
$$N! = \sqrt{2\pi N}\bigg(\frac{N}{e}\bigg)^Ne^{\frac{\theta}{12N}}$$ 
for some $\theta\in (0,1),$
we have 
\begin{equation}\label{stirling}
N!\bigg(\frac{x}{T}\bigg)^N = \sqrt{2\pi N} e^{N\log \big (N \cdot \frac {x}{T}\big )-N + \frac{\theta}{12N}}\leqslant \sqrt{2\pi} e^{\frac{13}{12}} \sqrt{\frac{T}{x}}e^{-\frac{T}{x}},    
\end{equation}
where the inequality comes from the choice of $N$ so that $\sqrt{N}\leqslant \sqrt{\frac{T}{x}},$ $N\cdot \frac{x}{T}\leqslant 1,$ $-N < 1 - \frac{T}{x}$ and $\frac{\theta}{12N}< \frac{1}{12}.$ 
Then (2) follows from (\ref{f-LBf}) and (\ref{stirling}) with 
$D=\sqrt{2\pi} e^{\frac{13}{12}}C\big(1+\frac{1}{1-\sigma T}\big).$
\end{proof}

\begin{proof}[Proof of Lemma \ref{Exp}]
For $n\geqslant 1$ and $k\geqslant 1,$ we let 
$$b_{n,k}=\sum_{\substack{
m_1,\ldots,m_n\geqslant0\\
m_1+\cdots+m_n=k\\
m_1+2m_2+\cdots+nm_n=n}}
\frac{k!}{m_1!\cdots m_n!}
\prod_{j=1}^{n}a_j^{m_j}$$
so that 
$b_n=e^{a_0}\sum_{k=0}^n \frac{b_{n,k}}{k!}.$
Then we have
\begin{equation}\label{induction}
b_{n,1}=a_n\quad\text{and}\quad b_{n,k+1}=\sum_{j=1}^{n-1}b_{j,k}a_{n-j}.    
\end{equation}

We will first prove that 
\begin{equation}\label{bnk}
 |b_{n,k}|\leqslant C^k\sigma^nn!   
\end{equation}
by an induction in $k.$ The base case $k=1$ follows directly from \ref{induction} and (\ref{gevrey-1}). Now suppose (\ref{bnk}) holds for $k$ and for all $n\geqslant 1.$ Then 
\begin{equation*}
\begin{split}
|b_{n,k+1}|\leqslant \sum_{j=1}^{n-1} |b_{j,k}||a_{n-j}|\leqslant C^{k+1} \sigma^n \sum_{j=1}^{n-1} j!(n-j)!\leqslant C^{k+1}\sigma^n n!,   
\end{split}    
\end{equation*}
where the last inequality comes from that 
$\sum_{j=1}^{n-1} j!(n-j)!\leqslant n!.$

Next we prove that 
\begin{equation}\label{fk}
\bigg|\big(f(x)-a_0)^k - \sum_{n=0}^{N-1}b_{n,k}x^n\bigg| \leqslant C^k\sigma^N N!x^N\sum_{j=0}^{k-1}(\sigma x)^j
\end{equation}
by an induction in $k.$ The base case $k=1$ follows directly from \ref{induction} and (\ref{gevrey-1}). 
Now suppose (\ref{fk}) holds for $k$ and for all $N\geqslant 1.$ Then by (\ref{induction}), we have 
$$\sum_{n=0}^{N-1}b_{n,k+1}x^n=\sum_{n=0}^{N-1}\sum_{j=1}^{n-1}b_{j,k}a_{n-j}x^n=\sum_{j=1}^{N-1}b_{j,k}x^j\sum_{l=1}^{N-j-1}a_lx^l;$$
and as a consequence, we have
\begin{equation*}
\begin{split}
&\big(f(x)-a_0\big)^{k+1} -\sum_{n=0}^{N-1}b_{n,k+1}x^n\\  
= & \big(f(x)-a_0\big)\bigg(\big(f(x)-a_0\big)^k-\sum_{n=0}^{N-1}b_{n,k}x^n\bigg) + \sum_{j=1}^{N-1}b_{j,k}x^j\bigg(f(x)-a_0-\sum_{l=1}^{N-j-1}a_lx^l\bigg).
\end{split}
\end{equation*}
By (\ref{gevrey-1}), we have
$$\bigg|f(x)-a_0-\sum_{l=1}^{N-j-1}a_lx^l\bigg|\leqslant C\sigma^{N-j}(N-j)!x^{N-j}$$
for each $j\leqslant N-1.$ Together with the induction hypothesis (\ref{fk}) for $k$ and (\ref{bnk}), we have 
\begin{equation*}
\begin{split}
&\bigg| \big(f(x)-a_0\big)^{k+1} -\sum_{n=0}^{N-1}b_{n,k+1}x^n\bigg|\\
\leqslant & C^{k+1}\sigma^N N!x^N \sum_{j=0}^{k-1} (\sigma x)^{j+1} + C^{k+1}\sigma^Nx^N \sum_{j=1}^{N-1} j!(N-j)!\\
\leqslant & C^{k+1}\sigma^N N!x^N \sum_{j=0}^{k} (\sigma x)^{j},
\end{split}    
\end{equation*}
where the last equality comes from that $\sum_{j=1}^{N-1} j!(N-j)!\leqslant N!.$

Finally, for $N=0,$ by (\ref{gevrey-1}), we have 
$$\big|e^{f(x)}\big| \leqslant e^{|a_0|}|e^{|f(x)-a_0|}\leqslant e^C e^{C\sigma x_0}.$$
For $N\geqslant 1,$ by (\ref{fk}), $|a_0|\leqslant C,$ $x<x_0,$ and the convention that $b_{n,k}=0$ for $k>n,$ we have 
\begin{equation*}
\begin{split}
\bigg|e^{f(x)}-\sum_{n=0}^{N-1}b_nx^n\bigg| = & \bigg|e^{a_0}\sum_{k=1}^\infty \frac{1}{k!}\bigg(\big(f(x)-a_0\big)^k-\sum_{n=0}^{N-1}b_{n,k}x^n\bigg)\bigg|\\
\leqslant& e^C\sigma^NN!x^N\sum_{k=1}^\infty \frac{C^k}{k!}\sum_{j=0}^{k-1}(\sigma x_0)^j,
\end{split}    
\end{equation*}
where $\sum_{k=1}^\infty \frac{C^k}{k!}\sum_{j=0}^{k-1}(\sigma x_0)^j=\frac{e^C-e^{C\sigma x_0}}{1-\sigma x_0}$ when $\sigma x_0\neq 1,$ and $\sum_{k=1}^\infty \frac{C^k}{k!}\sum_{j=0}^{k-1}(\sigma x_0)^j=Ce^C$ when $\sigma x_0 = 1,$ which converges in both cases. Then the result follows with 
$$C_1=e^C\max\bigg\{e^{C\sigma x_0}, \sum_{k=1}^\infty \frac{C^k}{k!}\sum_{j=0}^{k-1}(\sigma x_0)^j\bigg\},$$
which depends on $C,$ $\sigma$ and $x_0,$ and $\sigma_1=\sigma.$
\end{proof}

\begin{proof}[Proof of Lemma \ref{Int}]
We will construct intermediate terms $A_n(x),$ $B_n(x)$ and $C_n(x)$ that connect $I(x)$ with its desired Gevrey - $1$ expansion $\sum_{n=0}^\infty b_nx^n,$
and estimate the difference between the adjacent two of them.
\medskip

As a preparation, we will define an injective function $V$ of $u$ such  that $V(u)^2=u^2g(u),$ and describe its inverse function $U$ and estimate its derivative. To do these, we need the following estimates. For any $u\in\mathbb C$ with $|u|<\rho$ and  $s\in(0,\rho-|u|),$ the circle $|z-u|=s$ lies in the disk $|u|<\rho;$ and by Cauchy's Derivative Formula, we have
$$                
 \big|g'(u)\big|
 \leqslant \frac1{2\pi}\int_{|z-u|=s}
                   \frac{|g(z)|}{|z-u|^2}\,|dz|
 \leqslant \frac{G}{s}.
$$ 
Letting $s\to \rho-|u|,$ we get 
$$ \big|g'(u)\big|\leqslant \frac{G}{\rho-|u|}.$$
Let $d=\frac{\rho c}{8G}\leqslant \frac{\rho}{8}.$ If $|u|\leqslant d,$ then 
 $$|g'(u)|
 \leqslant \frac{G}{\rho-d}\leqslant \frac{8G}{7\rho};$$
and integrating along the straight segment from $0$ to $u,$ we have
$$
\big|g(u)-g(0)\big|
 \leqslant \frac{8Gd}{7\rho}=\frac{c}{7}.$$
As a consequence, 
$$\mathrm{Re} g(u)\geqslant \frac{6c}{7}\quad \text{and}\quad \big|ug'(u)\big|\leqslant\frac c7$$ 
on the disk $|u|< d.$
 In particular, $\mathrm{Re}g(u)>0,$ and $\sqrt{g(u)}$ is well defined, holomorphic, and with the argument in $\{-\frac{\pi}{4},\frac{\pi}{4}\}$ on the dist $|u|<d.$ This implies that $\mathrm{Re}\big(\sqrt{g(u)}\overline{\sqrt{g(0)}}\big)>0;$ and as a consequence, we have
$$\Big|\sqrt{g(u)}+\sqrt{g(0)}\Big|^2
 =\big|g(u)\big|+\big|g(0)\big|
          +2\mathrm{Re}\Big(\sqrt{g(u)}\,\overline{\sqrt{g(0)}}\Big)\geqslant\big |g(0)\big|,$$
$$ \Big|\sqrt{g(u)}-\sqrt{g(0)}\Big|
 =\frac{\big|g(u)-g(0)\big|}{\big|\sqrt{g(u)}+\sqrt{g(0)}\big|}
 \leqslant \frac{c}{7\big|\sqrt{g(0)}\big|}\leqslant \frac{1}{7}\big|\sqrt{g(0)}\big|,$$
 and
$$\Big|\sqrt{g(u)}\Big|\geqslant \Big|\sqrt{g(0)}\Big|-\Big|\sqrt{g(u)}-\sqrt{g(0)}\Big|
 \geqslant \frac{6}{7}\big|\sqrt{g(0)}\big|.$$
 
 Now define 
 $$V(u)\doteq u\sqrt{g(u)}$$ 
 so that  $u^2g(u)=V(u)^2.$
 
 We first prove that $V$ is injective on the disk $|u|<d.$ The estimates above imply
\begin{equation}\label{eq:v-derivative}
\Big|V'(u)-\sqrt{g(0)}\Big|
 \leqslant  \Big|\sqrt{g(u)}-\sqrt{g(0)}\Big|+\frac{|ug'(u)|}{2|\sqrt{g(u)}|} \leqslant  \frac{1}{7}\Big|\sqrt{g(0)}\Big|+\frac{1}{12}\frac{c}{\big|\sqrt{g(0)}\big|}
<\frac{1}{4}\Big|\sqrt{g(0)}\Big|.    
\end{equation}
Therefore, 
for any two points $u_1\neq u_2$ in the disk $|u|<d,$ integrating along the straight line segment
from $u_1$ to $u_2$ and using  \eqref{eq:v-derivative}, we have
\begin{equation}\label{diff}
\begin{split}
\Big|V(u_1)-V(u_2)-\sqrt{g(0)}(u_1-u_2)\Big|
\leqslant \int_{u_1}^{u_2}\Big|V'(u)-\sqrt{g(0)}\Big||du|
 \leqslant \frac{1}{4}\big|\sqrt{g(0)}\big||u_1-u_2|;    
\end{split}    
\end{equation}
and as a consequence,
\begin{equation}\label{v-v}
\big|V(u_1)-V(u_2)\big|\geqslant \frac{3}{4}\Big|\sqrt{g(0)}\Big||u_1-u_2|,   
\end{equation}
implying the injectivity of $V.$

Next, we describe the inverse function $U$ of $V,$ and estimate its derivative.  Let $l=\frac{\sqrt c\,d}{2}.$
For any $v_0$  with
$|v_0|<l$ and every $|u|=d,$ we have
\begin{equation*}
\Big|V(u)-v_0-\sqrt{g(0)} u\Big|\leqslant \Big|V(u)-\sqrt{g(0)} u\Big|+|v_0|\leqslant \frac{1}{4}\Big|\sqrt{g(0)}\Big|d+l <\Big|\sqrt{g(0)}\Big||u|,
\end{equation*}
where the second inequality comes from (\ref{diff}) with $u_1=u$ and $u_2=0,$ and the last inequality comes from $c<|g(0)|.$ Then by Rouch\'e's Theorem, for each $v_0$ with $|v_0|<l,$ there is exactly  one zero  $u_0$ of
$V(u)-v_0$ in the disk $|u|<d.$
Therefore, the map $U: \big\{v\in \mathbb C\ \big|\ |v|<l\big\} \to \big\{ u\in\mathbb C\ \big|\ |u|<d\big\}$ sending $v_0$ to $u_0$ is
well-defined, holomorphic, and  satisfies $U(V(u))=u.$ By (\ref{v-v}), we have
$$\big|V'(u)\big|\geqslant \frac{3}{4}\Big|\sqrt{g(0)}\Big|\geqslant \frac{3\sqrt{c}}{4}$$
for each $u$ with $|u|<d,$
and as a consequence, we have
\begin{equation}\label{U'est}
\big|U'(v)\big|\leqslant \frac{4}{3\sqrt{c}} < \frac{2}{\sqrt{c}}   \end{equation}
for each $v$ with $|v|<l.$ We have the following steps I - IV.
\medskip

\noindent Step I. For  $n\geqslant 0,$ we define
 $$A_n(x)=\frac1{\sqrt{2\pi x}}\int_{-\rho_0}^{\rho_0} e^{-\frac{u^2g(u)}{x}}a_n(u)du.$$
 Then for each $N\geqslant 1,$ we have
\begin{equation}\label{IA}
 \begin{split}
\bigg|I(x)-\sum_{n=0}^{N-1}A_n(x)x^n\bigg|
 &=\Bigg|\frac1{\sqrt{2\pi x}}\int_{-\rho_0}^{\rho_0} e^{-\frac{u^2g(u)}{x}}
       \bigg(f(u,x)-\sum_{n=0}^{N-1}a_n(u)x^n\bigg)\,du\Bigg|\\
 &\leqslant \frac{C\sigma^N N!x^N}{\sqrt{2\pi x}}\int_{\mathbb R} e^{-\frac{cu^2}{x}}\,du
 \leqslant  \frac{C}{\sqrt{2c}}\,\sigma^N N!x^N.
 \end{split}
\end{equation}
\medskip

\noindent  Step II. Let $\mu=\frac12\min\big\{\rho_0,d,\frac{\ell}{4\sqrt G}\big\},$
and for $n\geqslant 0,$ define 
 $$B_n(x)=\frac1{\sqrt{2\pi x}}\int_{-\mu}^{\mu} e^{-\frac{u^2g(u)}{x}}a_n(u)du.$$
By (\ref{anbound}), we have 
\begin{equation}\label{supan}
\sup_{|u|<\rho}|a_n(u)|\leqslant  C\sigma^n n!.    
\end{equation}
Together with that $e^{-\frac{cu^2}{x}}\leqslant 
 e^{-\frac{c\mu^2}{2x}}e^{-\frac{cu^2}{2x}}$ for all $u\in [-\rho_0,\rho_0]\setminus [-\mu,\mu],$ we have
 \begin{equation}\label{AB}
\begin{split}
    \big|A_n(x)-B_n(x)\big|
 =&\Bigg|\frac1{\sqrt{2\pi x}}\int_{[-\rho_0,\rho_0]\setminus[-\mu,\mu]}
                     e^{-\frac{u^2g(u)}{x}}a_n(u)\,du\Bigg|\\
 \leqslant& \frac{C}{\sqrt{2\pi x}}\sigma^n n!e^{-\frac{c\mu^2}{2x}}
                     \int_{\R}e^{-\frac{cu^2}{2x}}\,du\\
 =&\frac{C}{\sqrt c}\sigma^n n!e^{-\frac{c\mu^2}{2x}} < \frac{C}{\sqrt c}\sigma^n n!(N-n)! \bigg(\frac{2x}{c\mu^2}\bigg)^{N-n},
\end{split}
\end{equation}
where the last inequality comes from that 
\begin{equation} \label{e>}
    e^v>\frac{v^{N-n}}{(N-n)!}
\end{equation}
for all $v\geqslant 0.$
\medskip

\noindent Step III. For $n\geqslant 0,$ define the holomorphic function $F_n$ on the disk $|v|<l$ by 
$$F_n(v)\doteq a_n\big(U(v)\big)U'(v)=\frac{a_n(u)}{V'(u)}=\frac{a_n(u)}{\big(u\sqrt{g(u)}\big)'}$$
for $u$ satisfying $v=V(u)=u\sqrt{g(u)}.$
Then by (\ref{supan}) and (\ref{U'est}), we have
\begin{equation}\label{eq:F-bound}
 \big|F_n(v)\big|\leqslant \frac{2C}{\sqrt c}\sigma^n n!.
\end{equation}
Let $\gamma=\sqrt{c}\mu;$ and for $n\geqslant 0,$ define 
$$C_n(x)=\frac1{\sqrt{2\pi x}}
          \int_{-\gamma}^{\gamma}e^{-\frac{v^2}{x}}F_n(v)dv.$$
We will show that 
\begin{equation}\label{BC}
 \big|B_n(x)-C_n(x)\big|\leqslant \frac{4C G}{c\sqrt{2\pi e c}}\,
    \sigma^nn!(N-n)!\bigg(\frac{2G x}{c\gamma^2}\bigg)^{N-n}.
\end{equation}
To this end, let $\Gamma=V([-\mu,\mu]),$ 
and by the substitution  $u=U(v),$ we have
$$
 B_n(x)=\frac1{\sqrt{2\pi x}}\int_\Gamma e^{-\frac{v^2}{x}}F_n(v)dv.$$
For $u\in [-\mu,\mu],$ by the definition of $\mu$ and the assumption of the Lemma, we have
\begin{equation}\label{eq:Gamma-containment}
\big|V(u)\big|=|u|\Big|\sqrt{g(u)}\Big|
 \leqslant \mu\sqrt G \quad \text{and}\quad \big|V(\pm\mu)\big|\geqslant \mu\sqrt{c}=\gamma.
\end{equation}
Let $\theta=\frac12\arccos\left(\frac{c}{G}\right)<\frac{\pi}{4}.$
For $u\in[-\mu,\mu],$ since
$\cos\big(\arg g(u)\big)=\frac{\mathrm{Re} g(u)}{|g(u)|}\geqslant \frac cG,$ we have
$\big|\arg\sqrt{g(u)}\big|\leqslant \theta.$
Let $\Gamma_+$ and $\Gamma_-$ respectively be as shown in Figure \ref{Gammapm}.
\begin{figure}[htbp]
\centering
\includegraphics[scale=0.6]{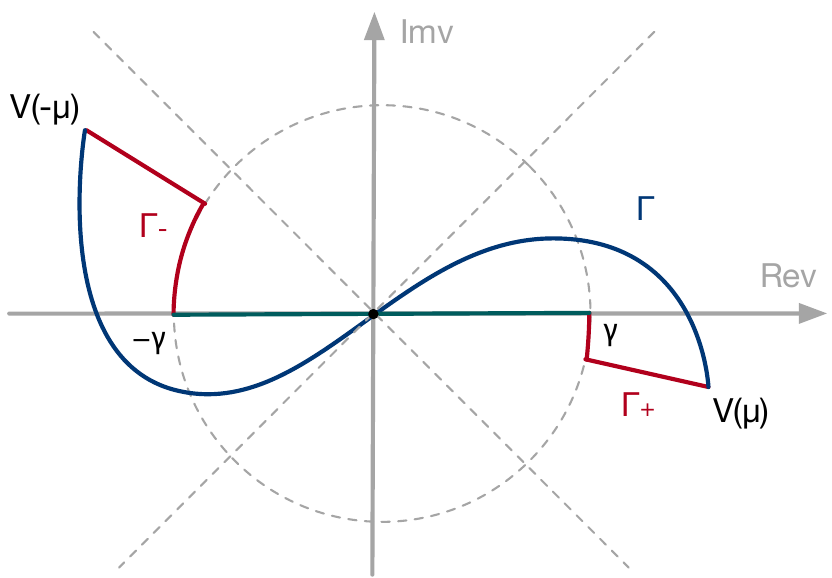}
\caption{The curves $\Gamma_+$ and $\Gamma_-.$}
\label{Gammapm} 
\end{figure}

Then we have the following:
\begin{enumerate}[(1)]
    \item $\Gamma$ is homotopic to $\Gamma_-\cup [-\gamma,\gamma]\cup\Gamma_+$ relative to the end points $V(\pm \mu).$

    \item The lengths 
   \begin{equation*}
      |\Gamma_+|+|\Gamma_-|\leqslant \big|V(\mu)\big|+\big|V(-\mu)\big|-2\gamma(1-\theta)\leqslant 2\mu\sqrt{G}. 
   \end{equation*}

    \item Each $v\in \Gamma_{\pm}$ satisfies 
$|v|\geqslant \gamma$ and either  $|\arg v|\leqslant  \theta$ or $|\arg(-v)|\leqslant  \theta;$ and as a consequence, we have 
\begin{equation*}
 \mathrm{Re}\big(v^2\big)\geqslant \gamma ^2\cos(2\theta)\geqslant \frac{c\gamma^2}{G}  
\end{equation*}
for each $v\in L.$
\end{enumerate}
Therefore, by (1), we have 
$$B_n(x)-C_n(x)=\frac1{\sqrt{2\pi x}}\int_{\Gamma_+\cup\Gamma_-} e^{-\frac{v^2}{x}}F_n(v)dv;$$
and by (2), (3), (\ref{eq:F-bound}) and $\gamma=\mu\sqrt{c}$,  we have
\begin{equation}
\label{B-C}
\begin{split}
 \big|B_n(x)-C_n(x)\big| \leqslant     \frac1{\sqrt{2\pi x}}\int_{\Gamma_+\cup\Gamma_-} e^{-\frac{\mathrm{Re}v^2}{x}}\big|F_n(v)\big||dv| \leqslant  \frac{4C\gamma\sqrt G}{c\sqrt{2\pi}}
             \sigma^n n! \frac{1}{\sqrt{x}}e^{-\frac{c\gamma^2}{Gx}}.
\end{split}    
\end{equation}
As $e^{\frac{c\gamma^2}{G x}-1}\geqslant \frac{c\gamma^2}{G x},$ we have $\frac{1}{\sqrt{x}}e^{-\frac{c\gamma^2}{2Gx}} \leqslant \frac{\sqrt{G}}{\gamma\sqrt{e c}}$, and as a consequence, 
\begin{equation}
    \label{B-C2}
    \frac{1}{\sqrt{x}}e^{-\frac{c\gamma^2}{Gx}}=\frac{1}{\sqrt{x}}e^{-\frac{c\gamma^2}{2Gx}} e^{-\frac{c\gamma^2}{2Gx}}< \frac{\sqrt{G}}{\gamma\sqrt{e c}} (N-n)!\bigg(\frac{2G x}{c\gamma^2}\bigg)^{N-n},
\end{equation}
where the inequality comes from (\ref{e>}). Then (\ref{BC}) follows from 
(\ref{B-C}) and (\ref{B-C2}). 
\medskip

\noindent Step IV. We will show that for each $n\geqslant 0,$
\begin{equation}\label{eq:Q-remainder}
 \Bigg|C_n(x)-\frac1{\sqrt2}\sum_{k=0}^{N-n-1}
       \frac{F_n^{(2k)}(0)}{4^k k!}x^k\Bigg|
 \leqslant \frac{20C}{3\sqrt{2c}}\sigma^n \gamma^{-2(N-n)}n!(N-n)!x^{N-n}.
\end{equation}
To see this, as $2\gamma<l,$ the function $F_n(v)$ is holomorphic in the disk $|v|<2\gamma,$ and we have
$$F_n(v)=\sum_{i=0}^\infty c_{n,i}v^i$$
with 
$$c_{n,i}= \frac{F_n^{(i)}(0)}{i!};$$
and by the Cauchy Integral Theorem and (\ref{eq:F-bound}), we have 
\begin{equation}\label{cni}
 |c_{n,i}|\leqslant \frac1{2\pi}\int_{|v|=2\gamma}
                 \frac{|F_n(v)|}{|v|^{i+1}}|dv|
 \leqslant \frac1{2\pi}\frac{4\pi \gamma}{(2\gamma)^{i+1}}\frac{2C}{\sqrt c}\sigma^n n!
 =\frac{1}{(2\gamma)^i}\frac{2C}{\sqrt c}\sigma^n n!.
\end{equation}
By the Gaussian Moment Formula that for every integer $k\geqslant 0$,  
\begin{equation}\label{Gauss}
 \frac1{\sqrt{2\pi x}}\int_{\R}e^{-\frac{v^2}{x}}v^{2k}\,dv
 =\frac1{4^k \sqrt2  }\binom{2k}{k}k!x^k;
\end{equation}
and by the fact that the integral of odd functions vanishes on $[-\gamma,\gamma],$ for any $m\geqslant 0,$ we can write 
\begin{equation}\label{GEE}
\begin{split}
&C_n(x)-\frac1{\sqrt2}\sum_{k=0}^{m-1}
       \frac{F_n^{(2k)}(0)}{4^k k!}x^k\\  
=& \frac1{\sqrt{2\pi x}}\int_{-\gamma}^{\gamma}e^{-\frac{v^2}{x}}
       \Big(F_n(v)-\sum_{i=0}^{2m-1}c_{n,i}v^i\Big)dv- 
\frac1{\sqrt{2\pi x}}\sum_{j=0}^{m-1}c_{n,2j}
       \int_{\mathbb R\setminus[-\gamma,\gamma]}e^{-\frac{v^2}{x}}v^{2{j}}dv.       
\end{split}
\end{equation}

On the one hand, for $v$ with $|v|\leqslant \gamma$, by
(\ref{cni}), we have
$$\bigg|F_n(v)-\sum_{i=0}^{2m-1}c_{n,i}v^i\bigg|
 \leqslant  \frac{2C}{\sqrt c}\sigma^n n!\sum_{i=2m}^{\infty}\left(\frac{|v|}{2\gamma}\right)^i
 \leqslant\frac{4C}{\sqrt c}\sigma^n n!\left(\frac{|v|}{2\gamma}\right)^{2m}.$$
As a consequence, we obtain
\begin{equation}\label{E1}
 \begin{aligned}
\bigg|\frac1{\sqrt{2\pi x}}\int_{-\gamma}^{\gamma}e^{-\frac{v^2}{x}}
       \Big(F_n(v)-\sum_{i=0}^{2m-1}c_{n,i}v^i\Big)dv\bigg|
 &\leqslant \frac{(2\gamma)^{-2m}}{\sqrt{2\pi x}} \frac{4C}{\sqrt c}\sigma^n n!
       \int_{\R}e^{-\frac{v^2}{x}}v^{2m}\,dv\\
 &\leqslant \frac{4^{-m+1}C}{\sqrt {2c}}\sigma^n n!\gamma^{-2m}m!x^m,
 \end{aligned}
\end{equation}
where the last inequality comes from (\ref{Gauss}) for $k=m$ and $\binom{2m}{m}\leqslant 4^{m}.$

On the other hand, for  $v$ with $v> \gamma$ and for $j< m,$ we have
$$|v|^{2j}=|v|^{2m}|v|^{-2(m-j)}
 \leqslant  |v|^{2m}\gamma^{-2(m-j)};$$
 and by (\ref{cni}), we have
 $$|c_{n,2j}|\,|v|^{2j}
 \leqslant  \frac{2C}{\sqrt c}\sigma^n n!\frac{1}{(2\gamma)^{2j}}\gamma^{-2(m-j)}|v|^{2m}
 =\frac{2C}{\sqrt c}\sigma^n n!4^{-j}\gamma^{-2m}|v|^{2m}.$$
As a consequence, 
\begin{equation}\label{E2}
 \begin{aligned}
\bigg|\frac1{\sqrt{2\pi x}}\sum_{j=0}^{m-1}c_{n,2j}
       \int_{\mathbb R\setminus[-\gamma,\gamma]}e^{-\frac{v^2}{x}}v^{2{j}}dv\bigg|
 &\leqslant \frac1{\sqrt{2\pi x}}\sum_{j=0}^{m-1}|c_{n,2j}|
       \int_{\mathbb R\setminus[-\gamma,\gamma]}e^{-\frac{v^2}{x}}|v|^{2j}dv\\
 &\leqslant \frac{\gamma^{-2m}}{\sqrt{2\pi x}}\frac{2C}{\sqrt c}\sigma^n n!
       \bigg(\sum_{j=0}^{m-1}4^{-j}\bigg)
       \int_{\mathbb R\setminus[-\gamma,\gamma]}e^{-\frac{v^2}{x}}|v|^{2m}dv\\
 &\leqslant \frac{\gamma^{-2m}}{\sqrt{2\pi x}}
      \frac{2C}{\sqrt c}\sigma^n n! \bigg(\sum_{j=0}^{m-1}4^{-j}\bigg)
       \int_{\mathbb R}e^{-\frac{v^2}{x}}v^{2m}dv\\
 &\leqslant \frac1{\sqrt2}\frac{2C}{\sqrt c}\sigma^n n!\bigg(\sum_{j=0}^{m-1}4^{-j}\bigg)\gamma^{-2m}m!x^m\\
 &\leqslant \frac{8C}{3\sqrt {2c}}\sigma^n n!\gamma^{-2m}m!x^m,
 \end{aligned}
\end{equation}
where the penultimate inequality comes from (\ref{Gauss}) again. 

Then (\ref{eq:Q-remainder}) follows from (\ref{GEE}) with $m=N-n,$ (\ref{E1}), (\ref{E2}) and the computation that 
$4^{-m+1}+\frac{8}{3}\leqslant \frac{20}{3}$ for $m\geqslant 0.$ 
\medskip

Finally, for $n\geqslant 0,$ define
\begin{equation*}
 c_n=\frac1{\sqrt2}\sum\limits_{\substack{j,k\geqslant 0\\j+k=n}}
                   \frac{F_j^{(2k)}(0)}{4^k k!}.
\end{equation*}
Then $$
 \sum\limits_{n=0}^{N-1}c_nx^n
 =\frac1{\sqrt2}\sum\limits_{j=0}^{N-1}x^j
       \sum\limits_{k=0}^{N-j-1}\frac{F_j^{(2k)}(0)}{4^k k!}x^k.$$
By(\ref{IA}), (\ref{AB}), (\ref{BC}), (\ref{eq:Q-remainder}) and the Binomial Theorem, we have
\begin{equation*}
 \bigg|I(x)-\sum\limits_{n=0}^{N-1}c_nx^n\bigg|\leqslant C\sqrt{\frac{2}{c}}\Bigg(\frac{1}{2}\sigma^N + \frac{1}{\sqrt{2}}\bigg(\sigma+\frac{2}{c\mu^2}\bigg)^N+ \frac{2G}{c\sqrt{\pi e}}\bigg(\sigma+\frac{2G}{c\gamma^2}\bigg)^N+ \frac{10}{3}\bigg(\sigma+\frac{1}{\gamma^2}\bigg)^N\Bigg)N!x^N.
\end{equation*}
As a consequence,
 $$\bigg|I(x)-\sum\limits_{n=0}^{N-1}c_nx^n\bigg|\leqslant C_2\sigma_2^N N!x^N$$
 with 
$$C_2=C\sqrt{\frac{2}{c}}\bigg(\frac{1}{2}+\frac{1}{\sqrt 2}+ \frac{2G}{c\sqrt{\pi e}}+\frac{10}{3}\bigg)\quad\text{and}\quad \sigma_2=\max\bigg\{\sigma, \sigma+\frac{2}{c\mu^2},\sigma+\frac{2G}{c\gamma^2}, \sigma+\frac{1}{\gamma^2}\bigg\},$$
 which is the desired real Gevrey- $1$ estimate of $I(x).$
 \end{proof}

\end{document}